\documentclass[12pt,twoside]{article}

\usepackage{amsmath}
\usepackage{amssymb}
\usepackage{mathrsfs}
\usepackage{amsthm}
\usepackage{needspace}
\usepackage{cite}

\usepackage{latexsym}

\usepackage{indentfirst}
\usepackage{color}
\usepackage{txfonts}

\usepackage{anysize}

\usepackage[colorlinks=true,
  linkcolor=red,
  citecolor=blue,
  urlcolor=magenta]{hyperref}

\def\R{{\mathbb R}}

\def\D{\mathcal{D}}
\def\M{\mathcal{M}}

\newtheorem{theorem}{Theorem}[section]
\newtheorem{lemma}[theorem]{Lemma}
\newtheorem{corollary}[theorem]{Corollary}
\newtheorem{proposition}[theorem]{Proposition}

\theoremstyle{definition}
\newtheorem{remark}[theorem]{Remark}
\newtheorem{definition}[theorem]{Definition}
\newtheorem{assumption}[theorem]{Assumption}

\numberwithin{equation}{section}

\newcommand{\C}{\mathbb C}
\newcommand{\Z}{\mathbb Z}
\newcommand{\cS}{\mathcal S}
\newcommand{\cP}{\mathcal P}
\newcommand{\supp}{\operatorname{supp}}
\newcommand{\esssup}{\operatorname*{ess\,sup}}

\begin{document}
\title{\bf\Large Matrix-weighted $Z$-spaces and tent spaces:\newline
characterizations,
duality and interpolation
\footnotetext{\hspace{-0.35cm} 2020 {\it
Mathematics Subject Classification}. Primary 42B35; Secondary 42B25, 46B70, 46E30. \endgraf
{\it Key words and phrases}. matrix weight, $Z$-space, tent space, Whitney average, duality, interpolation.
\endgraf
The work is supported by the National Natural Science Foundation of China (Grant No.
12561002), the Science and Technology Project of Guangxi (Guike AD25069086) and Hainan Provincial Natural Science Foundation of China (Grant No. 126MS0135).
}}
\author{Tengfei Bai, Pengfei Guo and Jingshi Xu\footnote{Corresponding author,
E-mail: \texttt{jingshixu@126.com}}}
\date{}
\maketitle

\vspace{-0.8cm}

\begin{center}
\begin{minipage}{13cm}
{\small {\bf Abstract:}\quad
Let $\alpha\in\R$, $0<p<\infty$, $0<q,r\le\infty$, and
$W$ be a matrix $\mathcal A_p$ weight on $\R^n$.
We introduce the matrix-weighted $Z$-spaces $Z^{p,q,r}_\alpha(W)$
and tent spaces $T^{p,q,r}_\alpha(W)$ on the upper half-space
and establish their characterizations using reducing operators.
We then characterize homogeneous matrix-weighted Besov and
Triebel--Lizorkin spaces in terms of the corresponding $Z$-spaces
and tent spaces, respectively, using Schwartz kernels satisfying
Tauberian and finite moment conditions.
Furthermore, we establish duality for the matrix-weighted
$Z$-spaces and tent spaces and obtain real and complex
interpolation results for these spaces and the associated
homogeneous matrix-weighted Besov and Triebel--Lizorkin spaces.
}
\end{minipage}
\end{center}

\setcounter{tocdepth}{2}
\tableofcontents

\section{Introduction}\label{sec:introduction}

Besov and Triebel--Lizorkin spaces play a fundamental role in
the theory of function spaces, include many classical function
spaces, and have important applications to partial differential
equations and other areas of analysis.
Their classical theory originates in the work of Besov
\cite{Besov59,Besov61}, Lizorkin \cite{Lizorkin72,Lizorkin74},
Triebel \cite{Triebel73} and Peetre \cite{Peetre73,Peetre75}.
The decompositions and $\varphi$-transform of Frazier and Jawerth
\cite{FJ85,FJ90} provide sequence-space descriptions of these
spaces. For the classical theory, we refer to \cite{Sawano18,Triebel83}.

The matrix-weighted theory builds on the work of Treil and
Volberg \cite{TV97} and Goldberg \cite{Goldberg03} on matrix
$\mathcal A_p$ weights. Roudenko
\cite{Roudenko03,RoudenkoDuality} introduced matrix-weighted
Besov spaces and studied their duality. Frazier and Roudenko
\cite{FRBesov,FR} established characterizations of matrix-weighted
Besov and Triebel--Lizorkin spaces using reducing operators, allowing
integrability indices below one.

Recent works of Yang and his collaborators concern
matrix-weighted Besov-type and Triebel--Lizorkin-type spaces,
real-variable characterizations, logarithmic smoothness, and
extensions to spaces of homogeneous type and to settings with
multiple parameters or variable exponents; see
\cite{MultiMatrix,BHYY,BHYYII,BHYYIII,BHYYInfinity,BYYHomogeneous,BYYZ,LYYLog,MatrixPoincare,WYZ,VariableMatrixBesov,OptimalMatrixOperators}.
Related studies cover Hardy and Campanato spaces, variable and
local matrix weights, and singular and fractional integrals
\cite{MatrixHardy,MatrixInfinity,MatrixCampanato,VariableMatrixHardy,LocalMatrix,VariableMatrixWeights,MatrixSharpCZ,MatrixFractional}.
For matrix-weighted anisotropic smoothness spaces, see
\cite{BGX26,BFXYYAnisotropic}. For function spaces with
operator-valued $\mathcal A_p$ weights in Banach spaces, see
\cite{OperatorWeightedSpaces}. An overview is provided in
\cite{MatrixSurvey}.

Upper-half-space extensions provide a continuous approach to
the study of these smoothness spaces.
The theory of tent spaces was initiated by Coifman, Meyer and
Stein \cite{CMS85}, who placed non-tangential maximal functions
and square functions in a common framework closely related to
$L^p$ and Hardy spaces. In the study of elliptic operators with
rough coefficients, Hofmann, Mayboroda and McIntosh
\cite[Section 8.3]{HMM11} developed tent spaces with power
weights in the vertical variable and proved complex interpolation
results for these spaces. Huang \cite{HuangTent} introduced a further scale
based on Whitney averages, with a local averaging exponent
independent of the conical integrability exponent. His strong
factorization theorems led to complex interpolation in the
quasi-Banach range and to multiplier and duality results.

The $Z$-spaces have a related origin in elliptic boundary value
problems. Barton and Mayboroda \cite{BM16} introduced spaces
defined by local averages in Whitney regions to treat fractional
Besov boundary data and establish mapping properties of layer
potentials. They also proved real and complex interpolation
results for these averaged spaces. Amenta \cite{Amenta18}
studied this scale under the name of $Z$-spaces and identified
certain real interpolants of weighted tent spaces with
$Z$-spaces.

The present paper is mainly inspired by Auscher, Bechtel and
Haardt \cite{ABH}. They introduced a three-index scale of
$Z$-spaces, allowing the spatial integrability and scale
summability indices to vary independently, and characterized
homogeneous Besov spaces by Whitney averages. They obtained
the corresponding characterizations of homogeneous
Triebel--Lizorkin spaces using tent spaces with an independent
local averaging exponent. Their work also treats completeness,
duality, embeddings, and real and complex interpolation,
including the identification of real interpolants of tent
spaces with $Z$-spaces. For related continuous characterizations,
see \cite{Ullrich,HaardtTent}.

We introduce matrix-weighted $Z$-spaces and tent spaces associated
with $\mathcal A_p$ weights and establish their dyadic
characterizations using reducing operators. We then characterize
homogeneous matrix-weighted Besov and Triebel--Lizorkin spaces
in terms of the corresponding $Z$-spaces and tent spaces,
respectively, via continuous Whitney averages. These characterizations
complement the discrete Lusin-average descriptions in
\cite[Definition 4.6 and Theorem 4.7]{BYYZ}.
The kernels satisfy Tauberian and finite moment conditions,
with the required moment order reflecting both the smoothness
and the $\mathcal A_p$-dimension of the weight.

We also identify the continuous duals of the matrix-weighted
$Z$-spaces and tent spaces and determine their real and complex
interpolants.
The results include duality for quasi-Banach $Z$-spaces and
complex interpolation between spaces whose matrix weights
need not commute. The interpolation results also apply to homogeneous
matrix-weighted Besov and Triebel--Lizorkin spaces.

The paper is organized as follows.
Section~\ref{sec:preliminaries} contains the preliminaries.
Section~\ref{sec:ZT-spaces} introduces the upper-half-space spaces
and establishes their dyadic characterizations and basic properties.
Section~\ref{sec:main} establishes Whitney average characterizations
of homogeneous matrix-weighted Besov and Triebel--Lizorkin spaces
in terms of the $Z$-spaces and tent spaces, respectively.
Section~\ref{sec:Z} identifies the continuous duals of the
$Z$-spaces and tent spaces.
Section~\ref{sec:interpolation} proves the interpolation results
and their consequences for homogeneous matrix-weighted Besov
and Triebel--Lizorkin spaces.

\section{Preliminaries}\label{sec:preliminaries}

\subsection{Notation and matrix weights}

We write $\mathbb N=\{1,2,\ldots\}$, $\mathbb N_0=\mathbb N\cup\{0\}$,
and $\Z$ for the integers.
Throughout this paper, let $n,m\in\mathbb N$, $\alpha\in\R$,
$p\in(0,\infty)$, and $q,r\in(0,\infty]$. We use $|\cdot|$ for
the Euclidean norm and $\|\cdot\|$ for the operator norm of matrices.
For a complex matrix $M$, $M^*$ denotes its conjugate transpose;
$I_m$ is the $m\times m$ identity matrix, and $I$ denotes the
identity map on the space under consideration. For a multi-index
$\beta\in\mathbb N_0^n$, $|\beta|=\beta_1+\cdots+\beta_n$ and
$\partial^\beta$ denotes the corresponding partial derivative.
The indicator of a measurable set $E$ is denoted by $\mathbf1_E$.
For $0<u\le\infty$, put
\[
 u'=\begin{cases}
 \infty,&0<u\le1,\\
 u/(u-1),&1<u<\infty,\\
 1,&u=\infty.
 \end{cases}
\]
The notation $U\lesssim V$ means that $U\leq CV$ for a positive
constant $C$ independent of the function, the scale, and the spatial
points. Write $U\asymp V$ if $U\lesssim V\lesssim U$.
Unless another domain is specified, $L^p$ denotes $L^p(\R^n)$
with respect to Lebesgue measure.

\begin{definition}\label{def:Schwartz-spaces}
The Schwartz space $\cS=\cS(\R^n)$ consists of all complex-valued
functions $\psi\in C^\infty(\R^n)$ for which
$\rho_{K,L}(\psi)<\infty$ for every $K,L\in\mathbb N_0$.
We use the following equivalent form of the Schwartz seminorms in
\cite[Definition 2.2.1]{GrafakosCFA}:
\[
 \rho_{K,L}(\psi)
 :=\max_{|\beta|\le L}\sup_{x\in\R^n}
 (1+|x|)^K|\partial^\beta\psi(x)|,
 \qquad K,L\in\mathbb N_0.
\]
For $\vec\psi=(\psi_1,\ldots,\psi_m)\in(\cS)^m$, define
\[
 \rho_{K,L}(\vec\psi)
 :=\max_{|\beta|\le L}\sup_{x\in\R^n}
 (1+|x|)^K\left(\sum_{\ell=1}^m
 |\partial^\beta\psi_\ell(x)|^2\right)^{1/2},
 \qquad K,L\in\mathbb N_0.
\]
Let $\cS'$ be the continuous linear dual of $\cS$, and let
$\cP$ be the space of polynomials on $\R^n$.
\end{definition}
For
$\vec f\in(\cS'/\cP)^m$, scalar kernels act componentwise on $\vec f$.
For $h\in\cS$, we write
\[
 \widehat h(\xi)=\mathcal Fh(\xi)
 :=\int_{\R^n}h(x)e^{-2\pi i x\cdot\xi}\,\mathrm d x,
 \qquad \xi\in\R^n,
\]
and denote the inverse Fourier transform by $\mathcal F^{-1}$.
For $\vec g=(g_1,\ldots,g_m)\in(\cS')^m$ and
$\vec\eta\in(\cS)^m$, use the pairing
\[
 \langle\vec g,\vec\eta\rangle
 :=\sum_{\ell=1}^m g_\ell(\overline{\eta_\ell}).
\]
If $\Phi\in\cS$ and $\widehat\Phi$ vanishes in a neighborhood of the
origin, then $\Phi*\vec f$ is independent of the representative of
$\vec f$ and is a smooth function of at most polynomial growth.
\Needspace{5\baselineskip}
For a measurable set $E\subset\R^n$ of finite positive measure, put
\[
 \fint_E h(x)\,\mathrm d x=\frac1{|E|}\int_E h(x)\,\mathrm d x.
\]
For a measurable scalar function $h$ on $\R^n$, the
uncentered Hardy--Littlewood maximal operator is defined by
\[
 \M h(x)=\sup_{B\ni x}\fint_B|h(y)|\,\mathrm d y,
\]
where the supremum is taken over all balls $B\subset\R^n$ containing $x$.

For $j\in\Z$, let
\[
 \D_j=\{2^{-j}(\nu+[0,1)^n):\nu\in\Z^n\},\qquad
 \D=\bigcup_{j\in\Z}\D_j.
\]
Denote by $Q_j(x)$ the unique cube in $\D_j$ containing $x$.
For $Q\in\D$, write $\ell(Q)$ for its side length and $x_Q$ for
its center.

\begin{definition}\label{def:mixed-norm-spaces}
Let $0<p<\infty$ and $0<q\le\infty$. The mixed-norm spaces
$L^p(\ell^q)$ and $\ell^q(L^p)$ consist of all sequences of measurable
functions $h_j:\R^n\to\C$, $j\in\Z$, for which the respective
quasi-norm below is finite:
\[
 \left\|\{h_j\}\right\|_{L^p(\ell^q)}
 =\left\|\left(\sum_{j\in\Z}|h_j|^q\right)^{1/q}\right\|_{L^p},
 \qquad
 \left\|\{h_j\}\right\|_{\ell^q(L^p)}
 =\left(\sum_{j\in\Z}\|h_j\|_{L^p}^q\right)^{1/q}.
\]
When $q=\infty$, the $\ell^q$ quasi-norm is replaced by a supremum.
We use
$L^p_x(L^q(\mathrm d t/t))$ for taking the scale norm on
$(0,\infty)$ before the spatial $L^p(\R^n)$ norm, and
$L^q(\mathrm d t/t;L^p_x)$ for the reverse order.
\end{definition}

For a positive definite Hermitian $m\times m$ matrix $A$, by the
spectral theorem, we have
\[
 A=U\operatorname{diag}(\lambda_1,\ldots,\lambda_m)U^*,
 \qquad \lambda_1,\ldots,\lambda_m>0,
\]
where $U$ is an $m\times m$ unitary matrix, that is,
$U^*U=UU^*=I_m$; see
\cite[Theorems 2.5.6(c) and 4.1.8]{HJMatrix}.
For $t\in\R$, define
\[
 A^t:=U\operatorname{diag}(\lambda_1^t,\ldots,\lambda_m^t)U^*.
\]
This definition is independent of the choice of $U$ and the order
of the eigenvalues; see \cite[pp.~407--408]{HJTopics}.

\Needspace{14\baselineskip}
\begin{definition}\label{def:Ap}
A matrix weight $W$ is a measurable function on $\R^n$ with values
in the positive definite Hermitian $m\times m$ matrices almost
everywhere, whose entries are locally integrable.
If $p\in(1,\infty)$, define
\[
 [W]_{\mathcal A_p}:=\sup_Q\fint_Q
 \left[\fint_Q
 \|W^{1/p}(x)W^{-1/p}(y)\|^{p'}\,\mathrm d y\right]^{p/p'}\mathrm d x.
\]
If $p\in(0,1]$, define
\[
 [W]_{\mathcal A_p}:=\sup_Q\esssup_{y\in Q}\fint_Q
           \|W^{1/p}(x)W^{-1/p}(y)\|^p\,\mathrm d x.
\]
Here and below these suprema are over cubes with sides parallel to
the coordinate axes. We write $W\in\mathcal A_p$ if $[W]_{\mathcal A_p}<\infty$.
\end{definition}

\begin{remark}\label{rem:Ap-monotonicity}
By \cite[Proposition 2.12]{BHYY}, the matrix weight classes in
Definition~\ref{def:Ap} satisfy $\mathcal A_p\subset \mathcal A_q$ whenever
$0<p<q<\infty$.
Thus a fixed weight in $\mathcal A_{\min\{p_1,p_2\}}$ is admissible at
both exponents $p_1$ and $p_2$.
\end{remark}

\subsection{Reducing operators and matrix estimates}

\Needspace{5\baselineskip}
\begin{definition}\label{def:reducing}
Let $0<p<\infty$ and $W\in\mathcal A_p$. A sequence $\mathbb A:=\{A_Q\}_{Q\in\D}$ of positive
definite Hermitian matrices is a sequence of reducing operators of
order $p$ for $W$ if
\begin{equation}\label{eq:reducing}
 |A_Q\vec v|\asymp
 \left(\fint_Q|W^{1/p}(x)\vec v|^p\,\mathrm d x\right)^{1/p},
 \qquad Q\in\D,\quad \vec v\in\C^m,
\end{equation}
where the positive equivalence constants depend only on $m$ and $p$.
For $\mathbb A=\{A_Q\}_{Q\in\D}$, set
\[
 \mathbb A_j(x)=\sum_{Q\in\D_j}A_Q\mathbf1_Q(x),\qquad j\in\Z.
\]
For each continuous scale $t>0$, there is a unique $j\in\Z$
such that $2^{-j}\leq t<2^{-j+1}$. We define
\begin{equation}\label{eq:At-definition}
 \mathbb A(t,x):=A_{Q_j(x)}
       =\sum_{Q\in\D_j}A_Q\mathbf1_Q(x)=\mathbb A_j(x),
 \qquad 2^{-j}\leq t<2^{-j+1}.
\end{equation}
\end{definition}

Reducing operators in Definition~\ref{def:reducing} exist for all
$p\in(0,\infty)$; see \cite{FR,BHYY}.

\Needspace{8\baselineskip}
We recall the $\mathcal A_p$-dimensions of Bu--Hyt\"onen--Yang--Yuan
(Definition~2.22 and Proposition~2.23(vi), \cite{BHYY}).

\begin{definition}\label{def:Ap-dimension}
For $0<p<\infty$ and a matrix weight $W$, we say that $W$ has
$\mathcal A_p$-dimension $d\geq0$ if a constant $C$ independent of
the cube $Q$ and $\lambda\geq1$ satisfies
\[
 \esssup_{y\in\lambda Q}\fint_Q
 \|W^{1/p}(x)W^{-1/p}(y)\|^p\,\mathrm d x\leq C\lambda^d
 \quad\text{if }0<p\leq1,
\]
and
\[
 \fint_Q\left[\fint_{\lambda Q}
 \|W^{1/p}(x)W^{-1/p}(y)\|^{p'}\,\mathrm d y\right]^{p/p'}\mathrm d x
 \leq C\lambda^d
 \quad\text{if }1<p<\infty.
\]
Here $\lambda Q$ is the concentric dilate of $Q$.
\end{definition}

By \cite[Proposition 2.26]{BHYY}, for $W\in\mathcal A_p$, there
exists $d\in[0,n)$ such that $W$ has $\mathcal A_p$-dimension $d$.
For $p>1$, $\widetilde W:=W^{-1/(p-1)}\in\mathcal A_{p'}$.
Thus $\widetilde W$ has some $\mathcal A_{p'}$-dimension
$\widetilde d\in[0,n)$. For $p\leq1$, put $\widetilde d=0$.
Set
\begin{equation}\label{eq:dimension-exponents}
 \delta_1=\frac d p,\qquad
 \delta_2=\begin{cases}\widetilde d/p',&p>1,\\0,&0<p\leq1,
 \end{cases}
 \qquad \Delta=\delta_1+\delta_2.
\end{equation}
Write $u_+=\max\{u,0\}$ for $u\in\R$.

\Needspace{5\baselineskip}
The following estimates follow from \cite[Lemma 2.28 and Corollary 2.31]{BHYY}.

\begin{lemma}\label{lem:matrix-comparison}
Let $p\in(0,\infty)$, $W\in\mathcal A_p$, and let
$\mathbb A:=\{A_Q\}_{Q\in\D}$ satisfy
\eqref{eq:reducing}.
With the exponents in \eqref{eq:dimension-exponents},
there exists $C>0$ such that
\begin{equation}\label{eq:comparison}
 \|\mathbb A_j(x)\mathbb A_k(z)^{-1}\| \leq C2^{\delta_1(j-k)_++\delta_2(k-j)_+} (1+2^{\min(j,k)}|x-z|)^\Delta,\qquad j,k\in\Z,\quad x,z\in\R^n.
\end{equation}
\end{lemma}

\Needspace{5\baselineskip}
The following result comes from \cite[Lemmas 3.2 and 3.3]{FR}.

\begin{lemma}\label{lem:matrix-inverse}
Let $p\in(0,\infty)$, $W\in\mathcal A_p$, and let $\mathbb A$
be as in Definition~\ref{def:reducing}.

\smallskip
\noindent$\mathrm{(i)}$ There exists $C_1>0$, independent of $Q$, such that, for every $Q\in\D$,
\begin{equation}\label{eq:inverse}
 \fint_Q\|A_QW^{-1/p}(x)\|\,\mathrm d x\leq C_1.
\end{equation}
If $0<p\leq1$, the constant $C_1$ can be chosen so that
\[
 \esssup_{x\in Q}\|A_QW^{-1/p}(x)\|\leq C_1,
 \qquad Q\in\D.
\]

\smallskip
\noindent$\mathrm{(ii)}$ The set
\begin{equation}\label{eq:good-set}
 E_Q:=\{x\in Q:\|A_QW^{-1/p}(x)\|\leq2C_1\}
\end{equation}
satisfies $|E_Q|\geq|Q|/2$.
\end{lemma}

\begin{remark}\label{rem:good-set}
If $0<p\leq1$, the set $E_Q$ in \eqref{eq:good-set} satisfies
$|Q\setminus E_Q|=0$ for every $Q\in\D$.
\end{remark}

\Needspace{5\baselineskip}
The relation between inverse reducing operators and the dual weight
is given in \cite[Section 2 and the discussion before Lemma 3.3]{RoudenkoDuality}.

\begin{lemma}\label{lem:Z-dual-reducing}
Let $1<p<\infty$, $W\in\mathcal A_p$, and let
$\mathbb A=\{A_Q\}_{Q\in\D}$ be reducing operators of order $p$ for $W$. Set
\begin{equation*}
 \widetilde W(x):=W(x)^{-p'/p},\qquad
 \mathbb A^{-1}:=\{A_Q^{-1}\}_{Q\in\D}.
\end{equation*}
Then $\widetilde W\in\mathcal A_{p'}$ and
\begin{equation*}
 |A_Q^{-1}\vec v|
 \asymp\left(\fint_Q|W^{-1/p}(x)\vec v|^{p'}\,\mathrm d x\right)^{1/p'},
 \qquad Q\in\D,\quad\vec v\in\C^m.
\end{equation*}
The positive equivalence constants depend only on $m,p$ and
$[W]_{\mathcal A_p}$, and are independent of $Q$ and $\vec v$.
\end{lemma}

\Needspace{5\baselineskip}
The following multiplier estimates are from \cite[Corollary 3.8]{FR}.

\begin{lemma}\label{lem:matrix-multiplier}
Let $p\in(0,\infty)$, $q\in(0,\infty]$, $W\in\mathcal A_p$,
and let $\mathbb A$ be as in Definition~\ref{def:reducing}.
For each $j\in\Z$, let $c_j:\R^n\to[0,\infty)$ be constant
on every cube of $\D_j$. Then
\begin{align}
 \left\|\left\{\sum_{j\in\Z}
 [\|W^{1/p}(\cdot)\mathbb A_j(\cdot)^{-1}\|c_j]^q
 \right\}^{1/q}\right\|_{L^p(\R^n)}
 &\lesssim\left\|\left(\sum_{j\in\Z}c_j^q\right)^{1/q}\right\|_{L^p(\R^n)},
 \label{eq:constant-F}\\
 \left\{\sum_{j\in\Z}
 \|\|W^{1/p}(\cdot)\mathbb A_j(\cdot)^{-1}\|c_j\|_{L^p(\R^n)}^q
 \right\}^{1/q}
 &\lesssim\left(\sum_{j\in\Z}\|c_j\|_{L^p(\R^n)}^q\right)^{1/q}.
 \label{eq:constant-B}
\end{align}
\end{lemma}

\subsection{Homogeneous Besov and Triebel--Lizorkin spaces}

Choose $\varphi\in\cS$ such that
\begin{equation}\label{eq:partition}
 \begin{gathered}
 \supp\widehat\varphi\subset\{\xi:1/2\leq|\xi|\leq2\},\qquad
 \sum_{j\in\Z}\widehat\varphi(2^{-j}\xi)=1\quad(\xi\neq0),\\
 |\widehat\varphi(\xi)|\geq c>0
 \quad\hbox{if }3/5\leq|\xi|\leq5/3.
 \end{gathered}
\end{equation}
\begin{definition}\label{def:spaces}
Let $\alpha\in\R$, $p\in(0,\infty)$, $q\in(0,\infty]$,
$W\in\mathcal A_p$, and let $\mathbb A$ be a sequence of reducing
operators of order $p$ for $W$ as in Definition~\ref{def:reducing}.
Let $\varphi\in\cS$ satisfy \eqref{eq:partition}, and set
$\varphi_j(x)=2^{jn}\varphi(2^jx)$ for every $j\in\Z$.
The homogeneous matrix-weighted Triebel--Lizorkin space
$\dot F^{\alpha}_{p,q}(W)$ consists of all
$\vec f\in(\cS'/\cP)^m$ such that
\[
 \|\vec f\|_{\dot F^{\alpha}_{p,q}(W)}
 :=\left\|\left(\sum_{j\in\Z}
 [2^{j\alpha}|W^{1/p}(\cdot)(\varphi_j*\vec f)(\cdot)|]^q
 \right)^{1/q}\right\|_{L^p(\R^n)}<\infty.
\]
The homogeneous matrix-weighted Besov space $\dot B^{\alpha}_{p,q}(W)$
consists of all such $\vec f$ for which
\[
 \|\vec f\|_{\dot B^{\alpha}_{p,q}(W)}
 :=\left(\sum_{j\in\Z}
 [2^{j\alpha}\|W^{1/p}(\cdot)(\varphi_j*\vec f)(\cdot)\|_{L^p(\R^n)}]^q
 \right)^{1/q}<\infty.
\]
The corresponding spaces defined using the reducing operators
have the following quasi-norms:
\begin{align*}
 \|\vec f\|_{\dot F^\alpha_{p,q}(\mathbb A)}
 &:=\left\|\left(\sum_{j\in\Z}
 [2^{j\alpha}|\mathbb A_j(\cdot)(\varphi_j*\vec f)(\cdot)|]^q
 \right)^{1/q}\right\|_{L^p},\\
 \|\vec f\|_{\dot B^\alpha_{p,q}(\mathbb A)}
 &:=\left(\sum_{j\in\Z}
 [2^{j\alpha}\|\mathbb A_j(\cdot)(\varphi_j*\vec f)(\cdot)
 \|_{L^p}]^q\right)^{1/q}.
\end{align*}
We write $\vec f\in\dot F^\alpha_{p,q}(\mathbb A)$ or
$\vec f\in\dot B^\alpha_{p,q}(\mathbb A)$ when the respective quasi-norm
is finite.
\end{definition}

\Needspace{5\baselineskip}
The equivalence of the weighted and reducing-operator quasi-norms
is \cite[Theorem 3.1]{FR} for Triebel--Lizorkin spaces and
\cite[Lemma 5.7]{FRBesov} for Besov spaces.

\begin{lemma}\label{lem:known-equivalence}
Let $\alpha\in\R$, $0<p<\infty$, $0<q\le\infty$, and
$W\in\mathcal A_p$. Let $\mathbb A$ be as in
Definition~\ref{def:reducing}.
Let $\varphi\in\cS$ satisfy \eqref{eq:partition}, and use the
quasi-norms in Definition~\ref{def:spaces}.
Then, for every $\vec f\in(\cS'/\cP)^m$,
\begin{equation}\label{eq:known-equivalence}
 \|\vec f\|_{\dot F^{\alpha}_{p,q}(W)}
 \asymp\|\vec f\|_{\dot F^{\alpha}_{p,q}(\mathbb A)},\qquad
 \|\vec f\|_{\dot B^{\alpha}_{p,q}(W)}
 \asymp\|\vec f\|_{\dot B^{\alpha}_{p,q}(\mathbb A)}.
\end{equation}
\end{lemma}

\Needspace{13\baselineskip}
\subsection{Maximal estimates}\label{subsec:maximal}

The following bandlimited Peetre estimate is the scalar estimate used in
\cite[proof of Theorem 2.4]{FR}; see also
\cite[(3.39), proof of Theorem 3.14]{WYZ}.

\begin{lemma}\label{lem:bandlimited}
Let $j\in\Z$, $c,a\in(0,\infty)$, and $\theta\in(0,1]$. Let
$\vec u\in(\cS')^m$ have Fourier support in $B(0,c2^j)$ and let
$T\in\C^{m\times m}$ be a constant matrix. Then, for every $x\in\R^n$,
\[
 \left[\sup_{y\in\R^n}\frac{|T\vec u(y)|}{(1+2^j|x-y|)^a}\right]^\theta
 \lesssim\int_{\R^n}
 \frac{2^{jn}|T\vec u(z)|^\theta}{(1+2^j|x-z|)^{a\theta}}\,\mathrm d z.
\]
The constant is independent of $T$ and $\vec u$.
\end{lemma}

For $N>n$, $j\in\Z$, $x\in\R^n$, and any nonnegative
measurable function $h$ on $\R^n$, one has
\begin{equation}\label{eq:decay-max}
 \int_{\R^n}\frac{2^{jn}h(z)}{(1+2^j|x-z|)^N}\,\mathrm d z
 \lesssim\M h(x).
\end{equation}
This is the usual maximal estimate for an integrable radially
decreasing kernel; see also its dyadic-cube form in
\cite[Lemma 3.7]{WYZ}.

\Needspace{5\baselineskip}
The following estimate is a consequence of
Lemma~\ref{lem:bandlimited}, \eqref{eq:comparison}
and \eqref{eq:decay-max}.

\begin{lemma}\label{lem:peetre-point}
Let $0<p<\infty$, $W\in\mathcal A_p$, and let $\mathbb A$ be as in
Definition~\ref{def:reducing}. Let $\theta\in(0,1]$ and
$a>\Delta+n/\theta$, where $\Delta$ is as in
\eqref{eq:dimension-exponents}.
Let $\varphi\in\cS$ satisfy \eqref{eq:partition}, and set
$\varphi_j(x)=2^{jn}\varphi(2^jx)$ for $j\in\Z$.
Then, for every $j\in\Z$, $x\in\R^n$, and
$\vec f\in(\cS'/\cP)^m$,
\begin{equation}\label{eq:peetre-point}
 \sup_{y\in\R^n}
 \frac{|\mathbb A_j(x)(\varphi_j*\vec f)(y)|}{(1+2^j|x-y|)^a}
 \lesssim\left[\M\big(|\mathbb A_j(\cdot)
               (\varphi_j*\vec f)(\cdot)|^\theta\big)(x)
          \right]^{1/\theta}.
\end{equation}
\end{lemma}

\Needspace{5\baselineskip}
The following is the Fefferman--Stein inequality; see
\cite[Theorem 5.6.6 and Remark 5.6.7]{GrafakosCFA} for (i).
For its continuous-parameter version in (ii), see
\cite[Claim 1, p.~563]{Moritoh95}.

\begin{lemma}\label{lem:FS}
Let $p\in(1,\infty)$ and $q\in(1,\infty]$.

\smallskip
\noindent$\mathrm{(i)}$ For every sequence
of nonnegative measurable functions $\{h_j\}_{j\in\Z}$ on $\R^n$,
\[
 \left\|\left(\sum_{j\in\Z}[\M h_j]^q\right)^{1/q}
 \right\|_{L^p}
 \lesssim\left\|\left(\sum_{j\in\Z}h_j^q\right)^{1/q}
 \right\|_{L^p}.
\]

\smallskip
\noindent$\mathrm{(ii)}$ For every nonnegative measurable $h$ on $(0,\infty)\times\R^n$,
\begin{equation*}
 \left\|\left(\int_0^\infty
 [\M(h(t,\cdot))]^q\frac{\mathrm d t}{t}
 \right)^{1/q}\right\|_{L^p}
 \lesssim\left\|\left(\int_0^\infty
 h(t,\cdot)^q\frac{\mathrm d t}{t}\right)^{1/q}\right\|_{L^p}.
\end{equation*}
\end{lemma}

\section{Matrix-weighted \texorpdfstring{$Z$}{Z}-spaces and tent spaces}
\label{sec:ZT-spaces}

In this section, we introduce matrix-weighted $Z$-spaces and tent
spaces and establish their dyadic characterizations and basic
properties.

\subsection{Definitions}
\label{subsec:ZT-definitions}

The unweighted three-parameter $Z$-spaces were introduced by
Auscher, Bechtel and Haardt \cite{ABH}.
For the corresponding tent spaces with Whitney averages, see
\cite{HuangTent,ABH,HaardtTent}.

Let $\R^{n+1}_+=(0,\infty)\times\R^n$. Throughout
Sections~\ref{sec:ZT-spaces}--\ref{sec:interpolation}, $\vec F$ and
$\vec G$ denote measurable $\C^m$-valued functions on
$\R^{n+1}_+$, identified when they agree almost everywhere.

\Needspace{5\baselineskip}
\begin{definition}\label{def:Z}
Let $\alpha\in\R$, $p\in(0,\infty)$, and $q,r\in(0,\infty]$.
Let $\mathbb A:=\{A_Q\}_{Q\in\D}$ be a family of positive definite
Hermitian $m\times m$ matrices satisfying the comparison estimate
\eqref{eq:comparison} for some $\delta_1,\delta_2\geq0$ and
$\Delta=\delta_1+\delta_2$.
Define $\mathbb A_j(x)$ and $\mathbb A(t,x)$ by the formulas in
Definition~\ref{def:reducing}. The space $Z^{p,q,r}_\alpha(\mathbb A)$
consists of all measurable
$\vec F:\R^{n+1}_+\to\C^m$ such that
\begin{equation*}
 \|\vec F\|_{Z^{p,q,r}_\alpha(\mathbb A)}
 :=\left\{\int_0^\infty
 \left\|\left[\fint_{t/2}^{t}\fint_{B(\cdot,t)}
 |s^{-\alpha}\mathbb A(t,\cdot)\vec F(s,y)|^r\,\mathrm d y\,\mathrm d s
 \right]^{1/r}\right\|_{L^p(\R^n)}^q\frac{\mathrm d t}{t}\right\}^{1/q}
 <\infty.
\end{equation*}
If $W\in\mathcal A_p$, the space $Z^{p,q,r}_\alpha(W)$ consists of all
$\vec F$ such that
\begin{equation*}
 \|\vec F\|_{Z^{p,q,r}_\alpha(W)}
 :=\left\{\int_0^\infty
 \left\|\left[\fint_{t/2}^{t}\fint_{B(\cdot,t)}
 |s^{-\alpha}W^{1/p}(\cdot)\vec F(s,y)|^r\,\mathrm d y\,\mathrm d s
 \right]^{1/r}\right\|_{L^p(\R^n)}^q\frac{\mathrm d t}{t}\right\}^{1/q}
 <\infty,
\end{equation*}
with the usual modifications if $q$ and/or $r$ are infinite.
\end{definition}

\Needspace{5\baselineskip}
\begin{definition}\label{def:T}
Let $\alpha\in\R$, $p\in(0,\infty)$, $q,r\in(0,\infty]$,
and let $\mathbb A$ be as in Definition~\ref{def:Z}.
The tent space $T^{p,q,r}_\alpha(\mathbb A)$ consists of all
measurable $\vec F:\R^{n+1}_+\to\C^m$ such that
\[
 \|\vec F\|_{T^{p,q,r}_\alpha(\mathbb A)}
 :=\left\|\left\{\int_0^\infty
 \left[\fint_{t/2}^{t}\fint_{B(\cdot,t)}
 |s^{-\alpha}\mathbb A(t,\cdot)\vec F(s,y)|^r\,\mathrm d y\,\mathrm d s
 \right]^{q/r}\frac{\mathrm d t}{t}\right\}^{1/q}\right\|_{L^p(\R^n)}
 <\infty.
\]
If $W\in\mathcal A_p$, define $T^{p,q,r}_\alpha(W)$ by
\begin{equation*}
 \|\vec F\|_{T^{p,q,r}_\alpha(W)}
 :=\left\|\left\{\int_0^\infty
 \left[\fint_{t/2}^{t}\fint_{B(\cdot,t)}
 |s^{-\alpha}W^{1/p}(\cdot)\vec F(s,y)|^r\,\mathrm d y\,\mathrm d s
 \right]^{q/r}\frac{\mathrm d t}{t}\right\}^{1/q}\right\|_{L^p(\R^n)}
 <\infty,
\end{equation*}
with the usual modifications if $q$ and/or $r$ are infinite.
\end{definition}

\Needspace{4\baselineskip}
\subsection{Dyadic characterizations and basic properties}
\label{subsec:ZT-dyadic}

For $Q\in\D$, put
\begin{equation}\label{eq:Z-boxes}
 \widehat Q=[\ell(Q),2\ell(Q))\times Q,\qquad
 \mathrm d \nu_Q(s,y)=\frac{\mathbf1_{\widehat Q}(s,y)}{|Q|\log_e 2}
             \frac{\mathrm d y\,\mathrm d s}{s}.
\end{equation}
The Whitney boxes $\widehat Q$, $Q\in\D$, partition $\R^{n+1}_+$, and $\nu_Q$ is a probability measure on $\widehat Q$.

We first recall that replacing the time interval $(t/2,t)$ by $(at,bt)$
and the ball $B(x,t)$ by $B(x,ct)$, where $0<a<b<\infty$ and $c>0$
are fixed, yields equivalent quasi-norms for the unweighted
$Z$-spaces and tent spaces. The following part $\mathrm{(i)}$
follows from \cite[Observations 2.1 and 2.4]{HuangTent} and
\cite[Lemma 3.3]{ABH}; part $\mathrm{(ii)}$ is
\cite[Proposition 3.7]{ABH} at smoothness zero. These scalar results
are applied to $(s,y)\mapsto s^{-\alpha}|\vec G(s,y)|$.

\Needspace{4\baselineskip}
\begin{lemma}\label{lem:ZT-Whitney-parameters}
Let $\alpha\in\R$, $0<p<\infty$, $0<q,r\le\infty$,
$0<a<b<\infty$, $c>0$, and let $\mathbb I=\{I_m\}_{Q\in\D}$.
For every measurable $\vec G:\R^{n+1}_+\to\C^m$, the following
equivalences hold.

\smallskip
\noindent$\mathrm{(i)}$ For the tent space,
\[
 \left\|\left\{\int_0^\infty
 \left[\fint_{at}^{bt}\fint_{B(\cdot,ct)}
 |s^{-\alpha}\vec G(s,y)|^r\,\mathrm d y\,\mathrm d s\right]^{q/r}
 \frac{\mathrm d t}{t}\right\}^{1/q}\right\|_{L^p}
 \asymp\|\vec G\|_{T^{p,q,r}_\alpha(\mathbb I)}.
\]

\smallskip
\noindent$\mathrm{(ii)}$ For the $Z$-space,
\[
 \left\{\int_0^\infty
 \left\|\left[\fint_{at}^{bt}\fint_{B(\cdot,ct)}
 |s^{-\alpha}\vec G(s,y)|^r\,\mathrm d y\,\mathrm d s\right]^{1/r}
 \right\|_{L^p}^{q}\frac{\mathrm d t}{t}\right\}^{1/q}
 \asymp\|\vec G\|_{Z^{p,q,r}_\alpha(\mathbb I)}.
\]
The equivalence constants depend only on $n,p,q,r,a,b,c$.
\end{lemma}

\Needspace{8\baselineskip}
We next recall the unweighted dyadic characterizations. Part
$\mathrm{(ii)}$ follows from \cite[Proposition 3.12]{ABH}, after rescaling
the time variable by $2$ and using
Lemma~\ref{lem:ZT-Whitney-parameters}~$\mathrm{(ii)}$. For the tent-space
result, see \cite[Proposition 3.13]{HaardtTent}; we give an alternative short proof below.

\begin{lemma}\label{lem:ZT-unweighted-dyadic}
Let $\alpha\in\R$, $0<p<\infty$, $0<q,r\le\infty$, and let
$\mathbb I=\{I_m\}_{Q\in\D}$. Let $\widehat Q$ and $\nu_Q$ be
as in \eqref{eq:Z-boxes}. For every measurable
$\vec G:\R^{n+1}_+\to\C^m$, the following equivalences hold.

\smallskip
\noindent$\mathrm{(i)}$ For the tent space,
\[
 \|\vec G\|_{T^{p,q,r}_\alpha(\mathbb I)}
 \asymp\left\|\left\{\sum_{j\in\Z}\sum_{Q\in\D_j}\mathbf1_Q
 [2^{j\alpha}\|\vec G\|_{L^r(\widehat Q,\nu_Q;\C^m)}]^q
 \right\}^{1/q}\right\|_{L^p}.
\]

\smallskip
\noindent$\mathrm{(ii)}$ For the $Z$-space,
\[
 \|\vec G\|_{Z^{p,q,r}_\alpha(\mathbb I)}
 \asymp\left\{\sum_{j\in\Z}2^{j\alpha q}
 \left[\sum_{Q\in\D_j}|Q|
 \|\vec G\|_{L^r(\widehat Q,\nu_Q;\C^m)}^p\right]^{q/p}
 \right\}^{1/q}.
\]
\end{lemma}

\begin{proof}
\noindent$\mathrm{(i)}$ Let $0<\theta<\min\{1,p,q,r\}$. For $Q\in\D$, recall that
\[
 \widehat Q=[\ell(Q),2\ell(Q))\times Q,\qquad
 \mathrm d\nu_Q(s,y)=\frac{\mathrm d y\,\mathrm d s}{s|Q|\log_e2}
 \quad\text{on }\widehat Q.
\]
If $k\in\Z$ and $2^{-k}\le t<2^{-k+1}$, then
\[
 (t/2,t)\times B(x,t)
 \subset\bigcup_{j=k}^{k+1}
 \bigcup_{\substack{Q\in\D_j\\Q\cap B(x,t)\ne\varnothing}}\widehat Q.
\]
For each cube $Q$ in this union, $Q\subset B(x,(1+\sqrt n)t)$
and $|Q|\asymp t^n$. Then we have
\[
 [2^{j\alpha}\|\vec G\|_{L^r(\widehat Q,\nu_Q;\C^m)}]^\theta
 \lesssim\M\left(\sum_{R\in\D_j}\mathbf1_R
 [2^{j\alpha}\|\vec G\|_{L^r(\widehat R,\nu_R;\C^m)}]^\theta\right)(x).
\]
Since $s\asymp t$ on these Whitney boxes, we have
\[
 \frac{1}{(t/2)|B(x,t)|}\asymp\frac{1}{s|Q|\log_e2}.
\]
Using the preceding estimates and $\theta<r$, we obtain
\[
 \left[\fint_{t/2}^t\fint_{B(x,t)}
 |s^{-\alpha}\vec G(s,y)|^r\,\mathrm d y\,\mathrm d s\right]^{\theta/r}
 \lesssim\sum_{j=k}^{k+1}\M\left(
 \sum_{Q\in\D_j}\mathbf1_Q
 [2^{j\alpha}\|\vec G\|_{L^r(\widehat Q,\nu_Q;\C^m)}]^\theta\right)(x).
\]
Using Lemma~\ref{lem:FS} with $p/\theta>1$ and $q/\theta>1$, we obtain
\begin{align*}
 \|\vec G\|_{T^{p,q,r}_\alpha(\mathbb I)}
 &\lesssim\left\|\left\{\sum_{j\in\Z}
 \left[\M\left(\sum_{Q\in\D_j}\mathbf1_Q
 [2^{j\alpha}\|\vec G\|_{L^r(\widehat Q,\nu_Q;\C^m)}]^\theta
 \right)\right]^{q/\theta}\right\}^{1/q}\right\|_{L^p}\\
 &\lesssim\left\|\left\{\sum_{j\in\Z}\sum_{Q\in\D_j}\mathbf1_Q
 [2^{j\alpha}\|\vec G\|_{L^r(\widehat Q,\nu_Q;\C^m)}]^q
 \right\}^{1/q}\right\|_{L^p}.
\end{align*}
Conversely, for $Q\in\D$, $x\in Q$ and $\ell(Q)\le t<2\ell(Q)$,
$\widehat Q\subset(t/2,2t)\times B(x,2\sqrt n\,t)$. Thus
\begin{align*}
 &\left\|\left\{\sum_{j\in\Z}\sum_{Q\in\D_j}\mathbf1_Q
 [2^{j\alpha}\|\vec G\|_{L^r(\widehat Q,\nu_Q;\C^m)}]^q
 \right\}^{1/q}\right\|_{L^p}\\
 &\quad\lesssim\left\|\left\{\int_0^\infty
 \left[\fint_{t/2}^{2t}\fint_{B(\cdot,2\sqrt n\,t)}
 |s^{-\alpha}\vec G(s,y)|^r\,\mathrm d y\,\mathrm d s\right]^{q/r}
 \frac{\mathrm d t}{t}\right\}^{1/q}\right\|_{L^p}
 \asymp\|\vec G\|_{T^{p,q,r}_\alpha(\mathbb I)}.
\end{align*}
The last equivalence follows from
Lemma~\ref{lem:ZT-Whitney-parameters}~$\mathrm{(i)}$, with $a=1/2$,
$b=2$ and $c=2\sqrt n$.
\end{proof}

The next lemma transfers these characterizations to matrix families
satisfying \eqref{eq:comparison}.

\begin{lemma}\label{lem:Z-dyadic-A}
Let $\alpha\in\R$, $0<p<\infty$, $0<q,r\le\infty$, and
let $\mathbb A$ be as in Definition~\ref{def:Z}.
For every measurable $\vec F:\R^{n+1}_+\to\C^m$, the following equivalences hold.

\smallskip
\noindent$\mathrm{(i)}$ For the tent space,
\begin{equation}\label{eq:T-dyadic-A}
 \|\vec F\|_{T^{p,q,r}_\alpha(\mathbb A)}
\asymp\left\|\left\{\sum_{j\in\Z}\sum_{Q\in\D_j}
\mathbf1_Q\,2^{j\alpha q}
\left(\frac1{|Q|\log_e2}\!\int_{2^{-j}}^{2^{-j+1}}\int_Q
|A_Q\vec F(s,y)|^r\,\mathrm d y\,\frac{\mathrm d s}{s}\right)^{q/r}
\right\}^{1/q}\right\|_{L^p}.
\end{equation}

\smallskip
\noindent$\mathrm{(ii)}$ For the $Z$-space,
\begin{equation}\label{eq:Z-dyadic-A}
 \|\vec F\|_{Z^{p,q,r}_\alpha(\mathbb A)}
\asymp\left\{\sum_{j\in\Z}2^{j\alpha q}
\left[\sum_{Q\in\D_j}|Q|
\left(\frac1{|Q|\log_e2}\!\int_{2^{-j}}^{2^{-j+1}}\int_Q
|A_Q\vec F(s,y)|^r\,\mathrm d y\,\frac{\mathrm d s}{s}\right)^{p/r}
\right]^{q/p}\right\}^{1/q}.
\end{equation}
\end{lemma}

\begin{proof}
Let $\mathbb I=\{I_m\}_{Q\in\D}$. For $t>0$, $x\in\R^n$,
$t/2<s<t$ and $y\in B(x,t)$, the dyadic levels of $s$ and $t$
differ by at most one. Using \eqref{eq:comparison}, we obtain
\begin{equation}\label{eq:local-matrix-comparison}
 \|\mathbb A(t,x)\mathbb A(s,y)^{-1}\|
 +\|\mathbb A(s,y)\mathbb A(t,x)^{-1}\|\lesssim1.
\end{equation}
Set $\vec H(s,y):=\mathbb A(s,y)\vec F(s,y)$. Then, for $t/2<s<t$
and $y\in B(x,t)$, we have
\[
 |\mathbb A(t,x)\vec F(s,y)|\asymp|\vec H(s,y)|.
\]
Taking the corresponding averages and quasi-norms yields
\[
 \|\vec F\|_{T^{p,q,r}_\alpha(\mathbb A)}
 \asymp\|\vec H\|_{T^{p,q,r}_\alpha(\mathbb I)},\qquad
 \|\vec F\|_{Z^{p,q,r}_\alpha(\mathbb A)}
 \asymp\|\vec H\|_{Z^{p,q,r}_\alpha(\mathbb I)}.
\]
Since $\vec H=A_Q\vec F$ on $\widehat Q$,
Lemma~\ref{lem:ZT-unweighted-dyadic} proves $\mathrm{(i)}$ and $\mathrm{(ii)}$.
\end{proof}

\Needspace{5\baselineskip}
The following lemma shows that a matrix weight and its reducing
operators define equivalent quasi-norms on tent spaces and $Z$-spaces.

\begin{lemma}\label{lem:Z-W-A}
Let $\alpha\in\R$, $p\in(0,\infty)$, $q,r\in(0,\infty]$,
$W\in\mathcal A_p$, and let $\mathbb A$ be reducing operators of order $p$ for
$W$. For every measurable $\vec F:\R^{n+1}_+\to\C^m$, the following equivalences hold.

\smallskip
\noindent$\mathrm{(i)}$ For the tent space,
\begin{equation}\label{eq:T-W-A}
 \|\vec F\|_{T^{p,q,r}_\alpha(W)}
 \asymp\|\vec F\|_{T^{p,q,r}_\alpha(\mathbb A)}.
\end{equation}

\smallskip
\noindent$\mathrm{(ii)}$ For the $Z$-space,
\begin{equation}\label{eq:Z-W-A}
 \|\vec F\|_{Z^{p,q,r}_\alpha(W)}
 \asymp\|\vec F\|_{Z^{p,q,r}_\alpha(\mathbb A)}.
\end{equation}
\end{lemma}

\begin{proof}
\noindent$\mathrm{(i)}$ \emph{Step 1.}
Fix $k\in\Z$. For $Q\in\D_k$ and $j\in\{k,k+1\}$, let
\[
 \mathcal N_{k,j}(Q)=\{R\in\D_j:
 \operatorname{dist}(Q,R)\le2^{-k+1}\}.
\]
For $x\in Q\in\D_k$ and
$t\in [2^{-k},2^{-k+1})$,
\[
 |W^{1/p}(x)\vec F(s,y)|
 \leq\|W^{1/p}(x)A_Q^{-1}\|\,|A_Q\vec F(s,y)|.
\]
For $R\in\mathcal N_{k,j}(Q)$, $\|A_QA_R^{-1}\|\lesssim1$ by
\eqref{eq:comparison}, and
$B(x,t)\subset\bigcup_{R\in\mathcal N_{k,j}(Q)}R$.
Splitting $(t/2,t)$ at $2^{-k}$ and using
$s\asymp t\asymp2^{-j}$ and $|R|\asymp t^n$, we obtain
\begin{align}
 &\left[\fint_{t/2}^{t}\fint_{B(x,t)} |s^{-\alpha}W^{1/p}(x)\vec F(s,y)|^r\,\mathrm d y\,\mathrm d s \right]^{1/r} \lesssim\sum_{Q\in\D_k}\mathbf1_Q(x) \|W^{1/p}(x)A_Q^{-1}\|\notag\\
 &\quad \times\sum_{j=k}^{k+1} \sum_{R\in\mathcal N_{k,j}(Q)}2^{j\alpha} \left(\frac1{|R|\log_e 2}\!\int_{2^{-j}}^{2^{-j+1}}\int_R |A_R\vec F(s,y)|^r\,\mathrm d y\,\frac{\mathrm d s}{s}\right)^{1/r}.\label{eq:Z-cover-W}
 \end{align}
Choose an integer $L\ge0$ with $\sqrt n\,2^{-L}<1/4$, and
let $\mathcal C_L(Q)=\{R\in\D_{j+L}:R\subset Q\}$ for
$Q\in\D_j$. Divide $[2^{-j},2^{-j+1})$ into $S_{j,1}$ and
$S_{j,2}$, and define the auxiliary intervals $I_{j,1}$ and
$I_{j,2}$ as follows:
\[
 \begin{gathered}
 S_{j,1}=[2^{-j},3\cdot2^{-j}/2),\quad
 S_{j,2}=[3\cdot2^{-j}/2,2^{-j+1}),\\
 I_{j,1}=(7\cdot2^{-j}/4,15\cdot2^{-j}/8),\quad
 I_{j,2}=(9\cdot2^{-j}/4,5\cdot2^{-j}/2).
 \end{gathered}
\]
For $i\in\{1,2\}$, $t\in I_{j,i}$,
$R\in\mathcal C_L(Q)$ and $z\in R$, we have
\begin{equation}\label{eq:Z-box-inclusion}
 S_{j,i}\subset(t/2,t),\qquad R\subset B(z,t).
\end{equation}
With $C_1$ from \eqref{eq:inverse}, for each
$R\in\mathcal C_L(Q)$ put
\[
 E_R=\{z\in R:\|A_RW^{-1/p}(z)\|\leq2C_1\},
 \qquad |E_R|\geq|R|/2,
\]
For $z\in E_R$, using $R\subset Q$ with the fixed level difference
$L$, \eqref{eq:comparison} and the definition of $E_R$, we obtain
\[
 |A_Q\vec F(s,y)| \leq\|A_QA_R^{-1}\|\,\|A_RW^{-1/p}(z)\| |W^{1/p}(z)\vec F(s,y)| \lesssim|W^{1/p}(z)\vec F(s,y)|.
\]
For $t\in I_{j,i}$, using \eqref{eq:Z-box-inclusion}, we obtain,
for $z\in E_R$,
\begin{equation}\label{eq:Z-small-W}
 2^{j\alpha}\left(\int_{S_{j,i}}\fint_R
 |A_Q\vec F(s,y)|^r\,\frac{\mathrm d y\,\mathrm d s}{s\log_e 2}\right)^{1/r}
 \lesssim\left[\fint_{t/2}^{t}\fint_{B(z,t)}
 |s^{-\alpha}W^{1/p}(z)\vec F(s,y)|^r\,\mathrm d y\,\mathrm d s\right]^{1/r}.
\end{equation}

\medskip
\noindent\emph{Step 2.}
Fix $\theta\in(0,\min\{1,p,q,r\})$.
The sum after $\|W^{1/p}(x)A_Q^{-1}\|$ in
\eqref{eq:Z-cover-W} is constant on $Q$ and independent of $t$.
Integrate over $[2^{-k},2^{-k+1})$, sum in $k$, and apply
\eqref{eq:constant-F} to obtain
\begin{align*}
 &\|\vec F\|_{T^{p,q,r}_\alpha(W)}\\
 &\quad\lesssim\left\|\left\{\sum_{k\in\Z}
 \left[\sum_{Q\in\D_k}\mathbf1_Q
 \sum_{j=k}^{k+1}\sum_{R\in\mathcal N_{k,j}(Q)}2^{j\alpha}
 \left(\frac1{|R|\log_e 2}\!\int_{2^{-j}}^{2^{-j+1}}\int_R
 |A_R\vec F(s,y)|^r\,\mathrm d y\,\frac{\mathrm d s}{s}\right)^{1/r}
 \right]^q\right\}^{1/q}\right\|_{L^p}.
\end{align*}
For $k\in\Z$, $x\in Q\in\D_k$, $j\in\{k,k+1\}$, and
$R\in\mathcal N_{k,j}(Q)$, we have
\[
 \left[2^{j\alpha}\|A_R\vec F\|_{L^r(\widehat R,\nu_R;\C^m)}\right]^\theta
 \lesssim\M\left(\sum_{P\in\D_j}\mathbf1_P
 \left[2^{j\alpha}\|A_P\vec F\|_{L^r(\widehat P,\nu_P;\C^m)}\right]^\theta\right)(x).
\]
Summing in $R$ and using Lemma~\ref{lem:FS} with exponents
$p/\theta$ and $q/\theta$, together with \eqref{eq:T-dyadic-A}, we obtain
\begin{align*}
 &\|\vec F\|_{T^{p,q,r}_\alpha(W)} \lesssim\left\|\left\{\sum_{j\in\Z}\sum_{Q\in\D_j}
 \mathbf1_Q\left[2^{j\alpha}
 \left(\frac1{|Q|\log_e 2}\!\int_{2^{-j}}^{2^{-j+1}}\int_Q
 |A_Q\vec F(s,y)|^r\,\mathrm d y\,\frac{\mathrm d s}{s}\right)^{1/r}
 \right]^q\right\}^{1/q}\right\|_{L^p}\\
 &\quad\asymp\|\vec F\|_{T^{p,q,r}_\alpha(\mathbb A)}.
\end{align*}

\medskip
\noindent\emph{Step 3.}
Use the sets $E_R\subset R$, with $|E_R|\geq|R|/2$, and
\eqref{eq:Z-small-W} for $z\in E_R$ and $t\in I_{j,i}$.
Using the quasi-triangle inequality for the finitely many
subcubes and integrating the $\theta$th power of
\eqref{eq:Z-small-W} over $E_R$, we obtain, for $Q\in\D_j$,
\begin{align}\label{eq:T-good-set-mean}
 &\left[2^{j\alpha}
 \left(\frac1{|Q|\log_e 2}\int_{S_{j,i}}\int_Q
 |A_Q\vec F(s,y)|^r\,\mathrm d y\,\frac{\mathrm d s}{s}\right)^{1/r}
 \right]^\theta\notag\\
 &\quad\lesssim\sum_{R\in\mathcal C_L(Q)}\frac{|R|}{|Q|}
 \left[2^{j\alpha}\left(\frac1{|R|\log_e2}
 \int_{S_{j,i}}\int_R|A_Q\vec F(s,y)|^r
 \,\mathrm d y\,\frac{\mathrm d s}{s}\right)^{1/r}\right]^\theta\notag\\
 &\quad\lesssim\frac1{|Q|}
 \sum_{R\in\mathcal C_L(Q)}\frac{|R|}{|E_R|}
 \int_{E_R}\left[\fint_{t/2}^{t}\fint_{B(z,t)}
 |s^{-\alpha}W^{1/p}(z)\vec F(s,y)|^r\,\mathrm d y\,\mathrm d s
 \right]^{\theta/r}\mathrm d z\notag\\
 &\quad\lesssim\fint_Q
 \left[\fint_{t/2}^{t}\fint_{B(z,t)}
 |s^{-\alpha}W^{1/p}(z)\vec F(s,y)|^r\,\mathrm d y\,\mathrm d s
 \right]^{\theta/r}\mathrm d z.
\end{align}
Using \eqref{eq:T-good-set-mean}, we obtain, for $t\in I_{j,i}$,
\begin{align}\label{eq:T-lower-W-point}
 &\sum_{Q\in\D_j}\mathbf1_Q(x)\left[2^{j\alpha}
 \left(\frac1{|Q|\log_e 2}\int_{S_{j,i}}\int_Q
 |A_Q\vec F(s,y)|^r\,\mathrm d y\,\frac{\mathrm d s}{s}\right)^{1/r}
 \right]^\theta\notag\\
 &\quad\lesssim\M\left(
 \left[\fint_{t/2}^{t}\fint_{B(\cdot,t)}
 |s^{-\alpha}W^{1/p}(\cdot)\vec F(s,y)|^r\,\mathrm d y\,\mathrm d s
 \right]^{\theta/r}\right)(x).
\end{align}
Take the $L^{q/\theta}(I_{j,i},\mathrm d t/t)$ norm in
\eqref{eq:T-lower-W-point}. Since
\[
 S_{j,1}\cup S_{j,2}=[2^{-j},2^{-j+1}),\qquad
 \int_{I_{j,i}}\frac{\mathrm d t}{t}\asymp1,\qquad
 \sum_{i=1}^2\sum_{j\in\Z}\mathbf1_{I_{j,i}}\le2,
\]
using \eqref{eq:T-dyadic-A} and the continuous version of
Lemma~\ref{lem:FS}, we obtain
\[
 \begin{aligned}
 &\|\vec F\|_{T^{p,q,r}_\alpha(\mathbb A)} \lesssim\left\|\left\{\int_0^\infty \left[\M\left( \left[\fint_{t/2}^t\fint_{B(\cdot,t)} |s^{-\alpha}W^{1/p}(\cdot)\vec F(s,y)|^r\,\mathrm d y\,\mathrm d s \right]^{\theta/r}\right)\right]^{q/\theta} \frac{\mathrm d t}{t}\right\}^{1/q}\right\|_{L^p}\\
 &\quad\lesssim\|\vec F\|_{T^{p,q,r}_\alpha(W)}.
 \end{aligned}
\]
Together with Step 2, this proves \eqref{eq:T-W-A}.

\medskip
\noindent$\mathrm{(ii)}$ The proof is similar to that of (i), and we omit the details.
\end{proof}

For unweighted $Z$-spaces, the following properties follow from
\cite[Proposition 3.9 and Theorem 3.11]{ABH}; for unweighted tent
spaces, see \cite[Section 1, pp.~916--917, and Theorem 5.4]{HuangTent}.

\Needspace{4\baselineskip}
\begin{proposition}\label{prop:Z-basic}
Let $\alpha\in\R$, $0<p<\infty$, $0<q,r\le\infty$, and
let $\mathbb A=\{A_Q\}_{Q\in\D}$ be as in Definition~\ref{def:Z}.
Then the following statements hold.
\begin{enumerate}
\item[$\mathrm{(i)}$] Both $Z^{p,q,r}_\alpha(\mathbb A)$ and
$T^{p,q,r}_\alpha(\mathbb A)$ are quasi-Banach spaces.
They are Banach spaces if $p,q,r\geq1$.
\item[$\mathrm{(ii)}$] If $1<p,q,r<\infty$, both
$Z^{p,q,r}_\alpha(\mathbb A)$ and $T^{p,q,r}_\alpha(\mathbb A)$
are reflexive.
\item[$\mathrm{(iii)}$] If $q,r<\infty$, then
$C_c^\infty(\R^{n+1}_+;\C^m)$ is dense in
$Z^{p,q,r}_\alpha(\mathbb A)$ and $T^{p,q,r}_\alpha(\mathbb A)$.
\item[$\mathrm{(iv)}$] If $W\in\mathcal A_p$, then $\mathrm{(i)}$--$\mathrm{(iii)}$
also hold for $Z^{p,q,r}_\alpha(W)$ and $T^{p,q,r}_\alpha(W)$.
\end{enumerate}
\end{proposition}

\begin{proof}
\noindent$\mathrm{(i)}$ Let $\mathbb I=\{I_m\}_{Q\in\D}$ and define
\[
 (\mathcal U\vec F)(s,y)=A_Q\vec F(s,y),
 \qquad (s,y)\in\widehat Q,\quad Q\in\D.
\]
The map $\mathcal U$ is invertible, with inverse $A_Q^{-1}$ on
$\widehat Q$. By Lemma~\ref{lem:Z-dyadic-A},
\[
 \|\mathcal U\vec F\|_{T^{p,q,r}_\alpha(\mathbb I)}
 \asymp\|\vec F\|_{T^{p,q,r}_\alpha(\mathbb A)},\qquad
 \|\mathcal U\vec F\|_{Z^{p,q,r}_\alpha(\mathbb I)}
 \asymp\|\vec F\|_{Z^{p,q,r}_\alpha(\mathbb A)}.
\]
Thus completeness follows from the cited unweighted results,
applied to each coordinate. For $p,q,r\ge1$, the continuous
expressions in Definitions~\ref{def:Z} and~\ref{def:T} are norms
by Minkowski's inequality.

\medskip
\noindent$\mathrm{(ii)}$ For $1<p,q,r<\infty$, the unweighted spaces
are reflexive by \cite[Theorem 3.11]{ABH} and
\cite[Theorem 5.4]{HuangTent}. The same holds for the matrix spaces
by the isomorphisms in (i).

\medskip
\noindent$\mathrm{(iii)}$ Let $q,r<\infty$. By the cited unweighted
density results and $\mathcal U$, it suffices to approximate
$\vec F$ supported in $\bigcup_{Q\in\mathcal F}\widehat Q$, where
$\mathcal F\subset\D$ is finite and
$A_Q\vec F\in L^r(\widehat Q,\nu_Q;\C^m)$ for $Q\in\mathcal F$.
Choose $\vec\psi_{Q,\ell}\in
C_c^\infty(\operatorname{int}\widehat Q;\C^m)$, extended by zero,
such that
\[
 \|A_Q\vec F-\vec\psi_{Q,\ell}\|_{L^r(\widehat Q,\nu_Q;\C^m)}
 \longrightarrow0\qquad(\ell\to\infty).
\]
Then $\sum_{Q\in\mathcal F}A_Q^{-1}\vec\psi_{Q,\ell}
\in C_c^\infty(\R^{n+1}_+;\C^m)$. Using
Lemma~\ref{lem:Z-dyadic-A}, we obtain
\begin{align*}
 &\left\|\vec F-\sum_{Q\in\mathcal F}A_Q^{-1}\vec\psi_{Q,\ell}\right\|
 _{T^{p,q,r}_\alpha(\mathbb A)}
 +\left\|\vec F-\sum_{Q\in\mathcal F}A_Q^{-1}\vec\psi_{Q,\ell}\right\|
 _{Z^{p,q,r}_\alpha(\mathbb A)}\\
 &\quad\lesssim_{\mathcal F}\max_{Q\in\mathcal F}
 \ell(Q)^{-\alpha}|Q|^{1/p}
 \|A_Q\vec F-\vec\psi_{Q,\ell}\|_{L^r(\widehat Q,\nu_Q;\C^m)}
 \longrightarrow0.
\end{align*}

\medskip
\noindent$\mathrm{(iv)}$ For $W\in\mathcal A_p$, choose a sequence
$\mathbb A=\{A_Q\}_{Q\in\D}$ of reducing operators of order $p$
for $W$. By Lemma~\ref{lem:Z-W-A}, we have
\[
 \|\vec F\|_{T^{p,q,r}_\alpha(W)}
 \asymp\|\vec F\|_{T^{p,q,r}_\alpha(\mathbb A)},\qquad
 \|\vec F\|_{Z^{p,q,r}_\alpha(W)}
 \asymp\|\vec F\|_{Z^{p,q,r}_\alpha(\mathbb A)}.
\]
Hence $\mathrm{(i)}$--$\mathrm{(iii)}$ also hold for
$Z^{p,q,r}_\alpha(W)$ and $T^{p,q,r}_\alpha(W)$.
\end{proof}

\Needspace{9\baselineskip}
\section{Whitney average characterizations of matrix-weighted Besov and Triebel--Lizorkin spaces}\label{sec:main}

In this section, we characterize homogeneous matrix-weighted
Besov and Triebel--Lizorkin spaces via Whitney averages.
The kernels satisfy a Tauberian condition and a finite moment
condition.

\subsection{Kernel conditions and the main theorem}
\label{subsec:Whitney-main}

For kernels with finitely many vanishing moments, we use the
following realization of homogeneous distributions.

\begin{definition}\label{def:moment-spaces}
For $R\in\mathbb N_0$, set
\[
 \cS_R=\left\{\eta\in\cS: \int_{\R^n}x^\beta\eta(x)\,\mathrm d x=0 \text{ for all }\beta\in\mathbb N_0^n\text{ with }|\beta|\le R\right\},\qquad \cS_\infty=\bigcap_{R\ge0}\cS_R.
\]
Let $\cS_R'$ and $\cS_\infty'$ be their continuous duals with respect
to the subspace topology of $\cS$, and let $\cP_R$ be the polynomials
of degree at most $R$.
\end{definition}
We use the standard identifications and
convolution definition from \cite[Section 2]{ABH}:
\[
 \cS_R'\simeq\cS'/\cP_R,\qquad
 \cS_\infty'\simeq\cS'/\cP,\qquad
 (\Phi_s*\vec g)(y):=\vec g\big(\Phi_s(y-\cdot)\big),
\]
where $\Phi\in\cS_R$, $\Phi_s(x)=s^{-n}\Phi(x/s)$,
$\vec g\in(\cS_R')^m$, $s>0$, and $y\in\R^n$.

\begin{assumption}\label{ass:general-Phi}
Let $\alpha\in\R$, $0<p<\infty$, $W\in\mathcal A_p$,
$R\in\mathbb N_0$, and $\Phi\in\cS_R$. Assume that there is
$\varepsilon>0$ such that
\begin{equation}\label{eq:general-Tauberian}
 \widehat\Phi(\xi)\ne0
 \qquad\text{if }\varepsilon/2<|\xi|<\varepsilon.
\end{equation}
In Theorem~\ref{thm:general-Whitney}, we assume
\begin{equation}\label{eq:general-moment-bound}
 \int_{\R^n}x^\beta\Phi(x)\,\mathrm d x=0 \quad\text{for all }\beta\in\mathbb N_0^n\text{ with }|\beta|\le R,\qquad R+1>\alpha+\delta_1=\alpha+d/p.
\end{equation}
Here $d\in[0,n)$ is an $\mathcal A_p$-dimension of $W$,
and $\delta_1,\delta_2,\Delta$ are as in \eqref{eq:dimension-exponents}.
\end{assumption}

\Needspace{8\baselineskip}
We first construct a finite-moment realization of homogeneous distributions
that is common to the Besov and Triebel--Lizorkin spaces.
For the unweighted setting, see \cite[Remark 2.3]{ABH}.

\begin{lemma}\label{lem:general-realization}
Let $\alpha\in\R$, $p\in(0,\infty)$, $q\in(0,\infty]$, and
$W\in\mathcal A_p$. Let $R\in\mathbb N_0$ satisfy
$R+1>\alpha+d/p$, where $d\in[0,n)$ is the $\mathcal A_p$-dimension of
$W$ fixed in \eqref{eq:dimension-exponents}. Fix $\varphi\in\cS$ satisfying
\eqref{eq:partition}, with $\varphi_j(x)=2^{jn}\varphi(2^jx)$
for every $j\in\Z$. Let
$\vec f\in\dot B^\alpha_{p,q}(W)\cup
\dot F^\alpha_{p,q}(W)$.

\smallskip
\noindent$\mathrm{(i)}$ For every $\vec\eta\in(\cS_R)^m$,
\[
 \sum_{j\in\Z}\left|\left\langle\varphi_j*\vec f,\vec\eta\right\rangle\right|<\infty.
\]

\smallskip
\noindent$\mathrm{(ii)}$ Set
\begin{equation}\label{eq:general-realization}
 \mathcal I_R\vec f
       :=\sum_{j\in\Z}\varphi_j*\vec f
       \quad\text{in }(\cS_R')^m.
\end{equation}
The maps
\[
 \mathcal I_R:\dot B^\alpha_{p,q}(W)\to(\cS_R')^m,
 \qquad
 \mathcal I_R:\dot F^\alpha_{p,q}(W)\to(\cS_R')^m
\]
are continuous and linear.
Moreover, $\bigl.(\mathcal I_R\vec f)\bigr|_{(\cS_\infty)^m}
=\vec f\in(\cS_\infty')^m$.

\end{lemma}

% Allow long proof displays to continue across pages.
\allowdisplaybreaks[2]
\displaywidowpenalty=1000

\begin{proof}
\noindent$\mathrm{(i)}$
We give the proof for $\vec f\in\dot F^\alpha_{p,q}(W)$.
The Besov case is similar. Let $d,\delta_1,\delta_2$ and
$\Delta$ be the weight exponents in
\eqref{eq:dimension-exponents}, and let $\mathbb A=\{A_Q\}_{Q\in\D}$ be a sequence of reducing
operators of order $p$ for $W$.

Fix $\varphi\in\cS$ satisfying \eqref{eq:partition}, and set
$\varphi_j(x)=2^{jn}\varphi(2^jx)$ for $j\in\Z$.
Choose $\Psi\in\cS$ with $\widehat\Psi\in C_c^\infty$ such that
\[
 \supp\widehat\Psi\subset\{\xi:1/3<|\xi|<3\},\qquad
 \widehat\Psi(\xi)=1\quad\text{if }1/2\le|\xi|\le2.
\]
Put $\Psi_j(x)=2^{jn}\Psi(2^jx)$. Then
$\varphi_j*\vec f=\Psi_j*\varphi_j*\vec f$.

\medskip
\noindent\emph{Step 1.}
By \eqref{eq:known-equivalence}, for every $j\in\Z$,
\[
 \|\mathbb A_j(\cdot)(\varphi_j*\vec f)(\cdot)\|_{L^p(\R^n)}
 \lesssim2^{-j\alpha}\|\vec f\|_{\dot F^\alpha_{p,q}(W)}.
\]
Fix the reference matrix $\mathbb A_0(0)$. By
\eqref{eq:comparison},
\begin{equation}\label{eq:general-inverse-growth}
 \|\mathbb A_j(z)^{-1}\|\lesssim
 \begin{cases}
  2^{-j\delta_1}(1+2^j|z|)^\Delta,&j\le0,\\
  2^{j\delta_2}(1+|z|)^\Delta,&j\ge0.
 \end{cases}
\end{equation}
Choose $0<\theta<\min\{1,p\}$ and
\[
 L>\max\{\Delta\theta+n,\ n/p+\delta_2-\alpha,\ 0\}.
\]
The Fourier support of $\varphi_j*\vec f$ is contained in
$\{\xi:2^{j-1}\le|\xi|\le2^{j+1}\}$.
Using Lemma~\ref{lem:bandlimited} with $a=L/\theta$, we obtain
\[
 |(\varphi_j*\vec f)(y)|^\theta \lesssim\int_{\R^n} \frac{2^{jn}|(\varphi_j*\vec f)(z)|^\theta} {(1+2^j|y-z|)^L}\,\mathrm d z \le\int_{\R^n} \frac{2^{jn}\|\mathbb A_j(z)^{-1}\|^\theta |\mathbb A_j(z)(\varphi_j*\vec f)(z)|^\theta} {(1+2^j|y-z|)^L}\,\mathrm d z.
\]
Note that
\[
 1+2^j|z|\le(1+2^j|y|)(1+2^j|y-z|),\quad j\in\Z,\qquad 1+|z|\le(1+|y|)(1+2^j|y-z|),\quad j\ge0.
\]
Using H\"older's inequality with $p/\theta>1$ and
$(L-\Delta\theta)p/(p-\theta)>n$, we obtain
\begin{equation}\label{eq:general-holder}
 \begin{aligned}
 \int_{\R^n}\frac{2^{jn}|\mathbb A_j(z)(\varphi_j*\vec f)(z)|^\theta}
 {(1+2^j|y-z|)^{L-\Delta\theta}}\,\mathrm d z
 &\le\|\mathbb A_j(\cdot)(\varphi_j*\vec f)(\cdot)\|_{L^p(\R^n)}^\theta
 \left\|\frac{2^{jn}}
 {(1+2^j|y-\cdot|)^{L-\Delta\theta}}\right\|_{L^{p/(p-\theta)}(\R^n)}\\
 &\lesssim\|\mathbb A_j(\cdot)(\varphi_j*\vec f)(\cdot)\|_{L^p(\R^n)}^\theta
 2^{jn\theta/p}.
 \end{aligned}
\end{equation}
Using \eqref{eq:general-inverse-growth} and
\eqref{eq:general-holder}, we obtain
\begin{equation}\label{eq:general-point-growth}
 |(\varphi_j*\vec f)(y)|
 \lesssim\|\mathbb A_j(\cdot)(\varphi_j*\vec f)(\cdot)\|_{L^p(\R^n)}
 \begin{cases}
  2^{j(n/p-\delta_1)}(1+2^j|y|)^\Delta,&j\le0,\\
  2^{j(n/p+\delta_2)}(1+|y|)^\Delta,&j\ge0.
 \end{cases}
\end{equation}

\medskip
\noindent\emph{Step 2.}
Fix $K>n+\Delta$, keep $L$ from Step 1, and set
$N=\lfloor K\rfloor+1$. For $\vec\eta\in(\cS_R)^m$, by the moment conditions, we have
\[
 \partial^\beta\widehat{\vec\eta}(0)=0\qquad(|\beta|\le R).
\]
For $j\le0$, using Taylor's formula and the moment conditions, we obtain,
for $|\gamma|\le R+1$ and $\xi\in\supp\widehat\Psi$,
\begin{align*}
 2^{j|\gamma|}|\partial^\gamma\widehat{\vec\eta}(2^j\xi)|
 &\lesssim 2^{j|\gamma|}|2^j\xi|^{R+1-|\gamma|}
 \max_{|\mu|=R+1}\|\partial^\mu\widehat{\vec\eta}\|_{L^\infty}\\
 &\lesssim 2^{j(R+1)}
 \rho_{R+N+n+2,\,\lceil L\rceil+n+1}(\vec\eta).
\end{align*}
For $R+1<|\gamma|\le N$, the same bound follows from
$2^{j|\gamma|}\le2^{j(R+1)}$. By Leibniz's formula, we have
\begin{equation}\label{eq:general-test-fourier}
\begin{aligned}
 &\max_{|\beta|\le N}
 \|\partial^\beta[\overline{\widehat\Psi}\,
 \widehat{\vec\eta}(2^j\cdot)]\|_{L^1(\R^n)}\\
 &\quad\lesssim\sum_{|\gamma|\le N}2^{j|\gamma|}
 \int_{\supp\widehat\Psi}|\partial^\gamma\widehat{\vec\eta}(2^j\xi)|\,\mathrm d\xi
 \lesssim 2^{j(R+1)}
 \rho_{R+N+n+2,\,\lceil L\rceil+n+1}(\vec\eta).
\end{aligned}
\end{equation}
For $j\ge0$, using the Schwartz decay, we obtain
\[
 \max_{|\gamma|\le N}|\partial^\gamma\widehat{\vec\eta}(\xi)|
 \lesssim(1+|\xi|)^{-\lceil L\rceil-n-1}
 \rho_{R+N+n+2,\,\lceil L\rceil+n+1}(\vec\eta).
\]
Since
\[
 \supp[\widehat\Psi(2^{-j}\cdot)]
 \subset\{\xi:2^j/3<|\xi|<3\cdot2^j\},
\]
Leibniz's formula implies
\begin{equation}\label{eq:general-test-fourier-high}
\begin{aligned}
 &\max_{|\beta|\le N}
 \|\partial^\beta[\overline{\widehat\Psi}(2^{-j}\cdot)
 \,\widehat{\vec\eta}]\|_{L^1(\R^n)}\\
 &\quad\lesssim\sum_{|\gamma|\le N}2^{-j|\gamma|}
 \int_{2^j/3<|\xi|<3\cdot2^j}
 \max_{|\mu|\le N}|\partial^\mu\widehat{\vec\eta}(\xi)|\,\mathrm d\xi\\
 &\quad\lesssim2^{jn}2^{-j(\lceil L\rceil+n+1)}
 \rho_{R+N+n+2,\,\lceil L\rceil+n+1}(\vec\eta)
 \lesssim2^{-jL}\rho_{R+N+n+2,\,\lceil L\rceil+n+1}(\vec\eta).
\end{aligned}
\end{equation}
Using \eqref{eq:general-test-fourier} and
\eqref{eq:general-test-fourier-high} with $N>K$, Fourier inversion
and integration by parts, we obtain
\begin{equation}\label{eq:general-test-decay}
 |(\overline{\Psi_j(-\cdot)}*\vec\eta)(y)|\lesssim
 \rho_{R+N+n+2,\,\lceil L\rceil+n+1}(\vec\eta)
 \begin{cases}
  2^{j(n+R+1)}(1+2^j|y|)^{-K},&j\le0,\\
  2^{-jL}(1+|y|)^{-K},&j\ge0.
 \end{cases}
\end{equation}

\medskip
\noindent\emph{Step 3.}
Since $\varphi_j*\vec f=\Psi_j*\varphi_j*\vec f$,
\[
 \langle\varphi_j*\vec f,\vec\eta\rangle
 =\langle\varphi_j*\vec f,\overline{\Psi_j(-\cdot)}*\vec\eta\rangle.
\]
Since $K>n+\Delta$,
\begin{equation}\label{eq:general-spatial-integral}
 \int_{\R^n}(1+2^j|y|)^{\Delta-K}\,\mathrm d y
 =2^{-jn}\int_{\R^n}(1+|y|)^{\Delta-K}\,\mathrm d y
 \lesssim2^{-jn},\qquad j\in\Z.
\end{equation}
Using \eqref{eq:general-point-growth}, \eqref{eq:general-test-decay}
and \eqref{eq:general-spatial-integral}, we obtain
\begin{equation}\label{eq:general-realization-pairing}
 \left|\langle\varphi_j*\vec f,\vec\eta\rangle\right|
 \lesssim\|\vec f\|_{\dot F^\alpha_{p,q}(W)}\rho_{R+N+n+2,\,\lceil L\rceil+n+1}(\vec\eta)
 \begin{cases}
  2^{j[R+1-\alpha+(n-d)/p]},&j\le0,\\
  2^{-j[L+\alpha-n/p-\delta_2]},&j\ge0.
 \end{cases}
\end{equation}
Using \eqref{eq:general-realization-pairing} with
$R+1-\alpha+(n-d)/p>0$ and $L+\alpha-n/p-\delta_2>0$, we obtain
\begin{align}\label{eq:general-realization-continuity}
 \sum_{j\in\Z}|\langle\varphi_j*\vec f,\vec\eta\rangle|
 &\lesssim\|\vec f\|_{\dot F^\alpha_{p,q}(W)}
 \rho_{R+N+n+2,\,\lceil L\rceil+n+1}(\vec\eta) \times\left[\sum_{j\le0}2^{j[R+1-\alpha+(n-d)/p]}
 +\sum_{j\ge1}2^{-j[L+\alpha-n/p-\delta_2]}\right]\notag\\
 &\lesssim\|\vec f\|_{\dot F^\alpha_{p,q}(W)}
 \rho_{R+N+n+2,\,\lceil L\rceil+n+1}(\vec\eta)<\infty.
 \end{align}

\medskip
\noindent$\mathrm{(ii)}$
Using \eqref{eq:general-realization-continuity},
\eqref{eq:general-realization} defines
continuous linear maps into $(\cS_R')^m$. By
\eqref{eq:partition}, we have
$\bigl.(\mathcal I_R\vec f)\bigr|_{(\cS_\infty)^m}
=\vec f\in(\cS_\infty')^m$.\qedhere
\end{proof}

\Needspace{10\baselineskip}
\begin{theorem}
\label{thm:general-Whitney}
Let $\alpha\in\R$, $0<p<\infty$, $0<q,r\le\infty$,
$W\in\mathcal A_p$, and let $\mathbb A$ be as in
Definition~\ref{def:reducing}. Let $R\in\mathbb N_0$ and
$\Phi\in\cS_R$ satisfy Assumption~\ref{ass:general-Phi},
set $\Phi_s(x)=s^{-n}\Phi(x/s)$ for $s>0$, and use
$\mathcal I_R$ from \eqref{eq:general-realization}.

\smallskip
\noindent$\mathrm{(i)}$ For $\vec f\in\dot F^\alpha_{p,q}(W)$,
\begin{equation}\label{eq:general-Whitney-F}
 \begin{aligned}
 &\|\vec f\|_{\dot F^\alpha_{p,q}(W)} \asymp\left\|\left\{\int_0^\infty \left[\fint_{t/2}^t\fint_{B(\cdot,t)} |s^{-\alpha}\mathbb A(t,\cdot)(\Phi_s*\mathcal I_R\vec f)(y)|^r\,\mathrm d y\,\mathrm d s \right]^{q/r}\frac{\mathrm d t}{t}\right\}^{1/q}\right\|_{L^p(\R^n)}\\
 &\quad \asymp\left\|\left\{\int_0^\infty \left[\fint_{t/2}^t\fint_{B(\cdot,t)} |s^{-\alpha}W^{1/p}(\cdot)(\Phi_s*\mathcal I_R\vec f)(y)|^r\,\mathrm d y\,\mathrm d s \right]^{q/r}\frac{\mathrm d t}{t}\right\}^{1/q}\right\|_{L^p(\R^n)}.
 \end{aligned}
\end{equation}

\smallskip
\noindent$\mathrm{(ii)}$ For $\vec f\in\dot B^\alpha_{p,q}(W)$,
\begin{equation}\label{eq:general-Whitney-B}
 \begin{aligned}
 &\|\vec f\|_{\dot B^\alpha_{p,q}(W)} \asymp\left\{\int_0^\infty \left\|\left[\fint_{t/2}^t\fint_{B(\cdot,t)} |s^{-\alpha}\mathbb A(t,\cdot)(\Phi_s*\mathcal I_R\vec f)(y)|^r\,\mathrm d y\,\mathrm d s \right]^{1/r}\right\|_{L^p(\R^n)}^q\frac{\mathrm d t}{t}\right\}^{1/q}\\
 &\quad \asymp\left\{\int_0^\infty \left\|\left[\fint_{t/2}^t\fint_{B(\cdot,t)} |s^{-\alpha}W^{1/p}(\cdot)(\Phi_s*\mathcal I_R\vec f)(y)|^r\,\mathrm d y\,\mathrm d s \right]^{1/r}\right\|_{L^p(\R^n)}^q\frac{\mathrm d t}{t}\right\}^{1/q}.
 \end{aligned}
\end{equation}
\end{theorem}

\Needspace{5\baselineskip}
\begin{remark}\phantomsection\label{rem:Z-extension}
\noindent$\mathrm{(i)}$ Under the assumptions of Theorem~\ref{thm:general-Whitney},
suppose further that $\widehat\Phi\in C_c^\infty(\R^n\setminus\{0\})$.
For every $\vec f\in(\cS'/\cP)^m$,
\begin{align*}
 &\|\vec f\|_{\dot B^\alpha_{p,q}(W)} \asymp\bigl\|(s,y)\mapsto(\Phi_s*\vec f)(y)\bigr\|_{Z^{p,q,r}_\alpha(W)} \asymp\bigl\|(s,y)\mapsto(\Phi_s*\vec f)(y)\bigr\|_{Z^{p,q,r}_\alpha(\mathbb A)},\\
 &\|\vec f\|_{\dot F^\alpha_{p,q}(W)} \asymp\bigl\|(s,y)\mapsto(\Phi_s*\vec f)(y)\bigr\|_{T^{p,q,r}_\alpha(W)} \asymp\bigl\|(s,y)\mapsto(\Phi_s*\vec f)(y)\bigr\|_{T^{p,q,r}_\alpha(\mathbb A)}.
 \end{align*}
Indeed, all moments of $\Phi$ vanish, and $R$ may be chosen
to satisfy \eqref{eq:general-moment-bound}.
For $\vec f\in\dot B^\alpha_{p,q}(W)\cup\dot F^\alpha_{p,q}(W)$,
the realization \eqref{eq:general-realization} satisfies
\[
 \Phi_s*\mathcal I_R\vec f=\Phi_s*\vec f.
\]
Thus \eqref{eq:general-Whitney-F} and \eqref{eq:general-Whitney-B}
hold with $\Phi_s*\vec f$ in place of $\Phi_s*\mathcal I_R\vec f$.
By Proposition~\ref{prop:general-reverse} and
\eqref{eq:known-equivalence}, these equivalences extend to every
$\vec f\in(\cS'/\cP)^m$.

\smallskip
\noindent$\mathrm{(ii)}$ The whole $Z$ and $T$ spaces allow arbitrary measurable
functions of $(s,y)$; the extension characterizations concern
only the functions generated by convolution.

\smallskip
\Needspace{9\baselineskip}
\noindent$\mathrm{(iii)}$ When $m=1$ and $W\equiv1$, we may take
$A_Q=1$ for every $Q\in\D$, so that
\[
 \mathbb A_j(x)=\mathbb A(t,x)=1,\qquad
 d=\widetilde d=\delta_1=\delta_2=\Delta=0.
\]
Thus \eqref{eq:general-moment-bound} becomes $R+1>\alpha$.
On the common range $0<p<\infty$, $0<q,r\le\infty$, and
$R\in\mathbb N_0$, the Besov characterization
\eqref{eq:general-Whitney-B}, with the realization
\eqref{eq:general-realization} in $\cS_R'$, agrees with the
scalar unweighted result of \cite[Theorem 2.2]{ABH}.
The Triebel--Lizorkin characterization
\eqref{eq:general-Whitney-F} likewise agrees with
\cite[Theorem 2.4]{ABH}.

\smallskip
\noindent$\mathrm{(iv)}$ Discrete Lusin-average characterizations with matrix weights were
established in \cite[Definition 4.6 and Theorem 4.7, pp.~52--53]{BYYZ}.
Theorem~\ref{thm:general-Whitney} uses continuous scale integration
and Whitney averages over $(t/2,t)\times B(x,t)$.
\end{remark}

\subsection{Proof of Theorem~\ref{thm:general-Whitney}}
\label{subsec:Whitney-proofs}

\Needspace{4\baselineskip}
The following estimates follow from
\cite[Lemma A.3(i)]{Ullrich} by dilation, with the roles of the kernels
interchanged in the second case.

\begin{lemma}\label{lem:general-kernel}
Let $R\in\mathbb N_0$, $\Phi\in\cS_R$, and let
$\Psi\in\cS$ have Fourier support in
$\{\xi\in\R^n:1/3<|\xi|<3\}$. Set
$\Phi_s(x)=s^{-n}\Phi(x/s)$ for $s>0$ and
$\Psi_k(x)=2^{kn}\Psi(2^kx)$ for $k\in\Z$.
For every $L,H>0$
there is $C_{L,H}>0$ such that, for all $j,k\in\Z$,
$v\in\R^n$, and $s\in(2^{-j-1},2^{-j+1})$,
\begin{equation}\label{eq:general-kernel-decay}
 |(\Phi_s*\Psi_k)(v)|\le C_{L,H}2^{kn}(1+2^k|v|)^{-L}
 \begin{cases}
  2^{-(R+1)(j-k)},&k\le j,\\
  2^{-H(k-j)},&k>j.
 \end{cases}
\end{equation}
\end{lemma}

Here we use $(\Phi_s*\Psi_k)(v)=2^{kn}(\Phi_{s2^k}*\Psi)(2^kv)$
with $s2^k\asymp2^{k-j}$.

\Needspace{4\baselineskip}
The following Calder\'on reproducing formula uses the dilation ratio
$\sqrt2$; see \cite[(2.39)--(2.47)]{Ullrich}.

\begin{lemma}\label{lem:general-reproduction}
Let $\Phi\in\cS$ satisfy \eqref{eq:general-Tauberian} for
some $\varepsilon>0$. There exists $\Lambda\in\cS$ such that
$\widehat\Lambda$ is supported in a compact subset of
$\{\xi\in\R^n:\varepsilon/2<|\xi|<\varepsilon\}$ and
\begin{equation}\label{eq:general-annular-partition}
 \sum_{l\in\Z}\widehat\Lambda(2^{-l/2}\xi)
                  \widehat\Phi(2^{-l/2}\xi)=1,
 \qquad \xi\in\R^n\setminus\{0\}.
\end{equation}
\end{lemma}

For the construction, choose a nonnegative radial $\eta\in C_c^\infty(\R^n)$ supported
in $\{\xi:\varepsilon/2<|\xi|<\varepsilon\}$ such that
$\sum_{l\in\Z}\eta(2^{-l/2}\xi)>0$ for $\xi\ne0$. Set
\[
 \widehat\Lambda(\xi):=
 \begin{cases}
  \displaystyle\frac{\eta(\xi)}
   {\widehat\Phi(\xi)\sum_{v\in\Z}\eta(2^{-v/2}\xi)},
       &\varepsilon/2<|\xi|<\varepsilon,\\[4pt]
  0,&\text{otherwise}.
 \end{cases}
\]
Then $\widehat\Lambda\in C_c^\infty(\R^n)$, and
\eqref{eq:general-annular-partition} follows by summation.

\Needspace{5\baselineskip}
The following estimate follows from \cite[(2.66)]{Ullrich}
and its homogeneous version in \cite[proof of Theorem 2.8]{Ullrich},
using the dilation ratio $\sqrt2$ and applying the scalar estimate
to each coordinate of $T\vec g$. The kernel estimates are uniform
for $u\in[1,\sqrt2]$.

\begin{lemma}\label{lem:general-submean}
Let $\Phi\in\cS$ satisfy \eqref{eq:general-Tauberian}, with
$\Phi_s(x)=s^{-n}\Phi(x/s)$. Let $0<\theta\le1$,
$a>0$, and $N\ge a$. Let $\vec g\in(\cS')^m$ be a tempered
distribution and $T\in\C^{m\times m}$ be a constant matrix.
For all $k\in\Z$, $x\in\R^n$, and $1\le u\le\sqrt2$, we have
\begin{equation}\label{eq:general-submean}
 \left[\sup_{y\in\R^n} \frac{|T(\Phi_{2^{-k/2}u}*\vec g)(y)|} {(1+2^{k/2}|x-y|)^a}\right]^\theta \lesssim\sum_{h\ge0}2^{-hN\theta/2}2^{(k+h)n/2} \int_{\R^n}\frac{|T(\Phi_{2^{-(k+h)/2}u}*\vec g)(z)|^\theta} {(1+2^{k/2}|x-z|)^{a\theta}}\,\mathrm d z.
\end{equation}
\end{lemma}

\Needspace{5\baselineskip}
We have the following lemma by adapting the averaging argument of
\cite[proof of Theorem 2.2]{ABH}, using Lemma~\ref{lem:general-submean}
and the matrix estimates in Section~\ref{sec:preliminaries}.

\begin{lemma}
\label{lem:general-averages}
Let $\alpha\in\R$, $0<p<\infty$, $0<q,r\le\infty$,
$W\in\mathcal A_p$, $R\in\mathbb N_0$, and
$\vec g\in(\cS_R')^m$. Let $\mathbb A=\{A_Q\}_{Q\in\D}$ be a sequence of reducing
operators of order $p$ for $W$,
and let $\Phi\in\cS_R$ satisfy \eqref{eq:general-Tauberian},
with $\Phi_s(x)=s^{-n}\Phi(x/s)$. Use $\delta_2,\Delta$ from
\eqref{eq:dimension-exponents}. Choose
\begin{equation}\label{eq:general-parameters}
 0<\theta<\min\{1,p,q,r\},\qquad a>\Delta+n/\theta,\qquad N\ge a,\qquad (N+\alpha-\delta_2)\theta-n>0.
\end{equation}
For every $k\in\Z$ and almost every $x\in\R^n$, the following estimates hold.

\smallskip
\noindent$\mathrm{(i)}$ For the reducing operators,
\begin{align}
 &2^{k\alpha\theta/2}\int_1^{\sqrt2}\left[\sup_{y\in\R^n}\frac{|\mathbb A(2^{-k/2},x)(\Phi_{2^{-k/2}u}*\vec g)(y)|}{(1+2^{k/2}|x-y|)^a}\right]^{\theta} \frac{\mathrm d u}{u\log_e\sqrt2} \lesssim\sum_{h\ge0}2^{-h[(N+\alpha-\delta_2)\theta-n]/2}\label{eq:general-fine-average}\\
 &\quad \times\M\!\left( 2^{(k+h)\alpha\theta/2}\left[\int_1^{\sqrt2}\fint_{B(\cdot,2^{-(k+h)/2})} |\mathbb A(2^{-(k+h)/2},\cdot)(\Phi_{2^{-(k+h)/2}u}*\vec g)(y)|^r\,\mathrm d y\, \frac{\mathrm d u}{u\log_e\sqrt2}\right]^{\theta/r} \right)(x).\notag
 \end{align}

\smallskip
\noindent$\mathrm{(ii)}$ For the matrix weight,
\begin{align}
 &2^{k\alpha\theta/2}\int_1^{\sqrt2}\left[\sup_{y\in\R^n}\frac{|\mathbb A(2^{-k/2},x)(\Phi_{2^{-k/2}u}*\vec g)(y)|}{(1+2^{k/2}|x-y|)^a}\right]^{\theta} \frac{\mathrm d u}{u\log_e\sqrt2} \lesssim\sum_{h\ge0}2^{-h[(N+\alpha-\delta_2)\theta-n]/2}\label{eq:general-fine-average-W}\\
 &\quad \times\M\!\left( 2^{(k+h)\alpha\theta/2}\left[\int_1^{\sqrt2}\fint_{B(\cdot,2^{-(k+h)/2})} |W^{1/p}(\cdot)(\Phi_{2^{-(k+h)/2}u}*\vec g)(y)|^r\,\mathrm d y\, \frac{\mathrm d u}{u\log_e\sqrt2}\right]^{\theta/r} \right)(x).\notag
 \end{align}
\end{lemma}

\begin{proof}
\noindent$\mathrm{(i)}$ Let $\mathbb A=\{A_Q\}_{Q\in\D}$ be a sequence of reducing
operators of order $p$ for $W$, and let
$\delta_2,\Delta$ be as in \eqref{eq:dimension-exponents}.
Fix $k\in\Z$, $x\in\R^n$, $h\in\mathbb N_0$ and
$u\in[1,\sqrt2]$. The parameters $\theta,a,N$ satisfy
\eqref{eq:general-parameters}.

\medskip
\noindent\emph{Step 1.}
Choose $\widetilde{\vec g}\in(\cS')^m$ such that
$\bigl.\widetilde{\vec g}\bigr|_{(\cS_R)^m}
=\vec g\in(\cS_R')^m$.
Since $\Phi\in\cS_R$, we have
$\Phi_s*\widetilde{\vec g}=\Phi_s*\vec g$ for $s>0$.
Since $2^{-(k+h)/2}\le2^{-k/2}$, if
$|y-z|<2^{-(k+h)/2}$, then
\[
 \frac{(1+2^{k/2}|x-y|)^{a\theta}}
      {(1+2^{k/2}|x-z|)^{a\theta}}\asymp1.
\]
Inserting the normalized average over $z\in B(y,2^{-(k+h)/2})$
and using Tonelli's theorem and the preceding comparison, we obtain
\[
 \int_{\R^n} \frac{|\mathbb A(2^{-k/2},x) (\Phi_{2^{-(k+h)/2}u}*\vec g)(y)|^\theta} {(1+2^{k/2}|x-y|)^{a\theta}}\,\mathrm d y \lesssim\int_{\R^n} \frac{\fint_{B(z,2^{-(k+h)/2})} |\mathbb A(2^{-k/2},x) (\Phi_{2^{-(k+h)/2}u}*\vec g)(y)|^\theta\,\mathrm d y} {(1+2^{k/2}|x-z|)^{a\theta}}\,\mathrm d z.
\]
For $k\in\Z$ and $x,z\in\R^n$, by \eqref{eq:At-definition}, we have
\[
 \mathbb A(2^{-k/2},x)=\mathbb A_{\lceil k/2\rceil}(x).
\]
Hence \eqref{eq:comparison} yields
\begin{equation}\label{eq:general-fine-comparison}
 \|\mathbb A(2^{-k/2},x)\mathbb A(2^{-(k+h)/2},z)^{-1}\|
 \lesssim2^{h\delta_2/2}(1+2^{k/2}|x-z|)^\Delta,
 \qquad h\in\mathbb N_0.
\end{equation}
Using $\mathbb A(2^{-k/2},x)
=\mathbb A(2^{-k/2},x)\mathbb A(2^{-(k+h)/2},z)^{-1}
\mathbb A(2^{-(k+h)/2},z)$ and
\eqref{eq:general-fine-comparison}, we obtain
\begin{equation}\label{eq:general-transport-integral}
 \begin{aligned}
 &\int_{\R^n}
 \frac{|\mathbb A(2^{-k/2},x)
          (\Phi_{2^{-(k+h)/2}u}*\vec g)(y)|^\theta}
      {(1+2^{k/2}|x-y|)^{a\theta}}\,\mathrm d y\\
 &\quad\lesssim2^{h\delta_2\theta/2}
 \int_{\R^n}
 \frac{\fint_{B(z,2^{-(k+h)/2})}
       |\mathbb A(2^{-(k+h)/2},z)
          (\Phi_{2^{-(k+h)/2}u}*\vec g)(y)|^\theta\,\mathrm d y}
      {(1+2^{k/2}|x-z|)^{(a-\Delta)\theta}}\,\mathrm d z.
 \end{aligned}
\end{equation}

\medskip
\noindent\emph{Step 2.}
Since $\int_1^{\sqrt2}\mathrm d u/(u\log_e\sqrt2)=1$,
using H\"older's inequality with $\theta<r$, we obtain
\begin{equation}\label{eq:general-average-holder}
 \begin{aligned}
 &\int_1^{\sqrt2}\fint_{B(z,2^{-(k+h)/2})}
 |\mathbb A(2^{-(k+h)/2},z)(\Phi_{2^{-(k+h)/2}u}*\vec g)(y)|^\theta\,\mathrm d y\,
                  \frac{\mathrm d u}{u\log_e\sqrt2}\\
 &\qquad\le \left[\int_1^{\sqrt2}\fint_{B(z,2^{-(k+h)/2})}
       |\mathbb A(2^{-(k+h)/2},z)(\Phi_{2^{-(k+h)/2}u}*\vec g)(y)|^r\,\mathrm d y\,
       \frac{\mathrm d u}{u\log_e\sqrt2}\right]^{\theta/r}.
 \end{aligned}
\end{equation}
Note that
\[
 2^{k\alpha\theta/2}2^{-hN\theta/2}2^{(k+h)n/2}2^{h\delta_2\theta/2}
 =2^{-h[(N+\alpha-\delta_2)\theta-n]/2}
  2^{kn/2}2^{(k+h)\alpha\theta/2}.
\]
Using Lemma~\ref{lem:general-submean} for $\widetilde{\vec g}$ with
$T=\mathbb A(2^{-k/2},x)$, together with
\eqref{eq:general-transport-integral} and
\eqref{eq:general-average-holder}, we obtain
\begin{align*}
 &2^{k\alpha\theta/2}\int_1^{\sqrt2}\left[\sup_{y\in\R^n}\frac{|\mathbb A(2^{-k/2},x)(\Phi_{2^{-k/2}u}*\vec g)(y)|}{(1+2^{k/2}|x-y|)^a}\right]^{\theta}
                  \frac{\mathrm d u}{u\log_e\sqrt2}\\
 &\quad\lesssim\sum_{h\ge0}2^{-h[(N+\alpha-\delta_2)\theta-n]/2}
 \int_{\R^n}
 \frac{2^{kn/2}2^{(k+h)\alpha\theta/2}}
      {(1+2^{k/2}|x-z|)^{(a-\Delta)\theta}}\\
 &\qquad\qquad\times\left[\int_1^{\sqrt2}\fint_{B(z,2^{-(k+h)/2})}
       |\mathbb A(2^{-(k+h)/2},z)(\Phi_{2^{-(k+h)/2}u}*\vec g)(y)|^r\,\mathrm d y\,
       \frac{\mathrm d u}{u\log_e\sqrt2}\right]^{\theta/r}\,\mathrm d z.
\end{align*}
Since $(a-\Delta)\theta>n$, \eqref{eq:decay-max} at scale
$2^{-k/2}$ proves \eqref{eq:general-fine-average}.

\medskip
\noindent$\mathrm{(ii)}$ Choose an integer $L_0\ge0$ with
$\sqrt n\,2^{-L_0}<1/2$ and partition $\R^n$ into cubes
$Q\in\D_{\lceil(k+h)/2\rceil+L_0}$. With $C_1$ from \eqref{eq:inverse}, set
\[
 E_Q=\{z\in Q:\|A_QW^{-1/p}(z)\|\le2C_1\},\qquad
 |E_Q|\ge|Q|/2.
\]
For $y\in Q$
and $z\in E_Q$, one has $|y-z|<2^{-(k+h)/2}/2$ and
\[
 \mathbb A(2^{-k/2},x)=
 \mathbb A(2^{-k/2},x)A_Q^{-1}
 A_QW^{-1/p}(z)W^{1/p}(z).
\]
Thus \eqref{eq:comparison} and the definition of $E_Q$ imply
\begin{equation}\label{eq:general-good-set-transport}
 |\mathbb A(2^{-k/2},x) (\Phi_{2^{-(k+h)/2}u}*\vec g)(y)|^\theta \lesssim2^{h\delta_2\theta/2} (1+2^{k/2}|x-z|)^{\Delta\theta} |W^{1/p}(z)(\Phi_{2^{-(k+h)/2}u}*\vec g)(y)|^\theta.
\end{equation}
For $y\in Q$ and $z\in E_Q$,
\[
 1+2^{k/2}|x-y|\asymp1+2^{k/2}|x-z|,\qquad Q\subset B(z,2^{-(k+h)/2}),\qquad |B(z,2^{-(k+h)/2})|\lesssim |E_Q|.
\]
Insert the normalized average over $E_Q$, apply
\eqref{eq:general-good-set-transport}, and sum over the disjoint sets
$E_Q$. We obtain
\begin{equation}\label{eq:general-transport-integral-W}
 \begin{aligned}
 &\int_{\R^n}
 \frac{|\mathbb A(2^{-k/2},x)
          (\Phi_{2^{-(k+h)/2}u}*\vec g)(y)|^\theta}
      {(1+2^{k/2}|x-y|)^{a\theta}}\,\mathrm d y\\
 &\quad=\sum_{Q\in\D_{\lceil(k+h)/2\rceil+L_0}}
 \frac1{|E_Q|}\int_{E_Q}\int_Q
 \frac{|\mathbb A(2^{-k/2},x)(\Phi_{2^{-(k+h)/2}u}*\vec g)(y)|^\theta}
 {(1+2^{k/2}|x-y|)^{a\theta}}\,\mathrm d y\,\mathrm d z\\
 &\quad\lesssim2^{h\delta_2\theta/2}
 \int_{\R^n}
 \frac{\fint_{B(z,2^{-(k+h)/2})}
       |W^{1/p}(z)(\Phi_{2^{-(k+h)/2}u}*\vec g)(y)|^\theta\,\mathrm d y}
      {(1+2^{k/2}|x-z|)^{(a-\Delta)\theta}}\,\mathrm d z.
 \end{aligned}
\end{equation}
Using H\"older's inequality with $\theta<r$, we obtain
\[
 \begin{aligned}
 &\int_1^{\sqrt2}\fint_{B(z,2^{-(k+h)/2})}
 |W^{1/p}(z)(\Phi_{2^{-(k+h)/2}u}*\vec g)(y)|^\theta\,\mathrm d y\,
 \frac{\mathrm d u}{u\log_e\sqrt2}\\
 &\qquad\le\left[\int_1^{\sqrt2}\fint_{B(z,2^{-(k+h)/2})}
 |W^{1/p}(z)(\Phi_{2^{-(k+h)/2}u}*\vec g)(y)|^r\,\mathrm d y\,
 \frac{\mathrm d u}{u\log_e\sqrt2}\right]^{\theta/r}.
 \end{aligned}
\]
Note that
\[
 2^{k\alpha\theta/2}2^{-hN\theta/2}2^{(k+h)n/2}2^{h\delta_2\theta/2}
 =2^{-h[(N+\alpha-\delta_2)\theta-n]/2}
  2^{kn/2}2^{(k+h)\alpha\theta/2}.
\]
Using \eqref{eq:general-submean}, \eqref{eq:general-transport-integral-W},
and \eqref{eq:decay-max} with
$(a-\Delta)\theta>n$, we obtain
\begin{align*}
 &2^{k\alpha\theta/2}\int_1^{\sqrt2}
 \left[\sup_{y\in\R^n}\frac{|\mathbb A(2^{-k/2},x)(\Phi_{2^{-k/2}u}*\vec g)(y)|}
 {(1+2^{k/2}|x-y|)^a}\right]^\theta\frac{\mathrm d u}{u\log_e\sqrt2}\\
 &\quad\lesssim\sum_{h\ge0}2^{-h[(N+\alpha-\delta_2)\theta-n]/2}
 \int_{\R^n}\frac{2^{kn/2}2^{(k+h)\alpha\theta/2}}
 {(1+2^{k/2}|x-z|)^{(a-\Delta)\theta}}\\
 &\qquad\qquad\times\left[\int_1^{\sqrt2}\fint_{B(z,2^{-(k+h)/2})}
 |W^{1/p}(z)(\Phi_{2^{-(k+h)/2}u}*\vec g)(y)|^r\,\mathrm d y\,
 \frac{\mathrm d u}{u\log_e\sqrt2}\right]^{\theta/r}\,\mathrm d z
 \\
 &\quad\lesssim\sum_{h\ge0}2^{-h[(N+\alpha-\delta_2)\theta-n]/2}\\
 &\qquad\quad\times\M\!\left(2^{(k+h)\alpha\theta/2}
 \left[\int_1^{\sqrt2}\fint_{B(\cdot,2^{-(k+h)/2})}
 |W^{1/p}(\cdot)(\Phi_{2^{-(k+h)/2}u}*\vec g)(y)|^r\,\mathrm d y\,
 \frac{\mathrm d u}{u\log_e\sqrt2}\right]^{\theta/r}\right)(x).
 \qedhere
\end{align*}
\end{proof}

\Needspace{5\baselineskip}
We next estimate averages at scales $2^{-k/2}$ in terms of Whitney averages.

\begin{lemma}
\label{lem:general-windows}
Let $\alpha\in\R$, $0<p<\infty$, $0<q,r\le\infty$,
$W\in\mathcal A_p$, $R\in\mathbb N_0$, and
$\vec g\in(\cS_R')^m$. Let $\mathbb A=\{A_Q\}_{Q\in\D}$ be a sequence of reducing
operators of order $p$ for $W$,
and let $\Phi\in\cS_R$, with $\Phi_s(x)=s^{-n}\Phi(x/s)$.
Then the following statements hold.

\smallskip
\noindent$\mathrm{(i)}$ For the reducing operators, we have
\begin{align}
 &\left\|\left\{\sum_{k\in\Z}2^{k\alpha q/2}\left[\int_1^{\sqrt2}\fint_{B(\cdot,2^{-k/2})}
       |\mathbb A(2^{-k/2},\cdot)(\Phi_{2^{-k/2}u}*\vec g)(y)|^r\,\mathrm d y\,
       \frac{\mathrm d u}{u\log_e\sqrt2}\right]^{q/r}\right\}^{1/q}
                           \right\|_{L^p(\R^n)}\notag\\
 &\qquad\lesssim\left\|\left\{\int_0^\infty
 \left[\fint_{t/2}^t\fint_{B(\cdot,t)}
       |s^{-\alpha}\mathbb A(t,\cdot)(\Phi_s*\vec g)(y)|^r\,\mathrm d y\,\mathrm d s
 \right]^{q/r}\frac{\mathrm d t}{t}\right\}^{1/q}\right\|_{L^p(\R^n)},
 \label{eq:general-window-F}
\end{align}

\begin{align}
 &\left\{\sum_{k\in\Z}2^{k\alpha q/2}\left\|\left[\int_1^{\sqrt2}\fint_{B(\cdot,2^{-k/2})}
       |\mathbb A(2^{-k/2},\cdot)(\Phi_{2^{-k/2}u}*\vec g)(y)|^r\,\mathrm d y\,
       \frac{\mathrm d u}{u\log_e\sqrt2}\right]^{1/r}\right\|_{L^p(\R^n)}^q
                           \right\}^{1/q}\notag\\
 &\qquad\lesssim\left\{\int_0^\infty
 \left\|\left[\fint_{t/2}^t\fint_{B(\cdot,t)}
       |s^{-\alpha}\mathbb A(t,\cdot)(\Phi_s*\vec g)(y)|^r\,\mathrm d y\,\mathrm d s
 \right]^{1/r}\right\|_{L^p(\R^n)}^q\frac{\mathrm d t}{t}\right\}^{1/q}.
 \label{eq:general-window-B}
\end{align}
\smallskip
\noindent$\mathrm{(ii)}$ For the matrix weight, we have
\begin{align}
 &\left\|\left\{\sum_{k\in\Z}2^{k\alpha q/2}\left[\int_1^{\sqrt2}\fint_{B(\cdot,2^{-k/2})}
       |W^{1/p}(\cdot)(\Phi_{2^{-k/2}u}*\vec g)(y)|^r\,\mathrm d y\,
       \frac{\mathrm d u}{u\log_e\sqrt2}\right]^{q/r}\right\}^{1/q}
                           \right\|_{L^p(\R^n)}\notag\\
 &\qquad\lesssim\left\|\left\{\int_0^\infty
 \left[\fint_{t/2}^t\fint_{B(\cdot,t)}
       |s^{-\alpha}W^{1/p}(\cdot)(\Phi_s*\vec g)(y)|^r\,\mathrm d y\,\mathrm d s
 \right]^{q/r}\frac{\mathrm d t}{t}\right\}^{1/q}\right\|_{L^p(\R^n)},
 \label{eq:general-window-F-W}
\end{align}

\begin{align}
 &\left\{\sum_{k\in\Z}2^{k\alpha q/2}\left\|\left[\int_1^{\sqrt2}\fint_{B(\cdot,2^{-k/2})}
       |W^{1/p}(\cdot)(\Phi_{2^{-k/2}u}*\vec g)(y)|^r\,\mathrm d y\,
       \frac{\mathrm d u}{u\log_e\sqrt2}\right]^{1/r}\right\|_{L^p(\R^n)}^q
                           \right\}^{1/q}\notag\\
 &\qquad\lesssim\left\{\int_0^\infty
 \left\|\left[\fint_{t/2}^t\fint_{B(\cdot,t)}
       |s^{-\alpha}W^{1/p}(\cdot)(\Phi_s*\vec g)(y)|^r\,\mathrm d y\,\mathrm d s
 \right]^{1/r}\right\|_{L^p(\R^n)}^q\frac{\mathrm d t}{t}\right\}^{1/q}.
 \label{eq:general-window-B-W}
\end{align}
\end{lemma}

\begin{proof}
\noindent$\mathrm{(i)}$ Use $\Phi_s(x)=s^{-n}\Phi(x/s)$ and the reducing family
$\mathbb A$ of order $p$ for $W$.
Fix $\sqrt2<c_1<c_2<2$. For $k\in\Z$ and $x\in\R^n$, whenever
$c_1\,2^{-k/2}<t<c_2\,2^{-k/2}$,
\[
 (2^{-k/2},2^{(1-k)/2})\subset(t/2,t),\qquad
 B(x,2^{-k/2})\subset B(x,t).
\]
For $2^{-k/2}<s<2^{(1-k)/2}$, we have
\[
 \frac{|B(x,t)|}{|B(x,2^{-k/2})|}
 =(2^{k/2}t)^n\asymp1,\qquad
 \frac{c_1}{\sqrt2}<\frac{t}{s}<c_2,\qquad
 s^{-\alpha}\asymp2^{k\alpha/2}.
\]
Moreover,
$\|\mathbb A(2^{-k/2},x)\mathbb A(t,x)^{-1}\|\lesssim1$
by \eqref{eq:comparison}. Consequently
\begin{equation}\label{eq:general-window-pointwise}
 \begin{aligned}
 &2^{k\alpha/2}\left[\int_1^{\sqrt2}\fint_{B(x,2^{-k/2})}
       |\mathbb A(2^{-k/2},x)(\Phi_{2^{-k/2}u}*\vec g)(y)|^r\,\mathrm d y\,
       \frac{\mathrm d u}{u\log_e\sqrt2}\right]^{1/r}\\
 &\qquad\lesssim\left[\fint_{t/2}^t\fint_{B(x,t)}
       |s^{-\alpha}\mathbb A(t,x)(\Phi_s*\vec g)(y)|^r\,\mathrm d y\,\mathrm d s
       \right]^{1/r},\qquad
 c_1\,2^{-k/2}<t<c_2\,2^{-k/2}.
 \end{aligned}
\end{equation}
The intervals $(c_1\,2^{-k/2},c_2\,2^{-k/2})$ have bounded
overlap and satisfy
\[
 \int_{c_1\,2^{-k/2}}^{c_2\,2^{-k/2}}\frac{\mathrm d t}{t}
 =\log_e(c_2/c_1)>0.
\]
Integrating the $q$th power of
\eqref{eq:general-window-pointwise} and summing over $k$ proves
\eqref{eq:general-window-F} pointwise before taking $L^p$.
Taking $L^p$ first in \eqref{eq:general-window-pointwise}
and then integrating in $t$ and summing over $k$ proves
\eqref{eq:general-window-B}.

\medskip
\noindent$\mathrm{(ii)}$ For $k\in\Z$, almost every $x\in\R^n$ and
$c_1\,2^{-k/2}<t<c_2\,2^{-k/2}$, using the change of variables
$s=2^{-k/2}u$ and the inclusions in (i), we obtain
\begin{align*}
 &2^{k\alpha/2}\left[\int_1^{\sqrt2}\fint_{B(x,2^{-k/2})}
 |W^{1/p}(x)(\Phi_{2^{-k/2}u}*\vec g)(y)|^r\,\mathrm d y\,
 \frac{\mathrm d u}{u\log_e\sqrt2}\right]^{1/r}\\
 &\quad\asymp\left[\int_{2^{-k/2}}^{2^{(1-k)/2}}
 \fint_{B(x,2^{-k/2})}
 |s^{-\alpha}W^{1/p}(x)(\Phi_s*\vec g)(y)|^r\,\mathrm d y\,
 \frac{\mathrm d s}{s}\right]^{1/r}\\
 &\quad\lesssim\left[\fint_{t/2}^t\fint_{B(x,t)}
 |s^{-\alpha}W^{1/p}(x)(\Phi_s*\vec g)(y)|^r\,\mathrm d y\,
 \mathrm d s\right]^{1/r}.
\end{align*}
Applying the summation and integration argument in (i) to this
inequality, we obtain \eqref{eq:general-window-F-W}.
Taking $L^p$ first, we obtain \eqref{eq:general-window-B-W}.
\end{proof}

\Needspace{4\baselineskip}
The following proposition proves one direction of the characterizations
of homogeneous matrix-weighted Besov and Triebel--Lizorkin spaces
via Whitney averages.

\begin{proposition}
\label{prop:general-reverse}
Let $R\in\mathbb N_0$ and let $\Phi\in\cS_R$ satisfy
\eqref{eq:general-Tauberian}, with $\Phi_s(x)=s^{-n}\Phi(x/s)$.
Let $\alpha\in\R$,
$0<p<\infty$, $0<q,r\le\infty$, and $W\in\mathcal A_p$.
Let $\mathbb A=\{A_Q\}_{Q\in\D}$ be a sequence of reducing
operators of order $p$ for $W$.
For $\vec g\in(\cS_R')^m$, set
$\vec f:=\bigl.\vec g\bigr|_{(\cS_\infty)^m}\in(\cS_\infty')^m$.

\smallskip
\noindent$\mathrm{(i)}$ For the Triebel--Lizorkin spaces, we have
\begin{align}
 \|\vec f\|_{\dot F^\alpha_{p,q}(\mathbb A)}
 &\lesssim\left\|\left\{\int_0^\infty
 \left[\fint_{t/2}^t\fint_{B(\cdot,t)}
       |s^{-\alpha}\mathbb A(t,\cdot)(\Phi_s*\vec g)(y)|^r\,\mathrm d y\,\mathrm d s
 \right]^{q/r}\frac{\mathrm d t}{t}\right\}^{1/q}\right\|_{L^p(\R^n)},\label{eq:general-reverse-F}\\
 \|\vec f\|_{\dot F^\alpha_{p,q}(\mathbb A)}
 &\lesssim\left\|\left\{\int_0^\infty
 \left[\fint_{t/2}^t\fint_{B(\cdot,t)}
       |s^{-\alpha}W^{1/p}(\cdot)(\Phi_s*\vec g)(y)|^r\,\mathrm d y\,\mathrm d s
 \right]^{q/r}\frac{\mathrm d t}{t}\right\}^{1/q}\right\|_{L^p(\R^n)}.\notag
\end{align}

\smallskip
\noindent$\mathrm{(ii)}$ For the Besov spaces, we have
\begin{align*}
 \|\vec f\|_{\dot B^\alpha_{p,q}(\mathbb A)}
 &\lesssim\left\{\int_0^\infty
 \left\|\left[\fint_{t/2}^t\fint_{B(\cdot,t)}
       |s^{-\alpha}\mathbb A(t,\cdot)(\Phi_s*\vec g)(y)|^r\,\mathrm d y\,\mathrm d s
 \right]^{1/r}\right\|_{L^p(\R^n)}^q\frac{\mathrm d t}{t}\right\}^{1/q},\\
 \|\vec f\|_{\dot B^\alpha_{p,q}(\mathbb A)}
 &\lesssim\left\{\int_0^\infty
 \left\|\left[\fint_{t/2}^t\fint_{B(\cdot,t)}
       |s^{-\alpha}W^{1/p}(\cdot)(\Phi_s*\vec g)(y)|^r\,\mathrm d y\,\mathrm d s
 \right]^{1/r}\right\|_{L^p(\R^n)}^q\frac{\mathrm d t}{t}\right\}^{1/q}.
\end{align*}
All four estimates hold without the restriction $R+1>\alpha+d/p$,
where $d\in[0,n)$ is the $\mathcal A_p$-dimension of $W$
fixed in \eqref{eq:dimension-exponents}.
\end{proposition}

\begin{proof}
\noindent$\mathrm{(i)}$ Let $\mathbb A=\{A_Q\}_{Q\in\D}$ be a sequence of reducing
operators of order $p$ for $W$,
and let $\delta_2,\Delta$ be as in
\eqref{eq:dimension-exponents}.

\medskip
\noindent\emph{Step 1.}
Fix $\varphi\in\cS$ satisfying \eqref{eq:partition}, with
$\varphi_j(x)=2^{jn}\varphi(2^jx)$ for $j\in\Z$.
Choose the parameters as in \eqref{eq:general-parameters}:
\[
 0<\theta<\min\{1,p,q,r\},\qquad a>\Delta+n/\theta,
 \qquad N\ge a,\quad (N+\alpha-\delta_2)\theta>n.
\]
Let $\Lambda$ be the annular reproducing kernel in
Lemma~\ref{lem:general-reproduction}. For $\xi\ne0$ and
$u\in[1,\sqrt2]$, using
\eqref{eq:general-annular-partition}, we have
\[
 \widehat\varphi(2^{-j}\xi)
 =\sum_{k\in\Z}\widehat\varphi(2^{-j}\xi)
   \widehat\Lambda(2^{-k/2}u\xi)
   \widehat\Phi(2^{-k/2}u\xi).
\]
The Fourier supports imply
that only $|k-2j|\le C_0$ occur, where $C_0\in\mathbb N_0$ is independent
of $j$ and $u\in[1,\sqrt2]$.
Choose $\widetilde{\vec g}\in(\cS')^m$ such that
$\bigl.\widetilde{\vec g}\bigr|_{(\cS_R)^m}=\vec g$.
Then $\bigl.\widetilde{\vec g}\bigr|_{(\cS_\infty)^m}
=\vec f\in(\cS_\infty')^m$.
For $s>0$, set $\Lambda_s(x)=s^{-n}\Lambda(x/s)$.
Using $\Phi_s*\widetilde{\vec g}=\Phi_s*\vec g$ for $s>0$ and
$\varphi_j*\widetilde{\vec g}=\varphi_j*\vec f$, we obtain
\begin{equation}\label{eq:general-finite-reconstruction}
 \varphi_j*\vec f=
 \sum_{|k-2j|\le C_0}
 (\varphi_j*\Lambda_{2^{-k/2}u})*(\Phi_{2^{-k/2}u}*\vec g).
\end{equation}
For these indices,
\[
 2^{-j}\asymp2^{-k/2},\qquad
 \|\mathbb A_j(x)\mathbb A(2^{-k/2},x)^{-1}\|\lesssim1,
\]
and
\[
 \int_{\R^n}|(\varphi_j*\Lambda_{2^{-k/2}u})(v)|
                  (1+2^{k/2}|v|)^a\,\mathrm d v\lesssim1.
\]
Consequently
\[
 2^{j\alpha}|\mathbb A_j(x)(\varphi_j*\vec f)(x)| \lesssim\sum_{|k-2j|\le C_0} 2^{k\alpha/2}\sup_{y\in\R^n} \frac{|\mathbb A(2^{-k/2},x)(\Phi_{2^{-k/2}u}*\vec g)(y)|} {(1+2^{k/2}|x-y|)^a}.
\]
Hence, for $\theta<1$, we have
\[
 \big[2^{j\alpha}|\mathbb A_j(x)(\varphi_j*\vec f)(x)|\big]^\theta \lesssim\sum_{|k-2j|\le C_0} 2^{k\alpha\theta/2}\int_1^{\sqrt2}\left[\sup_{y\in\R^n}\frac{|\mathbb A(2^{-k/2},x)(\Phi_{2^{-k/2}u}*\vec g)(y)|}{(1+2^{k/2}|x-y|)^a}\right]^{\theta} \frac{\mathrm d u}{u\log_e\sqrt2}.
\]

\medskip
\noindent\emph{Step 2.}
Using \eqref{eq:general-fine-average} and Lemma~\ref{lem:FS}
with $p/\theta>1$ and $q/\theta>1$, together with
\[
 \sum_{h\ge0}2^{-h[(N+\alpha-\delta_2)\theta-n]/2}<\infty,
\]
we obtain
\begin{align*}
 &\left\|\left\{\sum_{j\in\Z}
    \big[2^{j\alpha}|\mathbb A_j(\cdot)(\varphi_j*\vec f)(\cdot)|\big]^q
                  \right\}^{\theta/q}\right\|_{L^{p/\theta}(\R^n)}\\
 &\quad\lesssim\left\|\left\{\sum_{k\in\Z}\Bigg[
  2^{k\alpha\theta/2}\int_1^{\sqrt2}\left[\sup_{y\in\R^n}\frac{|\mathbb A(2^{-k/2},\cdot)(\Phi_{2^{-k/2}u}*\vec g)(y)|}{(1+2^{k/2}|\cdot-y|)^a}\right]^{\theta}
                  \frac{\mathrm d u}{u\log_e\sqrt2}
                  \Bigg]^{q/\theta}\right\}^{\theta/q}\right\|_{L^{p/\theta}(\R^n)}\\
 &\quad\lesssim\left\|\left\{\sum_{k\in\Z}
 2^{k\alpha q/2}\left[\int_1^{\sqrt2}\fint_{B(\cdot,2^{-k/2})}
       |\mathbb A(2^{-k/2},\cdot)(\Phi_{2^{-k/2}u}*\vec g)(y)|^r\,\mathrm d y\,
       \frac{\mathrm d u}{u\log_e\sqrt2}\right]^{q/r}
                  \right\}^{\theta/q}\right\|_{L^{p/\theta}(\R^n)}.
\end{align*}
Taking the $1/\theta$ power and using
\eqref{eq:general-window-F}, we obtain
\eqref{eq:general-reverse-F}.

Using \eqref{eq:general-fine-average-W}, Lemma~\ref{lem:FS}
with exponents $p/\theta$ and $q/\theta$, and
\eqref{eq:general-window-F-W}, we obtain
\begin{align*}
 &\|\vec f\|_{\dot F^\alpha_{p,q}(\mathbb A)} \lesssim\left\|\left\{\sum_{k\in\Z}2^{k\alpha q/2}
 \left[\int_1^{\sqrt2}\fint_{B(\cdot,2^{-k/2})}
 |W^{1/p}(\cdot)(\Phi_{2^{-k/2}u}*\vec g)(y)|^r\,\mathrm d y\,
 \frac{\mathrm d u}{u\log_e\sqrt2}\right]^{q/r}
 \right\}^{1/q}\right\|_{L^p(\R^n)}\\
 &\quad\lesssim\left\|\left\{\int_0^\infty
 \left[\fint_{t/2}^t\fint_{B(\cdot,t)}
 |s^{-\alpha}W^{1/p}(\cdot)(\Phi_s*\vec g)(y)|^r\,\mathrm d y\,
 \mathrm d s\right]^{q/r}\frac{\mathrm d t}{t}
 \right\}^{1/q}\right\|_{L^p(\R^n)}.
\end{align*}

\medskip
\noindent$\mathrm{(ii)}$ The proof is similar to that of (i), and we omit the details.
\end{proof}

\Needspace{5\baselineskip}
The following proposition proves the other direction of the characterizations
of homogeneous matrix-weighted Besov and Triebel--Lizorkin spaces
via Whitney averages.

\begin{proposition}
\label{prop:general-upper}
Let $\alpha\in\R$, $0<p<\infty$, $0<q,r\le\infty$,
$W\in\mathcal A_p$, and let $\mathbb A$ be reducing operators
of order $p$ for $W$. Let $R\in\mathbb N_0$ and
$\Phi\in\cS_R$ satisfy Assumption~\ref{ass:general-Phi}.
Set $\Phi_s(x)=s^{-n}\Phi(x/s)$ for $s>0$, and use
$\mathcal I_R$ from Lemma~\ref{lem:general-realization}.

\smallskip
\noindent$\mathrm{(i)}$ For every $\vec f\in\dot F^\alpha_{p,q}(W)$,
\begin{align}
 \left\|\left\{\int_0^\infty
 \left[\fint_{t/2}^t\fint_{B(\cdot,t)}
       |s^{-\alpha}\mathbb A(t,\cdot)(\Phi_s*\mathcal I_R\vec f)(y)|^r\,\mathrm d y\,\mathrm d s
 \right]^{q/r}\frac{\mathrm d t}{t}\right\}^{1/q}\right\|_{L^p(\R^n)}
 &\lesssim\|\vec f\|_{\dot F^\alpha_{p,q}(\mathbb A)},\label{eq:general-upper-F}\\
 \left\|\left\{\int_0^\infty
 \left[\fint_{t/2}^t\fint_{B(\cdot,t)}
       |s^{-\alpha}W^{1/p}(\cdot)(\Phi_s*\mathcal I_R\vec f)(y)|^r\,\mathrm d y\,\mathrm d s
 \right]^{q/r}\frac{\mathrm d t}{t}\right\}^{1/q}\right\|_{L^p(\R^n)}
 &\lesssim\|\vec f\|_{\dot F^\alpha_{p,q}(\mathbb A)}.\notag
\end{align}

\smallskip
\noindent$\mathrm{(ii)}$ For every $\vec f\in\dot B^\alpha_{p,q}(W)$,
\begin{align}
 \left\{\int_0^\infty
 \left\|\left[\fint_{t/2}^t\fint_{B(\cdot,t)}
       |s^{-\alpha}\mathbb A(t,\cdot)(\Phi_s*\mathcal I_R\vec f)(y)|^r\,\mathrm d y\,\mathrm d s
 \right]^{1/r}\right\|_{L^p(\R^n)}^q\frac{\mathrm d t}{t}\right\}^{1/q}
 &\lesssim\|\vec f\|_{\dot B^\alpha_{p,q}(\mathbb A)},\label{eq:general-upper-B}\\
 \left\{\int_0^\infty
 \left\|\left[\fint_{t/2}^t\fint_{B(\cdot,t)}
       |s^{-\alpha}W^{1/p}(\cdot)(\Phi_s*\mathcal I_R\vec f)(y)|^r\,\mathrm d y\,\mathrm d s
 \right]^{1/r}\right\|_{L^p(\R^n)}^q\frac{\mathrm d t}{t}\right\}^{1/q}
 &\lesssim\|\vec f\|_{\dot B^\alpha_{p,q}(\mathbb A)}.\notag
\end{align}
\end{proposition}

\begin{proof}
\noindent$\mathrm{(i)}$ Let $\mathbb A=\{A_Q\}_{Q\in\D}$ be a sequence of reducing
operators of order $p$ for $W$, and let
$\delta_1,\delta_2,\Delta$ be as in
\eqref{eq:dimension-exponents}. Use
$\varphi_k(x)=2^{kn}\varphi(2^kx)$ and
$\Psi_k(x)=2^{kn}\Psi(2^kx)$ for $k\in\Z$, where
$\varphi\in\cS$ satisfies \eqref{eq:partition} and
$\Psi\in\cS$ has compact Fourier support away from the origin
and satisfies $\widehat\Psi=1$ on $\supp\widehat\varphi$. Choose
\[
 0<\theta<\min\{1,p,q\},\qquad a>\Delta+n/\theta,
 \qquad H>\max\{0,2a+\delta_2+\Delta-\alpha\}.
\]

\medskip
\noindent\emph{Step 1.}
By \eqref{eq:peetre-point} and Lemma~\ref{lem:FS},
\begin{equation}\label{eq:general-Peetre-norms}
 \left\|\left\{\sum_{k\in\Z}\left[2^{k\alpha}\sup_{v\in\R^n}
 \frac{|\mathbb A_k(\cdot)(\varphi_k*\vec f)(v)|}
 {(1+2^k|\cdot-v|)^a}\right]^q\right\}^{1/q}\right\|_{L^p(\R^n)}
 \lesssim\|\vec f\|_{\dot F^\alpha_{p,q}(\mathbb A)}.
\end{equation}
By Lemma~\ref{lem:general-realization}, we have
\begin{equation*}
 \Phi_s*\mathcal I_R\vec f
 =\sum_{k\in\Z}(\Phi_s*\Psi_k)*(\varphi_k*\vec f).
\end{equation*}
The series converges absolutely at each $s,y$ by testing
$\mathcal I_R\vec f$ against $\Phi_s(y-\cdot)\in\cS_R$.

Fix $2^{-j}\le t<2^{-j+1}$, $t/2<s<t$ and $y\in B(x,t)$.
Choose $L>n+a$ in Lemma~\ref{lem:general-kernel}. Using its kernel
estimate and \eqref{eq:comparison}, we obtain
\begin{align}
 \int_{\R^n}|(\Phi_s*\Psi_k)(y-v)|(1+2^k|x-v|)^a\,\mathrm d v
 &\lesssim
 \begin{cases}
  2^{-(R+1)(j-k)},&k\le j,\\
  2^{-(H-a)(k-j)},&k>j,
 \end{cases}\label{eq:general-upper-kernel-integral}\\
 \|\mathbb A_j(x)\mathbb A_k(x)^{-1}\|
 &\lesssim
 \begin{cases}2^{\delta_1(j-k)},&k\le j,\\
               2^{\delta_2(k-j)},&k>j.
 \end{cases}\notag
\end{align}
Since $s^{-\alpha}2^{-k\alpha}\asymp2^{(j-k)\alpha}$ and
$\mathbb A(t,x)=\mathbb A_j(x)$, for $\theta<1$ we have
\begin{align*}
 &\left[\fint_{t/2}^t\fint_{B(x,t)}
 |s^{-\alpha}\mathbb A(t,x)(\Phi_s*\mathcal I_R\vec f)(y)|^r
 \,\mathrm d y\,\mathrm d s\right]^{\theta/r}\\
 &\quad\lesssim\sum_{k\le j}
 2^{-(R+1-\alpha-\delta_1)\theta(j-k)}
 \left[2^{k\alpha}\sup_{v\in\R^n}
 \frac{|\mathbb A_k(x)(\varphi_k*\vec f)(v)|}
 {(1+2^k|x-v|)^a}\right]^\theta\\
 &\qquad+\sum_{k>j}
 2^{-(H+\alpha-a-\delta_2)\theta(k-j)}
 \left[2^{k\alpha}\sup_{v\in\R^n}
 \frac{|\mathbb A_k(x)(\varphi_k*\vec f)(v)|}
 {(1+2^k|x-v|)^a}\right]^\theta.
\end{align*}
Since $R+1-\alpha-\delta_1>0$ and $H+\alpha-a-\delta_2>0$,
\[
 \sum_{k\le j}2^{-(R+1-\alpha-\delta_1)\theta(j-k)}
 +\sum_{k>j}2^{-(H+\alpha-a-\delta_2)\theta(k-j)}\lesssim1,
 \qquad
 \int_{2^{-j}}^{2^{-j+1}}\frac{\mathrm d t}{t}=\log_e2.
\]
By Young's inequality with $q/\theta>1$ and
\eqref{eq:general-Peetre-norms},
\begin{align*}
 &\left\|\left\{\int_0^\infty
 \left[\fint_{t/2}^t\fint_{B(\cdot,t)}
 |s^{-\alpha}\mathbb A(t,\cdot)(\Phi_s*\mathcal I_R\vec f)(y)|^r
 \,\mathrm d y\,\mathrm d s\right]^{q/r}
 \frac{\mathrm d t}{t}\right\}^{\theta/q}
 \right\|_{L^{p/\theta}(\R^n)}\\
 &\quad\lesssim\left\|\left\{\sum_{k\in\Z}
 \left[2^{k\alpha}\sup_{v\in\R^n}
 \frac{|\mathbb A_k(\cdot)(\varphi_k*\vec f)(v)|}
 {(1+2^k|\cdot-v|)^a}\right]^q
 \right\}^{\theta/q}\right\|_{L^{p/\theta}(\R^n)}
 \lesssim\|\vec f\|_{\dot F^\alpha_{p,q}(\mathbb A)}^\theta.
\end{align*}
Taking the $1/\theta$ power, we obtain \eqref{eq:general-upper-F}.

\medskip
\noindent\emph{Step 2.}
Let $x\in Q\in\D_j$. For every $z\in Q$, \eqref{eq:general-upper-kernel-integral}
holds with $z$ in place of $x$, since $|y-z|\lesssim2^{-j}$.
Using $\mathbb A_j(z)=A_Q$, we obtain
\begin{align*}
 &\left[\fint_{t/2}^t\fint_{B(x,t)}
 |s^{-\alpha}W^{1/p}(x)(\Phi_s*\mathcal I_R\vec f)(y)|^r\,\mathrm d y\,\mathrm d s
 \right]^{1/r}\\
 &\quad\lesssim\|W^{1/p}(x)A_Q^{-1}\|\Bigg[
 \sum_{k\le j}2^{-(R+1-\alpha-\delta_1)(j-k)}2^{k\alpha}
 \sup_{\substack{z\in Q\\v\in\R^n}}
 \frac{|\mathbb A_k(z)(\varphi_k*\vec f)(v)|}{(1+2^k|z-v|)^a}\\
 &\qquad\qquad+\sum_{k>j}2^{-(H+\alpha-a-\delta_2)(k-j)}2^{k\alpha}
 \sup_{\substack{z\in Q\\v\in\R^n}}
 \frac{|\mathbb A_k(z)(\varphi_k*\vec f)(v)|}{(1+2^k|z-v|)^a}
 \Bigg].
\end{align*}
For each $j\in\Z$, the expression in brackets is constant on every
$Q\in\D_j$. For $k\in\Z$, $x,z\in Q$ and $v\in\R^n$,
using \eqref{eq:comparison} and the triangle inequality, we have
\[
 \begin{aligned}
 \|\mathbb A_k(z)\mathbb A_k(x)^{-1}\|
 &\lesssim(1+2^k|x-z|)^\Delta
 \lesssim2^{\Delta(k-j)_+},\\
 \frac{1+2^k|x-v|}{1+2^k|z-v|}
 &\le1+2^k|x-z|\lesssim2^{(k-j)_+}.
 \end{aligned}
\]
Hence, for $x\in Q$,
\[
 \sup_{\substack{z\in Q\\v\in\R^n}}
 \frac{|\mathbb A_k(z)(\varphi_k*\vec f)(v)|}{(1+2^k|z-v|)^a}
 \lesssim2^{(a+\Delta)(k-j)_+}
 \sup_{v\in\R^n}
 \frac{|\mathbb A_k(x)(\varphi_k*\vec f)(v)|}{(1+2^k|x-v|)^a}.
\]
For $k>j$,
\[
 2^{-(H+\alpha-a-\delta_2)(k-j)}2^{(a+\Delta)(k-j)}
 =2^{-(H+\alpha-2a-\delta_2-\Delta)(k-j)}.
\]
Since $R+1-\alpha-\delta_1>0$ and
$H+\alpha-2a-\delta_2-\Delta>0$, applying
\eqref{eq:constant-F}, Young's inequality with exponent $q/\theta$,
and \eqref{eq:general-Peetre-norms}, we obtain
\begin{align*}
 &\left\|\left\{\int_0^\infty
 \left[\fint_{t/2}^t\fint_{B(\cdot,t)}
 |s^{-\alpha}W^{1/p}(\cdot)(\Phi_s*\mathcal I_R\vec f)(y)|^r
 \,\mathrm d y\,\mathrm d s\right]^{q/r}
 \frac{\mathrm d t}{t}\right\}^{1/q}\right\|_{L^p(\R^n)}\\
 &\quad\lesssim\left[
 \sum_{\ell=0}^\infty2^{-(R+1-\alpha-\delta_1)\theta\ell}
 +\sum_{\ell=1}^\infty
 2^{-(H+\alpha-2a-\delta_2-\Delta)\theta\ell}
 \right]^{1/\theta}\\
 &\qquad\times\left\|\left\{\sum_{k\in\Z}
 \left[2^{k\alpha}\sup_{v\in\R^n}
 \frac{|\mathbb A_k(\cdot)(\varphi_k*\vec f)(v)|}
 {(1+2^k|\cdot-v|)^a}\right]^q
 \right\}^{1/q}\right\|_{L^p(\R^n)}
 \lesssim\|\vec f\|_{\dot F^\alpha_{p,q}(\mathbb A)}.
\end{align*}

\medskip
\noindent$\mathrm{(ii)}$ The proof is similar to that of (i), and we omit the details.
\end{proof}

Theorem~\ref{thm:general-Whitney} follows from
Propositions~\ref{prop:general-reverse} and~\ref{prop:general-upper},
together with \eqref{eq:known-equivalence}.

\section{Duality}\label{sec:Z}

We identify the continuous duals of the spaces in
Definitions~\ref{def:Z} and~\ref{def:T}, using the dyadic
characterizations in Subsection~\ref{subsec:ZT-dyadic}.
For the corresponding scalar theory, see
\cite{ABH,HuangTent,HaardtTent}. We use inverse reducing
operators and, for $p>1$, the dual matrix weight, as in
\cite{RoudenkoDuality}.

For $\vec u,\vec v\in\C^m$, set
\[
 \langle\vec u,\vec v\rangle_{\C^m}
 =\sum_{\ell=1}^m u_\ell\overline{v_\ell}.
\]

Whenever the integral is absolutely convergent, put
\begin{equation}\label{eq:Z-pairing}
 \Lambda_{\vec G}(\vec F)
 :=\int_0^\infty\int_{\R^n}
 \langle\vec F(s,y),\vec G(s,y)\rangle_{\C^m}
 \,\mathrm d y\,\frac{\mathrm d s}{s}.
\end{equation}

\subsection{Duality of matrix-weighted \texorpdfstring{$Z$}{Z}-spaces: full range}
\label{subsec:Z-full-duality}

We use the unweighted duality theorem in
\cite[Theorem 3.16]{ABH} and the dyadic characterization to obtain
duality for matrix-weighted $Z$-spaces with $0<p,q<\infty$ and
$1\le r<\infty$.

Let $\mathbb A=\{A_Q\}_{Q\in\D}$ be as in
Definition~\ref{def:Z}. The inverse family
$\mathbb A^{-1}:=\{A_Q^{-1}\}_{Q\in\D}$ satisfies
\eqref{eq:comparison} with $\delta_1$ and $\delta_2$
interchanged. Indeed, since the matrices are Hermitian,
\begin{equation}\label{eq:Z-inverse-comparison}
 \|\mathbb A_j(x)^{-1}\mathbb A_k(z)\|
 =\|\mathbb A_k(z)\mathbb A_j(x)^{-1}\|
 \lesssim2^{\delta_2(j-k)_++\delta_1(k-j)_+}
              (1+2^{\min(j,k)}|x-z|)^\Delta.
\end{equation}
Here $\delta_1$, $\delta_2$, and $\Delta$ are as in
Definition~\ref{def:Z}.
Since $p'=\infty$ when $0<p\le1$, we first define the
$Z$-space with an infinite first index.

\begin{definition}\label{def:Z-infinite-spatial}
Let $\mathbb A=\{A_Q\}_{Q\in\D}$ be as in
Definition~\ref{def:Z}. Let $\widehat Q$ and $\nu_Q$ be the same
as in \eqref{eq:Z-boxes}.
For $\gamma\in\R$ and $0<b,c\le\infty$, let
$Z^{\infty,b,c}_\gamma(\mathbb A^{-1})$ consist of all measurable
$\vec G:\R^{n+1}_+\to\C^m$ such that
\[
 \|\vec G\|_{Z^{\infty,b,c}_\gamma(\mathbb A^{-1})}
 :=\left\{\sum_{j\in\Z}\left[\sup_{Q\in\D_j}2^{j\gamma}
 \|A_Q^{-1}\vec G\|_{L^c(\widehat Q,\nu_Q;\C^m)}
 \right]^b\right\}^{1/b}<\infty,
\]
with the usual modifications if $b$ and/or $c$ are infinite.
\end{definition}

Applying \eqref{eq:local-matrix-comparison} to $\mathbb A^{-1}$
and using \cite[Proposition 3.12]{ABH} with spatial index $\infty$, we obtain
\[
 \|\vec G\|_{Z^{\infty,b,c}_\gamma(\mathbb A^{-1})}
 \asymp\left\{\int_0^\infty
 \left[\operatorname*{ess\,sup}_{x\in\R^n}
 \left(\fint_{t/2}^{t}\fint_{B(x,t)}
 |s^{-\gamma}\mathbb A(t,x)^{-1}\vec G(s,y)|^c
 \,\mathrm d y\,\mathrm d s\right)^{1/c}\right]^b
 \frac{\mathrm d t}{t}\right\}^{1/b}.
\]

The following finite-dimensional form of unweighted $Z$-space duality
follows from \cite[Theorem 3.16]{ABH}, applied to each coordinate.

\Needspace{4\baselineskip}
\begin{lemma}\label{lem:Z-unweighted-dual}
Let $\alpha\in\R$, $0<p,q<\infty$, $1\le r<\infty$, and let
$\mathbb I=\{I_m\}_{Q\in\D}$. Under the pairing \eqref{eq:Z-pairing},
\[
 \bigl(Z^{p,q,r}_\alpha(\mathbb I)\bigr)^*
 \cong Z^{p',q',r'}_{-\alpha+n(1/p-1)_+}(\mathbb I).
\]
That is, for every $\vec F\in Z^{p,q,r}_\alpha(\mathbb I)$ and
$\vec G\in Z^{p',q',r'}_{-\alpha+n(1/p-1)_+}(\mathbb I)$,
\begin{align*}
 \int_0^\infty\int_{\R^n}
 |\langle\vec F(s,y),\vec G(s,y)\rangle_{\C^m}|
 \,\mathrm d y\,\frac{\mathrm d s}{s}
 &\lesssim\|\vec F\|_{Z^{p,q,r}_\alpha(\mathbb I)}
 \|\vec G\|_{Z^{p',q',r'}_{-\alpha+n(1/p-1)_+}(\mathbb I)},\\
 \|\Lambda_{\vec G}\|_{(Z^{p,q,r}_\alpha(\mathbb I))^*}
 &\asymp\|\vec G\|_{Z^{p',q',r'}_{-\alpha+n(1/p-1)_+}(\mathbb I)}.
\end{align*}
Every bounded linear functional has a unique representing function
modulo equality almost everywhere.
\end{lemma}

We transfer this result to matrix families by multiplying the functions
by $A_Q$ on each Whitney box $\widehat Q$.

\Needspace{4\baselineskip}
\begin{theorem}\label{thm:Z-full-dual}
Let $\alpha\in\R$, $0<p,q<\infty$, $1\le r<\infty$, and let
$\mathbb A$ be as in Definition~\ref{def:Z}.
Under the pairing \eqref{eq:Z-pairing}, the following identifications hold.

\smallskip
\noindent$\mathrm{(i)}$ For the matrix family $\mathbb A$,
\begin{equation*}
 \bigl(Z^{p,q,r}_\alpha(\mathbb A)\bigr)^*
 \cong Z^{p',q',r'}_{-\alpha+n(1/p-1)_+}(\mathbb A^{-1}).
\end{equation*}
That is, for every $\vec F\in Z^{p,q,r}_\alpha(\mathbb A)$ and
$\vec G\in Z^{p',q',r'}_{-\alpha+n(1/p-1)_+}(\mathbb A^{-1})$,
\begin{align*}
 \int_0^\infty\int_{\R^n}
 |\langle\vec F(s,y),\vec G(s,y)\rangle_{\C^m}|
 \,\mathrm d y\,\frac{\mathrm d s}{s}
 &\lesssim\|\vec F\|_{Z^{p,q,r}_\alpha(\mathbb A)}
 \|\vec G\|_{Z^{p',q',r'}_{-\alpha+n(1/p-1)_+}(\mathbb A^{-1})}<\infty,\\
 \|\Lambda_{\vec G}\|_{(Z^{p,q,r}_\alpha(\mathbb A))^*}
 &\asymp\|\vec G\|_{Z^{p',q',r'}_{-\alpha+n(1/p-1)_+}(\mathbb A^{-1})}.
\end{align*}

\smallskip
\noindent$\mathrm{(ii)}$ If $W\in\mathcal A_p$ and $\mathbb A$ is a sequence
of reducing operators of order $p$ for $W$, then
\begin{equation}\label{eq:Z-full-dual-W}
 \bigl(Z^{p,q,r}_\alpha(W)\bigr)^*
 \cong Z^{p',q',r'}_{-\alpha+n(1/p-1)_+}(\mathbb A^{-1}).
\end{equation}
In each case, every bounded linear functional has a unique
representing function modulo equality almost everywhere.
\end{theorem}

\begin{proof}
\noindent$\mathrm{(i)}$ Let $\mathbb I=\{I_m\}_{Q\in\D}$. For $Q\in\D$, recall that
$\widehat Q=[\ell(Q),2\ell(Q))\times Q$. Define
\[
 (\mathcal U\vec F)(s,y)=A_Q\vec F(s,y),\qquad
 (\mathcal U^{-1}\vec H)(s,y)=A_Q^{-1}\vec H(s,y),
 \quad(s,y)\in\widehat Q.
\]
By Lemma~\ref{lem:Z-dyadic-A}~(ii) and
Definition~\ref{def:Z-infinite-spatial},
\[
\|\mathcal U\vec F\|_{Z^{p,q,r}_\alpha(\mathbb I)}
 \asymp\|\vec F\|_{Z^{p,q,r}_\alpha(\mathbb A)},\qquad 
 \|\mathcal U^{-1}\vec G\|_{Z^{p',q',r'}_{-\alpha+n(1/p-1)_+}(\mathbb I)}
 \asymp\|\vec G\|_{Z^{p',q',r'}_{-\alpha+n(1/p-1)_+}(\mathbb A^{-1})}.
\]
Since each $A_Q$ is Hermitian, on $\widehat Q$ we have
\[
 \langle\mathcal U\vec F,\mathcal U^{-1}\vec G\rangle_{\C^m}
 =\langle A_Q\vec F,A_Q^{-1}\vec G\rangle_{\C^m}
 =\langle\vec F,\vec G\rangle_{\C^m}.
\]
Using these identities and Lemma~\ref{lem:Z-unweighted-dual}, we obtain
the pairing estimate and the equivalence of the functional norm in (i).
For $\Lambda\in(Z^{p,q,r}_\alpha(\mathbb A))^*$, Lemma~\ref{lem:Z-unweighted-dual}
provides a unique
$\vec K\in Z^{p',q',r'}_{-\alpha+n(1/p-1)_+}(\mathbb I)$
representing $\Lambda\circ\mathcal U^{-1}$. Set $\vec G=\mathcal U\vec K$.
Then, for every $\vec F\in Z^{p,q,r}_\alpha(\mathbb A)$,
\[
 \Lambda(\vec F)
 =\int_0^\infty\int_{\R^n}
 \langle\mathcal U\vec F(s,y),\vec K(s,y)\rangle_{\C^m}
 \,\mathrm d y\,\frac{\mathrm d s}{s}
 =\Lambda_{\vec G}(\vec F).
\]
The uniqueness of $\vec K$ and the invertibility of $\mathcal U$
imply the uniqueness of $\vec G$.

\medskip
\noindent$\mathrm{(ii)}$ By Lemma~\ref{lem:Z-W-A}~(ii) and (i),
\[
 \bigl(Z^{p,q,r}_\alpha(W)\bigr)^*
 \cong\bigl(Z^{p,q,r}_\alpha(\mathbb A)\bigr)^*
 \cong Z^{p',q',r'}_{-\alpha+n(1/p-1)_+}(\mathbb A^{-1}).\qedhere
\]
\end{proof}

\begin{remark}\label{rem:Z-weight-dual}
Under the assumptions of Theorem~\ref{thm:Z-full-dual}, if
$p>1$, then $-\alpha+n(1/p-1)_+=-\alpha$.
Lemma~\ref{lem:Z-dual-reducing} and \eqref{eq:Z-W-A} identify
\eqref{eq:Z-full-dual-W} with
\[
 \bigl(Z^{p,q,r}_\alpha(W)\bigr)^*
 \cong Z^{p',q',r'}_{-\alpha}(W^{-p'/p}),
 \qquad 0<q<\infty,\quad1\le r<\infty.
\]
\end{remark}

\Needspace{4\baselineskip}
\subsection{Duality of tent spaces in the Banach range}

By \cite[Lemma 3.3]{ABH}, applied to $|\vec F|$, our unweighted tent
quasi-norm is equivalent to Huang's Whitney-averaged tent quasi-norm.
Taking $\beta_0=0$ in \cite[Theorem 5.4, pp.~926--927]{HuangTent},
we obtain the following
unweighted duality. The scalar result is applied to each coordinate.

\Needspace{4\baselineskip}
\begin{lemma}\label{lem:T-unweighted-dual}
Let $\alpha\in\R$, $1<p,q,r<\infty$, and let
$\mathbb I=\{I_m\}_{Q\in\D}$. Under the pairing
\eqref{eq:Z-pairing},
\[
 \bigl(T^{p,q,r}_\alpha(\mathbb I)\bigr)^*
 \cong T^{p',q',r'}_{-\alpha}(\mathbb I).
\]
That is, for every $\vec F\in T^{p,q,r}_\alpha(\mathbb I)$ and
$\vec G\in T^{p',q',r'}_{-\alpha}(\mathbb I)$,
\begin{align*}
 \int_0^\infty\int_{\R^n}
 |\langle\vec F(s,y),\vec G(s,y)\rangle_{\C^m}|
 \,\mathrm d y\,\frac{\mathrm d s}{s}
 &\lesssim\|\vec F\|_{T^{p,q,r}_\alpha(\mathbb I)}
 \|\vec G\|_{T^{p',q',r'}_{-\alpha}(\mathbb I)},\\
 \|\Lambda_{\vec G}\|_{(T^{p,q,r}_\alpha(\mathbb I))^*}
 &\asymp\|\vec G\|_{T^{p',q',r'}_{-\alpha}(\mathbb I)}.
\end{align*}
Every bounded linear functional has a unique representing function
modulo equality almost everywhere.
\end{lemma}

We apply this unweighted result to the matrix families of
Definition~\ref{def:Z} using the dyadic characterization.

\Needspace{4\baselineskip}
\begin{theorem}\label{thm:Z-dual-A}\label{thm:Z-dual-W}
Let $\alpha\in\R$, $1<p,q,r<\infty$, and let $\mathbb A$ be as in
Definition~\ref{def:Z}.

\smallskip
\noindent$\mathrm{(i)}$ Under the pairing \eqref{eq:Z-pairing},
\begin{equation}\label{eq:T-dual-A}
 \bigl(T^{p,q,r}_\alpha(\mathbb A)\bigr)^*
 \cong T^{p',q',r'}_{-\alpha}(\mathbb A^{-1}).
\end{equation}
That is, for every $\vec F\in T^{p,q,r}_\alpha(\mathbb A)$ and
$\vec G\in T^{p',q',r'}_{-\alpha}(\mathbb A^{-1})$,
\begin{align*}
 \int_0^\infty\int_{\R^n}
 |\langle\vec F(s,y),\vec G(s,y)\rangle_{\C^m}|
 \,\mathrm d y\,\frac{\mathrm d s}{s}
 &\lesssim\|\vec F\|_{T^{p,q,r}_\alpha(\mathbb A)}
 \|\vec G\|_{T^{p',q',r'}_{-\alpha}(\mathbb A^{-1})}<\infty,\\
 \|\Lambda_{\vec G}\|_{(T^{p,q,r}_\alpha(\mathbb A))^*}
 &\asymp\|\vec G\|_{T^{p',q',r'}_{-\alpha}(\mathbb A^{-1})}.
\end{align*}

\smallskip
\noindent$\mathrm{(ii)}$ If $W\in\mathcal A_p$ and $\mathbb A$ is a sequence
of reducing operators of order $p$ for $W$, then
\begin{equation}\label{eq:T-dual-W}
 \bigl(T^{p,q,r}_\alpha(W)\bigr)^*
 \cong T^{p',q',r'}_{-\alpha}(W^{-p'/p}).
\end{equation}
In each case, every bounded linear functional has a unique
representing function modulo equality almost everywhere.
\end{theorem}

\begin{proof}
\noindent$\mathrm{(i)}$ Let $\mathbb I=\{I_m\}_{Q\in\D}$. For $Q\in\D$, recall that
$\widehat Q=[\ell(Q),2\ell(Q))\times Q$. Define
\[
 (\mathcal U\vec F)(s,y)=A_Q\vec F(s,y),\qquad
 (\mathcal U^{-1}\vec H)(s,y)=A_Q^{-1}\vec H(s,y),
 \quad(s,y)\in\widehat Q.
\]
By Lemma~\ref{lem:Z-dyadic-A}~(i),
\[
 \|\mathcal U\vec F\|_{T^{p,q,r}_\alpha(\mathbb I)}
 \asymp\|\vec F\|_{T^{p,q,r}_\alpha(\mathbb A)},\qquad
 \|\mathcal U^{-1}\vec G\|_{T^{p',q',r'}_{-\alpha}(\mathbb I)}
 \asymp\|\vec G\|_{T^{p',q',r'}_{-\alpha}(\mathbb A^{-1})}.
\]
Since each $A_Q$ is Hermitian, on $\widehat Q$ we have
\[
 \langle\mathcal U\vec F,\mathcal U^{-1}\vec G\rangle_{\C^m}
 =\langle A_Q\vec F,A_Q^{-1}\vec G\rangle_{\C^m}
 =\langle\vec F,\vec G\rangle_{\C^m}.
\]
Using these identities and Lemma~\ref{lem:T-unweighted-dual}, we obtain
the pairing estimate and the equivalence of the functional norm in (i).
For $\Lambda\in(T^{p,q,r}_\alpha(\mathbb A))^*$, Lemma~\ref{lem:T-unweighted-dual}
provides a unique
$\vec K\in T^{p',q',r'}_{-\alpha}(\mathbb I)$
representing $\Lambda\circ\mathcal U^{-1}$. Set $\vec G=\mathcal U\vec K$.
Then, for every $\vec F\in T^{p,q,r}_\alpha(\mathbb A)$,
\[
 \Lambda(\vec F)
 =\int_0^\infty\int_{\R^n}
 \langle\mathcal U\vec F(s,y),\vec K(s,y)\rangle_{\C^m}
 \,\mathrm d y\,\frac{\mathrm d s}{s}
 =\Lambda_{\vec G}(\vec F).
\]
The uniqueness of $\vec K$ and the invertibility of $\mathcal U$
imply the uniqueness of $\vec G$.

\medskip
\noindent$\mathrm{(ii)}$ Using Lemma~\ref{lem:Z-dual-reducing}, the dyadic
characterizations in Lemmas~\ref{lem:Z-dyadic-A}
and~\ref{lem:Z-W-A}, and (i), we obtain
\[
 (T^{p,q,r}_\alpha(W))^*
 \cong(T^{p,q,r}_\alpha(\mathbb A))^*
 \cong T^{p',q',r'}_{-\alpha}(\mathbb A^{-1})
 \cong T^{p',q',r'}_{-\alpha}(W^{-p'/p}).\qedhere
\]
\end{proof}

\section{Interpolation of matrix-weighted spaces}\label{sec:interpolation}

We recall the interpolation spaces used below; see
\cite[Sections C.1--C.3]{HNVW16}.
\begin{definition}\label{def:interpolation-spaces}
A pair $(X_1,X_2)$ of complex Banach spaces is called compatible
if both spaces are continuously embedded in a common Hausdorff
topological vector space. For such a pair, write
$X_1+X_2=\{\vec F_1+\vec F_2:\vec F_i\in X_i,\ i\in\{1,2\}\}$.

\smallskip
\noindent$\mathrm{(i)}$ For $\vec F\in X_1+X_2$ and $t>0$, define the $K$-functional by
\[
 K(t,\vec F;X_1,X_2)
 :=\inf_{\substack{\vec F=\vec F_1+\vec F_2\\
                   \vec F_i\in X_i,\ i\in\{1,2\}}}
 (\|\vec F_1\|_{X_1}+t\|\vec F_2\|_{X_2}).
\]
We equip $X_1+X_2$ with the norm
$\|\vec F\|_{X_1+X_2}:=K(1,\vec F;X_1,X_2)$.
For $0<\theta<1$ and $1<u<\infty$, the real interpolation space
$(X_1,X_2)_{\theta,u}$ consists of all $\vec F\in X_1+X_2$ for which
\[
 \|\vec F\|_{(X_1,X_2)_{\theta,u}}
 :=\left(\int_0^\infty
 [t^{-\theta}K(t,\vec F;X_1,X_2)]^u\frac{\mathrm d t}{t}\right)^{1/u}<\infty.
\]

\smallskip
\noindent$\mathrm{(ii)}$ For the complex method, put
$\mathbb S:=\{z\in\C:0<\operatorname{Re}z<1\}$.
Let $\mathcal H(X_1,X_2)$ consist of all bounded continuous maps
$\vec h:\overline{\mathbb S}\to X_1+X_2$ that are holomorphic in
$\mathbb S$ and whose boundary maps
$t\mapsto\vec h(it)$ and $t\mapsto\vec h(1+it)$ are bounded and
continuous from $\R$ into $X_1$ and $X_2$, respectively.
For $0<\theta<1$, the complex interpolation space is
\[
 [X_1,X_2]_\theta
 :=\{\vec h(\theta):\vec h\in\mathcal H(X_1,X_2)\},
\]
with norm
\[
 \|\vec F\|_{[X_1,X_2]_\theta}
 :=\inf_{\substack{\vec h\in\mathcal H(X_1,X_2)\\
                   \vec h(\theta)=\vec F}}
 \max\left\{\sup_{t\in\R}\|\vec h(it)\|_{X_1},
             \sup_{t\in\R}\|\vec h(1+it)\|_{X_2}\right\}.
\]
\end{definition}

\subsection{Geometric means and reducing operators}

The following definition is from \cite[Definition 2.20]{BCU}.

\begin{definition}\label{def:matrix-geometric-mean}
Let $A$ and $B$ be two positive definite Hermitian $m\times m$
matrices. For $0<t<1$, define the weighted geometric mean of
$A$ and $B$ by
\[
 A\#_t B:=A^{1/2}(A^{-1/2}BA^{-1/2})^t A^{1/2}.
\]
\end{definition}

Let $0<\theta<1$, $1<p_1,p_2<\infty$, and
$W_i\in\mathcal A_{p_i}$ for $i\in\{1,2\}$. Define
\begin{equation}\label{eq:intermediate-weight}
 \frac1{p_\theta}=\frac{1-\theta}{p_1}+\frac\theta{p_2},
 \qquad
 W_\theta
 =\bigl(W_1^{2/p_1}\#_\theta W_2^{2/p_2}\bigr)^{p_\theta/2}.
\end{equation}
The notation in \cite{BCU} uses multiplier weights: their
matrix weight $U$ corresponds to the weight $U^p$ in
Definition~\ref{def:Ap}. Thus their membership condition at
exponent $p$ is equivalent to $U^p\in\mathcal A_p$ here.

Let $h:\C^m\to[0,\infty)$ be continuous, with
$h(\lambda\vec v)=|\lambda|h(\vec v)$ for $\lambda\in\C$ and
$h(\vec v)>0$ for $\vec v\ne\vec0$. Define
\[
 h^*(\vec v)=\sup_{\vec u\in\C^m\setminus\{\vec0\}}
 \frac{\left|\langle\vec v,\vec u\rangle_{\C^m}\right|}{h(\vec u)},
 \qquad h^{**}(\vec v)=(h^*)^*(\vec v)
 =\sup_{\vec u\in\C^m\setminus\{\vec0\}}
 \frac{\left|\langle\vec v,\vec u\rangle_{\C^m}\right|}{h^*(\vec u)}.
\]
Then $h$ is comparable to the Euclidean norm,
$h^*$ and $h^{**}$ are norms, and $h^{**}\le h$. This follows
from the bipolar theorem applied to the closed absolutely convex
hull of the unit sublevel set of $h$; see also
\cite[Section 2]{BCU}.

\Needspace{4\baselineskip}
Parts (i) and (ii) of the following lemma follow from
\cite[Lemmas 2.18 and 2.21, Corollary 2.22]{BCU};
the identity in (iii) is recorded in
\cite[proof of Proposition 8.8]{BCU}.
The same arguments apply to complex Hermitian matrices,
with unitary diagonalization in place of orthogonal diagonalization.

\begin{lemma}\label{lem:ellipsoid-mean}
Let $P,Q$ be positive definite Hermitian $m\times m$ matrices,
$0<\theta<1$, and $B=(P^2\#_\theta Q^2)^{1/2}$. Define
$h(\vec v)=|P\vec v|^{1-\theta}|Q\vec v|^\theta$,
and let $h^{**}$ be its double dual norm defined above.

\smallskip
\noindent$\mathrm{(i)}$ There are an invertible matrix $S$ and a positive diagonal
matrix $D=\operatorname{diag}(d_1,\ldots,d_m)$, with
$d_k>0$ for $1\le k\le m$, such that, for every $\vec v\in\C^m$,
\begin{equation*}
 |PS\vec v|=|\vec v|,\qquad
 |QS\vec v|=|D\vec v|,\qquad
 |BS\vec v|=|D^\theta\vec v|.
\end{equation*}

\smallskip
\noindent$\mathrm{(ii)}$ For every $\vec v\in\C^m$,
\begin{equation}\label{eq:ellipsoid-envelope}
 |B\vec v|\asymp h^{**}(\vec v)\le h(\vec v).
\end{equation}
The equivalence constants depend only on $m$.

\smallskip
\noindent$\mathrm{(iii)}$ The inverse matrix satisfies
\begin{equation}\label{eq:inverse-mean}
 B^{-1}=(P^{-2}\#_\theta Q^{-2})^{1/2}.
\end{equation}

\end{lemma}

\Needspace{7\baselineskip}
Part (i) of the following lemma follows from
\cite[Proposition 8.8]{BCU}, applied to the multiplier weights
$W_i^{1/p_i}$ for $i\in\{1,2\}$; see also
\cite[Remark 8.10]{BCU} for the complex case.

\begin{lemma}\label{lem:intermediate-reducing}
Let $0<\theta<1$, $1<p_i<\infty$, and $W_i\in\mathcal A_{p_i}$
for $i\in\{1,2\}$. Define $p_\theta,W_\theta$
by \eqref{eq:intermediate-weight}. For $i\in\{1,2\}$ let
$\mathbb A_i=\{A_{i,Q}\}_{Q\in\D}$ be reducing operators of
order $p_i$ for $W_i$. Put
\begin{equation}\label{eq:intermediate-reducing}
 B_Q=(A_{1,Q}^2\#_\theta A_{2,Q}^2)^{1/2},
 \qquad \mathbb B=\{B_Q\}_{Q\in\D}.
\end{equation}

\smallskip
\noindent$\mathrm{(i)}$ The intermediate weight satisfies $W_\theta\in\mathcal A_{p_\theta}$.

\smallskip
\noindent$\mathrm{(ii)}$ The family $\mathbb B$ satisfies
\begin{equation}\label{eq:intermediate-reducing-average}
 |B_Q\vec v|\asymp
 \left(\fint_Q|W_\theta^{1/p_\theta}(x)\vec v|^{p_\theta}\,\mathrm d x
 \right)^{1/p_\theta},
 \qquad Q\in\D,\quad\vec v\in\C^m.
\end{equation}
The positive equivalence constants depend only on
$n,m,p_1,p_2,\theta,[W_1]_{\mathcal A_{p_1}}$ and
$[W_2]_{\mathcal A_{p_2}}$, and are independent of $Q$ and $\vec v$.
The family $\mathbb B$ also satisfies \eqref{eq:comparison}
with the exponents $\delta_1,\delta_2,\Delta$ in
\eqref{eq:dimension-exponents} for $W_\theta$ at order $p_\theta$.
\end{lemma}

\begin{proof}
As mentioned earlier, we only need to prove $\mathrm{(ii)}$.
Fix $Q\in\D$.
Using \cite[(8.6)--(8.11), proof of Proposition 8.8]{BCU},
Lemma~\ref{lem:ellipsoid-mean}, the reducing property and
Lemma~\ref{lem:Z-dual-reducing}, we obtain, for $\vec v\in\C^m$,
\[
 \begin{aligned}
 \left(\fint_Q|W_\theta^{1/p_\theta}(x)\vec v|^{p_\theta}
 \,\mathrm d x\right)^{1/p_\theta}
 &\lesssim|B_Q\vec v|,\\
 \left(\fint_Q|W_\theta^{-1/p_\theta}(x)\vec v|^{p_\theta'}
 \,\mathrm d x\right)^{1/p_\theta'}
 &\lesssim|B_Q^{-1}\vec v|.
 \end{aligned}
\]
By the second estimate and H\"older's inequality,
\[
 \begin{aligned}
 |B_Q\vec v|
 &=\sup_{\vec0\ne\vec u\in\C^m}
   \frac{|\langle\vec v,\vec u\rangle_{\C^m}|}{|B_Q^{-1}\vec u|}\\
 &\lesssim\sup_{\vec0\ne\vec u\in\C^m}
   \frac{|\langle\vec v,\vec u\rangle_{\C^m}|}
   {\left(\fint_Q|W_\theta^{-1/p_\theta}(x)\vec u|^{p_\theta'}
   \,\mathrm d x\right)^{1/p_\theta'}}
 \le\left(\fint_Q|W_\theta^{1/p_\theta}(x)\vec v|^{p_\theta}
 \,\mathrm d x\right)^{1/p_\theta}.
 \end{aligned}
\]
Comparing \eqref{eq:intermediate-reducing-average} with
\eqref{eq:reducing} and applying Lemma~\ref{lem:matrix-comparison},
we obtain \eqref{eq:comparison} for $\mathbb B$.
\end{proof}

\Needspace{6\baselineskip}
\subsection{Complex interpolation of \texorpdfstring{$Z$}{Z}-spaces and tent spaces}
\label{subsec:ZT-complex}

We use the following unweighted interpolation identities. Part (i)
follows from \cite[Theorem 4.3, pp.~923--924]{HuangTent} and
\cite[Lemma 3.3]{ABH}; part (ii) is the Banach-range
case of \cite[Theorem 5.1]{ABH}. The scalar results extend to
$\C^m$-valued functions by applying them to the finitely many coordinates.

\Needspace{4\baselineskip}
\begin{lemma}\label{lem:ZT-unweighted-complex}
Let $0<\theta<1$, $\alpha_i\in\R$, $1<p_i,q_i,r_i<\infty$
for $i\in\{1,2\}$, and let $\mathbb I=\{I_m\}_{Q\in\D}$.
Set $\alpha_\theta=(1-\theta)\alpha_1+\theta\alpha_2$ and
\[
 \frac1{p_\theta}=\frac{1-\theta}{p_1}+\frac\theta{p_2},\qquad
 \frac1{q_\theta}=\frac{1-\theta}{q_1}+\frac\theta{q_2},\qquad
 \frac1{r_\theta}=\frac{1-\theta}{r_1}+\frac\theta{r_2}.
\]
Then the following identities hold with equivalent norms.

\smallskip
\noindent$\mathrm{(i)}$ For tent spaces,
\[
 [T^{p_1,q_1,r_1}_{\alpha_1}(\mathbb I),
  T^{p_2,q_2,r_2}_{\alpha_2}(\mathbb I)]_\theta
 =T^{p_\theta,q_\theta,r_\theta}_{\alpha_\theta}(\mathbb I).
\]

\smallskip
\noindent$\mathrm{(ii)}$ For $Z$-spaces,
\[
 [Z^{p_1,q_1,r_1}_{\alpha_1}(\mathbb I),
  Z^{p_2,q_2,r_2}_{\alpha_2}(\mathbb I)]_\theta
 =Z^{p_\theta,q_\theta,r_\theta}_{\alpha_\theta}(\mathbb I).
\]
\end{lemma}

\Needspace{4\baselineskip}
\begin{theorem}\label{thm:Z-complex}\label{thm:T-complex}
Let $0<\theta<1$, $\alpha_i\in\R$, $1<p_i,q_i,r_i<\infty$,
and $W_i\in\mathcal A_{p_i}$ for $i\in\{1,2\}$.
Define $p_\theta,W_\theta$ by \eqref{eq:intermediate-weight}, and put
\begin{equation}\label{eq:interpolation-indices}
 \frac1{q_\theta}=\frac{1-\theta}{q_1}+\frac\theta{q_2},\qquad
 \frac1{r_\theta}=\frac{1-\theta}{r_1}+\frac\theta{r_2},\qquad
 \alpha_\theta=(1-\theta)\alpha_1+\theta\alpha_2.
\end{equation}
The following identities hold with equivalent norms.

\smallskip
\noindent$\mathrm{(i)}$ For tent spaces,
\begin{equation}\label{eq:T-complex}
 [T^{p_1,q_1,r_1}_{\alpha_1}(W_1),
  T^{p_2,q_2,r_2}_{\alpha_2}(W_2)]_\theta
 =T^{p_\theta,q_\theta,r_\theta}_{\alpha_\theta}(W_\theta).
\end{equation}

\smallskip
\noindent$\mathrm{(ii)}$ For $Z$-spaces,
\begin{equation}\label{eq:Z-complex}
 [Z^{p_1,q_1,r_1}_{\alpha_1}(W_1),
  Z^{p_2,q_2,r_2}_{\alpha_2}(W_2)]_\theta
 =Z^{p_\theta,q_\theta,r_\theta}_{\alpha_\theta}(W_\theta).
\end{equation}
\end{theorem}

\begin{proof}
\noindent$\mathrm{(i)}$ Let $\mathbb I=\{I_m\}_{Q\in\D}$, and let
$\mathbb A_i=\{A_{i,Q}\}_{Q\in\D}$ be sequences of reducing
operators of order $p_i$ for $W_i$, $i\in\{1,2\}$. Put
\[
 B_Q=(A_{1,Q}^2\#_\theta A_{2,Q}^2)^{1/2},\qquad
 \mathbb B=\{B_Q\}_{Q\in\D}.
\]
For each $Q\in\D$, choose an invertible matrix $S_Q$ and a
positive diagonal matrix
$D_Q=\operatorname{diag}(d_{Q,1},\ldots,d_{Q,m})$, with
$d_{Q,k}>0$ for $1\le k\le m$,
as in Lemma~\ref{lem:ellipsoid-mean}~(i). Thus, for $\vec v\in\C^m$,
\begin{equation}\label{eq:block-diagonalization}
 |A_{1,Q}\vec v|=|S_Q^{-1}\vec v|,\qquad
 |A_{2,Q}\vec v|=|D_QS_Q^{-1}\vec v|,\qquad
 |B_Q\vec v|=|D_Q^\theta S_Q^{-1}\vec v|.
\end{equation}
For $Q\in\D$, recall that $\widehat Q=[\ell(Q),2\ell(Q))\times Q$.
For $z\in\C$, define on $\widehat Q$
\[
 (\mathcal V_z\vec F)(s,y)=D_Q^zS_Q^{-1}\vec F(s,y),\qquad
 (\mathcal V_z^{-1}\vec H)(s,y)=S_QD_Q^{-z}\vec H(s,y),
\]
where $D_Q^z=\operatorname{diag}(e^{z\log_e d_{Q,1}},\ldots,
e^{z\log_e d_{Q,m}})$. Since $|D_Q^{it}\vec v|=|\vec v|$ for
$t\in\R$, using \eqref{eq:block-diagonalization} and
Lemmas~\ref{lem:Z-dyadic-A} and \ref{lem:Z-W-A}, we obtain
\begin{align*}
 \|\mathcal V_{it}\vec F\|_{T^{p_1,q_1,r_1}_{\alpha_1}(\mathbb I)}
 &\asymp\|\vec F\|_{T^{p_1,q_1,r_1}_{\alpha_1}(W_1)},\\
 \|\mathcal V_{1+it}\vec F\|_{T^{p_2,q_2,r_2}_{\alpha_2}(\mathbb I)}
 &\asymp\|\vec F\|_{T^{p_2,q_2,r_2}_{\alpha_2}(W_2)}.
\end{align*}
At $z=\theta$, by Lemma~\ref{lem:intermediate-reducing}, we have
\[
 \|\mathcal V_\theta\vec F\|_{T^{p_\theta,q_\theta,r_\theta}_{\alpha_\theta}(\mathbb I)}
 \asymp\|\vec F\|_{T^{p_\theta,q_\theta,r_\theta}_{\alpha_\theta}(\mathbb B)}
 \asymp\|\vec F\|_{T^{p_\theta,q_\theta,r_\theta}_{\alpha_\theta}(W_\theta)}.
\]

First suppose that $\vec F$ is supported in
$E_{\mathcal F}:=\bigcup_{Q\in\mathcal F}\widehat Q$, where
$\mathcal F\subset\D$ is finite. On this set, $\mathcal V_z$
and $\mathcal V_z^{-1}$ are bounded analytic matrix multipliers
for $0<\operatorname{Re}z<1$, uniformly bounded on the closed strip.
For
$(T^{p_1,q_1,r_1}_{\alpha_1}(W_1),
  T^{p_2,q_2,r_2}_{\alpha_2}(W_2))$ and
$(T^{p_1,q_1,r_1}_{\alpha_1}(\mathbb I),
  T^{p_2,q_2,r_2}_{\alpha_2}(\mathbb I))$,
respectively, apply the maps
\[
 \vec h(z)\longmapsto\mathcal V_z(\mathbf1_{E_{\mathcal F}}\vec h(z)),
 \qquad
 \vec h(z)\longmapsto\mathcal V_z^{-1}(\mathbf1_{E_{\mathcal F}}\vec h(z))
\]
to admissible analytic functions as in
Definition~\ref{def:interpolation-spaces}~(ii). Using the boundary
estimates and Lemma~\ref{lem:ZT-unweighted-complex}~(i), we obtain
\begin{align*}
 &\|\vec F\|_{[T^{p_1,q_1,r_1}_{\alpha_1}(W_1),
                T^{p_2,q_2,r_2}_{\alpha_2}(W_2)]_\theta}
 \asymp\|\mathcal V_\theta\vec F\|
 _{[T^{p_1,q_1,r_1}_{\alpha_1}(\mathbb I),
    T^{p_2,q_2,r_2}_{\alpha_2}(\mathbb I)]_\theta}\\
 &\quad\asymp\|\mathcal V_\theta\vec F\|
 _{T^{p_\theta,q_\theta,r_\theta}_{\alpha_\theta}(\mathbb I)}
 \asymp\|\vec F\|
 _{T^{p_\theta,q_\theta,r_\theta}_{\alpha_\theta}(W_\theta)}.
\end{align*}
The constants are independent of $\mathcal F$.
Choose finite sets $\mathcal F_N\uparrow\D$, $N\in\mathbb N$.
For $i\in\{1,2,\theta\}$ and
$\vec F\in T^{p_i,q_i,r_i}_{\alpha_i}(W_i)$, using the dyadic
characterization and dominated convergence, we obtain
\[
 \|\mathbf1_{E_{\mathcal F_N}^c}\vec F\|
 _{T^{p_i,q_i,r_i}_{\alpha_i}(W_i)}\longrightarrow0
 \qquad\text{as }N\to\infty.
\]
For $\vec H\in T^{p_1,q_1,r_1}_{\alpha_1}(W_1)
\cap T^{p_2,q_2,r_2}_{\alpha_2}(W_2)$,
using \cite[(C.1)]{HNVW16}, we have
\[
 \|\mathbf1_{E_{\mathcal F_N}^c}\vec H\|
 _{[T^{p_1,q_1,r_1}_{\alpha_1}(W_1),
    T^{p_2,q_2,r_2}_{\alpha_2}(W_2)]_\theta}
 \le\|\mathbf1_{E_{\mathcal F_N}^c}\vec H\|
 _{T^{p_1,q_1,r_1}_{\alpha_1}(W_1)}^{1-\theta}
 \|\mathbf1_{E_{\mathcal F_N}^c}\vec H\|
 _{T^{p_2,q_2,r_2}_{\alpha_2}(W_2)}^\theta\longrightarrow0.
\]
Multiplication by $\mathbf1_{E_{\mathcal F_N}^c}$ is contractive
on $T^{p_i,q_i,r_i}_{\alpha_i}(W_i)$ for $i\in\{1,2\}$ and hence on
\[
 [T^{p_1,q_1,r_1}_{\alpha_1}(W_1),
  T^{p_2,q_2,r_2}_{\alpha_2}(W_2)]_\theta
\]
by \cite[Theorem C.2.6]{HNVW16}. Since
\[
 T^{p_1,q_1,r_1}_{\alpha_1}(W_1)\cap
 T^{p_2,q_2,r_2}_{\alpha_2}(W_2)
 \quad\text{is dense in}\quad
 [T^{p_1,q_1,r_1}_{\alpha_1}(W_1),
  T^{p_2,q_2,r_2}_{\alpha_2}(W_2)]_\theta
\]
by
\cite[Corollary C.2.8]{HNVW16}, for every
$\vec F\in[T^{p_1,q_1,r_1}_{\alpha_1}(W_1),
T^{p_2,q_2,r_2}_{\alpha_2}(W_2)]_\theta$, we obtain
\[
 \|\vec F-\mathbf1_{E_{\mathcal F_N}}\vec F\|
 _{[T^{p_1,q_1,r_1}_{\alpha_1}(W_1),
    T^{p_2,q_2,r_2}_{\alpha_2}(W_2)]_\theta}
 \longrightarrow0\qquad\text{as }N\to\infty.
\]
Together with convergence in
$T^{p_\theta,q_\theta,r_\theta}_{\alpha_\theta}(W_\theta)$,
this extends the norm equivalence by completion.

\medskip
\noindent$\mathrm{(ii)}$ The proof is similar, with the $Z$-space
norms and Lemma~\ref{lem:ZT-unweighted-complex}~(ii), and is omitted.
\end{proof}

\subsection{Real interpolation of \texorpdfstring{$Z$}{Z}-spaces and tent spaces}
\label{subsec:ZT-real}

\Needspace{7\baselineskip}
For the real interpolation result with different weights, we use
the following consequence of \cite[Proposition 5.1]{LL}, with
counting measure and the same Banach space at both endpoints.

\begin{lemma}\label{lem:weighted-sequence-real}
Let $E$ be a complex Banach space, $\mathcal J$ a countable set,
$1<p_1,p_2<\infty$ and $0<\theta<1$. Put
$1/p_\theta=(1-\theta)/p_1+\theta/p_2$.
Let $a_{i,\nu}>0$ for $i\in\{1,2\}$ and $\nu\in\mathcal J$.
For an $E$-valued sequence $h=\{h_\nu\}_{\nu\in\mathcal J}$ and $t>0$, set
\[
 K(t,h):=\inf_{h=h^{(1)}+h^{(2)}}\left[
 \left(\sum_{\nu\in\mathcal J}
 a_{1,\nu}^{p_1}\|h^{(1)}_\nu\|_E^{p_1}\right)^{1/p_1}
 +t\left(\sum_{\nu\in\mathcal J}
 a_{2,\nu}^{p_2}\|h^{(2)}_\nu\|_E^{p_2}\right)^{1/p_2}\right].
\]
Then
\[
 \left(\int_0^\infty[t^{-\theta}K(t,h)]^{p_\theta}
 \frac{\mathrm d t}{t}\right)^{1/p_\theta}
 \asymp\left(\sum_{\nu\in\mathcal J}
 [a_{1,\nu}^{1-\theta}a_{2,\nu}^{\theta}\|h_\nu\|_E]^{p_\theta}
 \right)^{1/p_\theta}.
\]
\end{lemma}

\Needspace{4\baselineskip}
\begin{theorem}\label{thm:Z-real}\label{cor:T-real-diagonal}
Let $1<p_1,p_2<\infty$ with $p_1\ne p_2$, let $1<r<\infty$,
$\alpha_1,\alpha_2\in\R$, and $0<\theta<1$.
Let $W_i\in\mathcal A_{p_i}$ for $i\in\{1,2\}$, and define
$p_\theta,W_\theta$ by \eqref{eq:intermediate-weight} and
$\alpha_\theta=(1-\theta)\alpha_1+\theta\alpha_2$.

\smallskip
\noindent$\mathrm{(i)}$ For $Z$-spaces,
\begin{equation}\label{eq:Z-real}
 (Z^{p_1,p_1,r}_{\alpha_1}(W_1),
  Z^{p_2,p_2,r}_{\alpha_2}(W_2))_{\theta,p_\theta}
 =Z^{p_\theta,p_\theta,r}_{\alpha_\theta}(W_\theta).
\end{equation}

\smallskip
\noindent$\mathrm{(ii)}$ For tent spaces,
\begin{equation}\label{eq:T-real-diagonal}
 (T^{p_1,p_1,r}_{\alpha_1}(W_1), T^{p_2,p_2,r}_{\alpha_2}(W_2))_{\theta,p_\theta} =T^{p_\theta,p_\theta,r}_{\alpha_\theta}(W_\theta) =Z^{p_\theta,p_\theta,r}_{\alpha_\theta}(W_\theta).
\end{equation}
Both identities hold with equivalent norms.
\end{theorem}

\begin{proof}
\noindent$\mathrm{(i)}$ Let $\mathbb A_i=\{A_{i,Q}\}_{Q\in\D}$ be sequences
of reducing operators of order $p_i$ for $W_i$, $i\in\{1,2\}$.
Put $B_Q=(A_{1,Q}^2\#_\theta A_{2,Q}^2)^{1/2}$ and
$\mathbb B=\{B_Q\}_{Q\in\D}$.
For each $Q\in\D$, choose an invertible matrix $S_Q$ and a
positive diagonal matrix $D_Q=\operatorname{diag}(d_{Q,1},\ldots,d_{Q,m})$
as in Lemma~\ref{lem:ellipsoid-mean}~(i), so that
\eqref{eq:block-diagonalization} holds. Recall that
\[
 \widehat Q=[\ell(Q),2\ell(Q))\times Q,\qquad
 \mathrm d\nu_Q(s,y)=\frac{\mathrm d y\,\mathrm d s}{s|Q|\log_e2}
 \quad\text{on }\widehat Q.
\]
Let
\[
 \Omega=[1,2)\times[0,1)^n,\qquad
 \mathrm d\mu(\sigma,z)=\frac{\mathrm d\sigma\,\mathrm d z}{\sigma\log_e2}.
\]
For $j\in\Z$, $Q=2^{-j}(z_Q+[0,1)^n)\in\D_j$, $z_Q\in\Z^n$,
$(\sigma,z)\in\Omega$ and $1\le k\le m$, define
\begin{equation}\label{eq:real-coordinate-change}
 g_{Q,k}(\sigma,z)
 :=[S_Q^{-1}\vec F(2^{-j}\sigma,2^{-j}(z_Q+z))]_k,
\end{equation}
where $[\cdot]_k$ denotes the $k$th coordinate.
This is a common linear change of coordinates, with inverse
\[
 \vec F(2^{-j}\sigma,2^{-j}(z_Q+z))
 =S_Q(g_{Q,1}(\sigma,z),\ldots,g_{Q,m}(\sigma,z))^{\mathsf T}.
\]
Using \eqref{eq:block-diagonalization}, the change of variables
$(s,y)=(2^{-j}\sigma,2^{-j}(z_Q+z))$ and
Lemmas~\ref{lem:Z-dyadic-A} and \ref{lem:Z-W-A}, we obtain
\begin{align*}
 \|\vec F\|_{Z^{p_1,p_1,r}_{\alpha_1}(W_1)}
 &\asymp\left\{\sum_{j\in\Z}\sum_{Q\in\D_j}\sum_{k=1}^m
 [2^{j\alpha_1}|Q|^{1/p_1}\|g_{Q,k}\|_{L^r(\Omega,\mu)}]^{p_1}
 \right\}^{1/p_1},\\
 \|\vec F\|_{Z^{p_2,p_2,r}_{\alpha_2}(W_2)}
 &\asymp\left\{\sum_{j\in\Z}\sum_{Q\in\D_j}\sum_{k=1}^m
 [2^{j\alpha_2}|Q|^{1/p_2}d_{Q,k}\|g_{Q,k}\|_{L^r(\Omega,\mu)}]^{p_2}
 \right\}^{1/p_2}.
\end{align*}
Apply Lemma~\ref{lem:weighted-sequence-real} with
$\mathcal J=\D\times\{1,\ldots,m\}$, $E=L^r(\Omega,\mu)$ and
\[
 a_{1,(Q,k)}=2^{j\alpha_1}|Q|^{1/p_1},\quad
 a_{2,(Q,k)}=2^{j\alpha_2}|Q|^{1/p_2}d_{Q,k},\quad Q\in\D_j.
\]
By the change of coordinates in \eqref{eq:real-coordinate-change}, we have
\[
 K(t,\vec F;Z^{p_1,p_1,r}_{\alpha_1}(W_1),
                Z^{p_2,p_2,r}_{\alpha_2}(W_2))
 \asymp K(t,\{g_{Q,k}\}_{(Q,k)\in\mathcal J}),\qquad t>0.
\]
Since
\[
 (2^{j\alpha_1}|Q|^{1/p_1})^{1-\theta}
 (2^{j\alpha_2}|Q|^{1/p_2}d_{Q,k})^\theta
 =2^{j\alpha_\theta}|Q|^{1/p_\theta}d_{Q,k}^\theta,
\]
Using Lemma~\ref{lem:intermediate-reducing}, we obtain
\begin{align*}
 &\|\vec F\|_{(Z^{p_1,p_1,r}_{\alpha_1}(W_1),
                Z^{p_2,p_2,r}_{\alpha_2}(W_2))_{\theta,p_\theta}}
 \asymp\left\{\sum_{j\in\Z}\sum_{Q\in\D_j}\sum_{k=1}^m
 [2^{j\alpha_\theta}|Q|^{1/p_\theta}d_{Q,k}^\theta
 \|g_{Q,k}\|_{L^r(\Omega,\mu)}]^{p_\theta}\right\}^{1/p_\theta}\\
 &\quad\asymp\|\vec F\|_{Z^{p_\theta,p_\theta,r}_{\alpha_\theta}(W_\theta)}.
\end{align*}

\medskip
\noindent$\mathrm{(ii)}$ Let $1<p,r<\infty$, $\alpha\in\R$ and
$W\in\mathcal A_p$. By Tonelli's theorem,
\begin{align*}
 \|\vec F\|_{T^{p,p,r}_\alpha(W)}^p
 &=\int_{\R^n}\int_0^\infty
 \left[\fint_{t/2}^{t}\fint_{B(x,t)}
 |s^{-\alpha}W^{1/p}(x)\vec F(s,y)|^r\,\mathrm d y\,\mathrm d s
 \right]^{p/r}\frac{\mathrm d t}{t}\,\mathrm d x\\
 &=\int_0^\infty\int_{\R^n}
 \left[\fint_{t/2}^{t}\fint_{B(x,t)}
 |s^{-\alpha}W^{1/p}(x)\vec F(s,y)|^r\,\mathrm d y\,\mathrm d s
 \right]^{p/r}\mathrm d x\,\frac{\mathrm d t}{t} =\|\vec F\|_{Z^{p,p,r}_\alpha(W)}^p.
\end{align*}
Using (i), we therefore obtain, with equivalent norms,
\begin{align*}
 &(T^{p_1,p_1,r}_{\alpha_1}(W_1),
   T^{p_2,p_2,r}_{\alpha_2}(W_2))_{\theta,p_\theta} =(Z^{p_1,p_1,r}_{\alpha_1}(W_1),
         Z^{p_2,p_2,r}_{\alpha_2}(W_2))_{\theta,p_\theta}\\
 &\quad=Z^{p_\theta,p_\theta,r}_{\alpha_\theta}(W_\theta)
       =T^{p_\theta,p_\theta,r}_{\alpha_\theta}(W_\theta).\qedhere
\end{align*}
\end{proof}

\begin{remark}\label{rem:interpolation-special-cases}
If $W_1=W_2=W\in\mathcal A_{\min\{p_1,p_2\}}$, both endpoint
assumptions hold by Remark~\ref{rem:Ap-monotonicity}, and
by \eqref{eq:intermediate-weight}, we have $W_\theta=W$.
Using Theorems~\ref{thm:Z-complex} and \ref{thm:Z-real},
we also obtain interpolation with a fixed matrix weight and varying exponents.
\end{remark}

The following unweighted identities are
\cite[Theorem 4.1 and Proposition 4.4]{ABH}, restricted to the
range used below. They extend to $\C^m$-valued functions coordinatewise.

\Needspace{4\baselineskip}
\begin{lemma}\label{lem:ZT-unweighted-real}
Let $1<p,r,q_1,q_2,u<\infty$, $\alpha_1,\alpha_2\in\R$
with $\alpha_1\ne\alpha_2$, $0<\theta<1$, and
$\mathbb I=\{I_m\}_{Q\in\D}$.
Set $\alpha_\theta=(1-\theta)\alpha_1+\theta\alpha_2$.
Then the following identities hold with equivalent norms.

\smallskip
\noindent$\mathrm{(i)}$ For $Z$-spaces,
\[
 (Z^{p,q_1,r}_{\alpha_1}(\mathbb I),
  Z^{p,q_2,r}_{\alpha_2}(\mathbb I))_{\theta,u}
 =Z^{p,u,r}_{\alpha_\theta}(\mathbb I).
\]

\smallskip
\noindent$\mathrm{(ii)}$ For tent spaces,
\[
 (T^{p,q_1,r}_{\alpha_1}(\mathbb I),
  T^{p,q_2,r}_{\alpha_2}(\mathbb I))_{\theta,u}
 =Z^{p,u,r}_{\alpha_\theta}(\mathbb I).
\]
\end{lemma}

\Needspace{4\baselineskip}
\begin{theorem}\label{lem:Z-real-smoothness}\label{thm:T-real-smoothness}
Let $1<p,r,q_1,q_2,u<\infty$, $W\in\mathcal A_p$,
$\alpha_1,\alpha_2\in\R$ with $\alpha_1\ne\alpha_2$, and
$0<\theta<1$. Set
$\alpha_\theta=(1-\theta)\alpha_1+\theta\alpha_2$.

Then the following identities hold with equivalent norms.

\smallskip
\noindent$\mathrm{(i)}$ For $Z$-spaces,
\begin{equation}\label{eq:Z-real-smoothness}
 (Z^{p,q_1,r}_{\alpha_1}(W),
  Z^{p,q_2,r}_{\alpha_2}(W))_{\theta,u}
 =Z^{p,u,r}_{\alpha_\theta}(W).
\end{equation}

\smallskip
\noindent$\mathrm{(ii)}$ For tent spaces,
\begin{equation}\label{eq:T-real-smoothness}
 (T^{p,q_1,r}_{\alpha_1}(W),
  T^{p,q_2,r}_{\alpha_2}(W))_{\theta,u}
 =Z^{p,u,r}_{\alpha_\theta}(W).
\end{equation}
\end{theorem}

\begin{proof}
Let $\mathbb I=\{I_m\}_{Q\in\D}$, and let
$\mathbb A=\{A_Q\}_{Q\in\D}$ be a sequence of reducing
operators of order $p$ for $W$. For $Q\in\D$, recall that
$\widehat Q=[\ell(Q),2\ell(Q))\times Q$. Define
\[
 (\mathcal U\vec F)(s,y)=A_Q\vec F(s,y),\qquad
 (\mathcal U^{-1}\vec H)(s,y)=A_Q^{-1}\vec H(s,y),
 \quad(s,y)\in\widehat Q.
\]
For $\alpha\in\R$ and $1<q<\infty$, using
Lemmas~\ref{lem:Z-dyadic-A} and \ref{lem:Z-W-A}, we obtain
\begin{align*}
 \|\mathcal U\vec F\|_{Z^{p,q,r}_\alpha(\mathbb I)}
 &\asymp\|\vec F\|_{Z^{p,q,r}_\alpha(W)},\\
 \|\mathcal U\vec F\|_{T^{p,q,r}_\alpha(\mathbb I)}
 &\asymp\|\vec F\|_{T^{p,q,r}_\alpha(W)}.
\end{align*}

\smallskip
\noindent$\mathrm{(i)}$ Using Lemma~\ref{lem:ZT-unweighted-real}~(i), we obtain
\begin{align*}
 &\|\vec F\|_{(Z^{p,q_1,r}_{\alpha_1}(W),
                Z^{p,q_2,r}_{\alpha_2}(W))_{\theta,u}}
 \asymp\|\mathcal U\vec F\|
 _{(Z^{p,q_1,r}_{\alpha_1}(\mathbb I),
    Z^{p,q_2,r}_{\alpha_2}(\mathbb I))_{\theta,u}}\\
 &\quad\asymp\|\mathcal U\vec F\|_{Z^{p,u,r}_{\alpha_\theta}(\mathbb I)}
 \asymp\|\vec F\|_{Z^{p,u,r}_{\alpha_\theta}(W)}.
\end{align*}

\medskip
\noindent$\mathrm{(ii)}$ Using Lemma~\ref{lem:ZT-unweighted-real}~(ii), we obtain
\begin{align*}
 &\|\vec F\|_{(T^{p,q_1,r}_{\alpha_1}(W),
                T^{p,q_2,r}_{\alpha_2}(W))_{\theta,u}}
 \asymp\|\mathcal U\vec F\|
 _{(T^{p,q_1,r}_{\alpha_1}(\mathbb I),
    T^{p,q_2,r}_{\alpha_2}(\mathbb I))_{\theta,u}}\\
 &\quad\asymp\|\mathcal U\vec F\|_{Z^{p,u,r}_{\alpha_\theta}(\mathbb I)}
 \asymp\|\vec F\|_{Z^{p,u,r}_{\alpha_\theta}(W)}.\qedhere
\end{align*}
\end{proof}

\Needspace{4\baselineskip}
\subsection{Consequences for Besov and Triebel--Lizorkin spaces}
\label{subsec:BT-interpolation}

In this subsection, we study complex and real interpolation
of homogeneous matrix-weighted Besov and
Triebel--Lizorkin spaces.

Recall that $\cS_\infty$ consists of Schwartz functions with
all moments equal to zero. We use the canonical identification
$\cS_\infty'\cong\cS'/\cP$; see \cite[Section 1]{FR}.

\Needspace{4\baselineskip}
Fix $\varphi\in\cS$ satisfying \eqref{eq:partition}, with
$\varphi_j(x)=2^{jn}\varphi(2^jx)$ for $j\in\Z$.
Choose $\Phi\in\cS$ such that
\begin{equation}\label{eq:Phi}
 \supp\widehat\Phi\subset\{\xi:1/4<|\xi|<4\},\qquad
 \widehat\Phi(\xi)=1\quad\text{if }3/8\le|\xi|\le5/2,
\end{equation}
and set $\Phi_s(x)=s^{-n}\Phi(x/s)$. For every
$j\in\Z$ and $\vec f\in(\cS'/\cP)^m$,
\begin{equation}\label{eq:reproducing}
 \varphi_j*\Phi_s*\vec f=\varphi_j*\vec f,
 \qquad 9\cdot2^{-j}/10<s<11\cdot2^{-j}/10.
\end{equation}
Indeed, for these $s$ and $\xi\in\supp\widehat\varphi_j$,
$9/20\le|s\xi|\le11/5$, so by \eqref{eq:Phi}, we have
$\widehat\Phi(s\xi)=1$.
Choose a nonnegative function
$\rho\in C_c^\infty((9/10,11/10))$ such that
$\int_0^\infty\rho(u)\,\mathrm d u/u=1$. Define
\begin{equation}\label{eq:BT-extension-synthesis}
 (\mathcal E\vec f)(s,y):=(\Phi_s*\vec f)(y),\qquad \mathcal R_\rho\vec F :=\sum_{j\in\Z}\varphi_j*\left[\int_0^\infty\rho(2^js)\vec F(s,\cdot)\,\frac{\mathrm d s}{s}\right].
\end{equation}
\Needspace{5\baselineskip}
The next lemma proves boundedness of $\mathcal E$ and
$\mathcal R_\rho$ and the identity $\mathcal R_\rho\mathcal E=I$.

\begin{lemma}\label{lem:BT-retraction}
Let $\alpha\in\R$, $1<p,q,r<\infty$, and $W\in\mathcal A_p$.
Let $\varphi,\Phi\in\cS$ satisfy \eqref{eq:partition} and
\eqref{eq:Phi}, respectively, and let
$0\le\rho\in C_c^\infty((9/10,11/10))$ satisfy
$\int_0^\infty\rho(u)\,\mathrm d u/u=1$.

\smallskip
\noindent$\mathrm{(i)}$ The operators $\mathcal E,\mathcal R_\rho$ in
\eqref{eq:BT-extension-synthesis} are bounded
in both directions in the following diagrams:
\begin{equation}\label{eq:BT-retraction-bounds}
 \dot B^\alpha_{p,q}(W)
 \mathrel{\mathop{\rightleftarrows}^{\mathcal E}_{\mathcal R_\rho}}
 Z^{p,q,r}_\alpha(W),
 \qquad
 \dot F^\alpha_{p,q}(W)
 \mathrel{\mathop{\rightleftarrows}^{\mathcal E}_{\mathcal R_\rho}}
 T^{p,q,r}_\alpha(W).
\end{equation}

\smallskip
\noindent$\mathrm{(ii)}$ For every $\vec f\in\dot B^\alpha_{p,q}(W)\cup\dot F^\alpha_{p,q}(W)$,
\begin{equation}\label{eq:BT-retraction-identity}
 \mathcal R_\rho\mathcal E\vec f=\vec f
 \qquad\text{in }(\cS'/\cP)^m.
\end{equation}
\end{lemma}

\begin{proof}
\noindent$\mathrm{(i)}$ Boundedness of $\mathcal E$ in \eqref{eq:BT-retraction-bounds}
follows from Theorem~\ref{thm:general-Whitney}, since the kernel
in \eqref{eq:Phi} has all moments equal to zero.

\medskip
\noindent\emph{Step 1.} First, we prove that $\mathcal R_\rho$ is well defined.
Let $\mathbb A=\{A_Q\}_{Q\in\D}$ be a sequence of reducing
operators of order $p$ for $W$. For $l\in\Z$ and $Q\in\D_l$, put
\[
 \widehat Q=[2^{-l},2^{-l+1})\times Q,\qquad
 \mathrm d\nu_Q(s,y)=\frac{\mathrm d y\,\mathrm d s}{s|Q|\log_e2}.
\]
The dyadic characterizations give
\begin{equation}\label{eq:BT-block-control}
 \sup_{l\in\Z}2^{l\alpha}
 \left(\sum_{Q\in\D_l}|Q|\|A_Q\vec F\|_{L^r(\widehat Q,\nu_Q;\C^m)}^p\right)^{1/p}
 \lesssim
 \begin{cases}
 \|\vec F\|_{Z^{p,q,r}_\alpha(W)},&
       \vec F\in Z^{p,q,r}_\alpha(W),\\
 \|\vec F\|_{T^{p,q,r}_\alpha(W)},&
       \vec F\in T^{p,q,r}_\alpha(W).
 \end{cases}
\end{equation}
For $j\in\Z$, since $\operatorname{supp}\rho(2^j\cdot)
\subset[2^{-j-1},2^{-j+1})$, only $l\in\{j,j+1\}$ occur.
H\"older's inequality yields
\begin{equation}\label{eq:BT-synthesis-holder}
 \int_{2^{-l}}^{2^{-l+1}}\int_Q
 |A_Q\vec F(s,y)|\,\mathrm d y\,\frac{\mathrm d s}{s}
 \le |Q|\log_e2\,\|A_Q\vec F\|_{L^r(\widehat Q,\nu_Q;\C^m)}.
\end{equation}

Let $\delta_1,\delta_2,\Delta$ be as in
\eqref{eq:dimension-exponents}, and choose an integer
$M>\Delta+n/p'$. By \eqref{eq:comparison},
\[
 \|\mathbb A_l(y)^{-1}\|
 \lesssim2^{\delta_1(-l)_++\delta_2l_+}(1+|y|)^\Delta.
\]
For $y,z\in Q\in\D_l$,
$(1+|y|)^{\Delta-M}
 \lesssim2^{M|l|}(1+|z|)^{\Delta-M}$.
Thus, by \eqref{eq:BT-synthesis-holder} and H\"older's inequality,
\begin{align*}
 &\int_{\R^n}\left|\int_0^\infty\rho(2^js)\vec F(s,y)\,\frac{\mathrm d s}{s}\right|(1+|y|)^{-M}\,\mathrm d y\\
 &\quad\lesssim\sum_{l\in\{j,j+1\}}
 2^{\delta_1(-l)_++\delta_2l_++M|l|}
 \sum_{Q\in\D_l}\|A_Q\vec F\|_{L^r(\widehat Q,\nu_Q;\C^m)}
       \int_Q(1+|z|)^{\Delta-M}\,\mathrm d z,\\
 &\sum_{Q\in\D_l}\|A_Q\vec F\|_{L^r(\widehat Q,\nu_Q;\C^m)}\int_Q(1+|z|)^{\Delta-M}\,\mathrm d z\\
 &\quad\le\left(\sum_{Q\in\D_l}|Q|\|A_Q\vec F\|_{L^r(\widehat Q,\nu_Q;\C^m)}^p\right)^{1/p}
 \|(1+|\cdot|)^{\Delta-M}\|_{L^{p'}(\R^n)}.
\end{align*}
Put $C_0=M+\max\{\delta_1,\delta_2\}+|\alpha|+1$.
Together with \eqref{eq:BT-block-control}, these estimates yield
\begin{equation}\label{eq:BT-tempered-time}
 \begin{aligned}
 &\int_{\R^n}\left|\int_0^\infty\rho(2^js)\vec F(s,y)\,\frac{\mathrm d s}{s}\right|(1+|y|)^{-M}\,\mathrm d y\\
 &\quad\lesssim2^{C_0|j|}
 \sup_{l\in\Z}2^{l\alpha}
 \left(\sum_{Q\in\D_l}|Q|\|A_Q\vec F\|_{L^r(\widehat Q,\nu_Q;\C^m)}^p\right)^{1/p},
 \end{aligned}
\end{equation}
For every integer $L>C_0$, the test-function estimates
\eqref{eq:general-test-fourier} and \eqref{eq:general-test-decay} give
\begin{equation}\label{eq:BT-test-decay}
 \sup_y(1+|y|)^M
 |(\overline{\varphi_j(-\cdot)}*\vec\psi)(y)|
 \lesssim2^{-L|j|}\rho_{2M+L+n+2,\,L+n+1}(\vec\psi),
 \quad \vec\psi\in(\cS_\infty)^m.
\end{equation}
Consequently,
\begin{align*}
 &\sum_{j\in\Z}\left|\left\langle\varphi_j*\left[\int_0^\infty\rho(2^js)\vec F(s,\cdot)\,\frac{\mathrm d s}{s}\right],\vec\psi\right\rangle\right|\\
 &\quad\lesssim\rho_{2M+L+n+2,\,L+n+1}(\vec\psi)
 \sup_{l\in\Z}2^{l\alpha}
 \left(\sum_{Q\in\D_l}|Q|\|A_Q\vec F\|_{L^r(\widehat Q,\nu_Q;\C^m)}^p\right)^{1/p}
 \sum_{j\in\Z}2^{-(L-C_0)|j|}<\infty.
\end{align*}
Hence $\mathcal R_\rho\vec F$ in \eqref{eq:BT-extension-synthesis}
is well defined in $(\cS_\infty')^m$ for
$\vec F\in Z^{p,q,r}_\alpha(W)\cup T^{p,q,r}_\alpha(W)$.

\medskip
\noindent\emph{Step 2.} Then we prove the boundedness of $\mathcal R_\rho$.
By the Fourier support of $\varphi$, for $k\in\Z$,
\begin{equation}\label{eq:BT-synthesis-localization}
 \varphi_k*\mathcal R_\rho\vec F
 =\sum_{|j-k|\leq2}
       \varphi_k*\varphi_j*\left[\int_0^\infty\rho(2^js)\vec F(s,\cdot)\,\frac{\mathrm d s}{s}\right].
\end{equation}
Choose $N>n$. By the Schwartz kernel estimate, we have
\begin{equation}\label{eq:BT-synthesis-kernel}
 |(\varphi_k*\varphi_j)(x-y)|
 \lesssim2^{kn}(1+2^k|x-y|)^{-N-\Delta},
 \qquad |j-k|\leq2.
\end{equation}
On $[2^{-l},2^{-l+1})\times Q$, $Q\in\D_l$, \eqref{eq:comparison} implies
\begin{equation}\label{eq:BT-synthesis-matrix}
 |\mathbb A_k(x)\vec F(s,y)|
 \lesssim(1+2^k|x-y|)^\Delta|A_Q\vec F(s,y)|,
 \qquad y\in Q,\quad |k-l|\leq3.
\end{equation}
For $|k-l|\leq3$,
\[
 1+2^k|x-y|\asymp1+2^k|x-z|,\qquad y,z\in Q.
\]
Using \eqref{eq:BT-synthesis-localization}, \eqref{eq:BT-synthesis-kernel},
\eqref{eq:BT-synthesis-matrix}, \eqref{eq:BT-synthesis-holder}, and
\eqref{eq:decay-max} with $N>n$, we obtain

\begin{align*}
 &2^{k\alpha}|\mathbb A_k(x) (\varphi_k*\mathcal R_\rho\vec F)(x)| \lesssim \sum_{|l-k|\leq3}\int_{\R^n} \frac{2^{kn}}{(1+2^k|x-z|)^N} \sum_{Q\in\D_l}\mathbf1_Q(z)2^{l\alpha} \|A_Q\vec F\|_{L^r(\widehat Q,\nu_Q;\C^m)}\,\mathrm d z\\
 &\quad \lesssim\sum_{|l-k|\leq3}\M\left( \sum_{Q\in\D_l}\mathbf1_Q\,2^{l\alpha} \|A_Q\vec F\|_{L^r(\widehat Q,\nu_Q;\C^m)}\right)(x),
 \end{align*}
Using Lemma~\ref{lem:FS} with $p,q>1$,
\eqref{eq:known-equivalence}, and \eqref{eq:T-dyadic-A}, we obtain
\[
 \|\mathcal R_\rho\vec F\|_{\dot F^\alpha_{p,q}(W)} \lesssim\left\|\left\{\sum_{l\in\Z}\sum_{Q\in\D_l}\mathbf1_Q [2^{l\alpha}\|A_Q\vec F\|_{L^r(\widehat Q,\nu_Q;\C^m)}]^q \right\}^{1/q}\right\|_{L^p(\R^n)} \asymp\|\vec F\|_{T^{p,q,r}_\alpha(W)}.
\]
Taking the $\ell^q(L^p)$ norm in the same maximal-function estimate,
using boundedness of $\mathcal M$ on $L^p(\R^n)$, gives
\[
 \|\mathcal R_\rho\vec F\|_{\dot B^\alpha_{p,q}(W)}
 \lesssim
 \left\{\sum_{l\in\Z}
 \left[2^{l\alpha}
 \left(\sum_{Q\in\D_l}|Q|\|A_Q\vec F\|_{L^r(\widehat Q,\nu_Q;\C^m)}^p\right)^{1/p}\right]^q
 \right\}^{1/q}
 \asymp\|\vec F\|_{Z^{p,q,r}_\alpha(W)}.
\]

\medskip
\noindent$\mathrm{(ii)}$
For $\vec f\in\dot B^\alpha_{p,q}(W)\cup\dot F^\alpha_{p,q}(W)$,
\eqref{eq:reproducing} and the normalization of $\rho$ give
\[
 \varphi_j*\left[\int_0^\infty
 \rho(2^js)(\Phi_s*\vec f)\,\frac{\mathrm d s}{s}\right]
 =\int_0^\infty\rho(2^js)\varphi_j*\vec f\,\frac{\mathrm d s}{s}
 =\varphi_j*\vec f.
\]
By \eqref{eq:partition} and the argument of
\eqref{eq:BT-test-decay}, applied also to derivatives of $\vec\psi$,
$\sum_{j\in\Z}\overline{\varphi_j(-\cdot)}*\vec\psi=\vec\psi$
in $(\cS_\infty)^m$ for every $\vec\psi\in(\cS_\infty)^m$.
Therefore,
\[
 \mathcal R_\rho\mathcal E\vec f
 =\sum_{j\in\Z}\varphi_j*\vec f=\vec f
 \quad\text{in }(\cS_\infty')^m\cong(\cS'/\cP)^m.\qedhere
\]
\end{proof}

\Needspace{4\baselineskip}
\begin{corollary}\label{cor:BT-complex}
Let $0<\theta<1$, $\alpha_i\in\R$, $1<p_i,q_i<\infty$,
and $W_i\in\mathcal A_{p_i}$ for $i\in\{1,2\}$.
Let $p_\theta,W_\theta$ be as in
\eqref{eq:intermediate-weight}, and put
\[
 \frac1{q_\theta}=\frac{1-\theta}{q_1}+\frac\theta{q_2},
 \qquad \alpha_\theta=(1-\theta)\alpha_1+\theta\alpha_2.
\]

Then the following identities hold with equivalent norms.

\smallskip
\noindent$\mathrm{(i)}$ For the Triebel--Lizorkin spaces,
\begin{equation}
 [\dot F^{\alpha_1}_{p_1,q_1}(W_1),
  \dot F^{\alpha_2}_{p_2,q_2}(W_2)]_\theta
 =\dot F^{\alpha_\theta}_{p_\theta,q_\theta}(W_\theta),
 \label{eq:F-complex}
\end{equation}

\smallskip
\noindent$\mathrm{(ii)}$ For the Besov spaces,
\begin{equation}
 [\dot B^{\alpha_1}_{p_1,q_1}(W_1),
  \dot B^{\alpha_2}_{p_2,q_2}(W_2)]_\theta
 =\dot B^{\alpha_\theta}_{p_\theta,q_\theta}(W_\theta).
 \label{eq:B-complex}
\end{equation}
\end{corollary}

\begin{proof}
Take $r=2$ in Lemma~\ref{lem:BT-retraction}.
Both identities follow from Theorem~\ref{thm:Z-complex}, using
\cite[Theorem C.2.6]{HNVW16}.
\end{proof}

\Needspace{4\baselineskip}
\begin{corollary}\label{cor:BT-real}
Let $0<\theta<1$, $\alpha_1,\alpha_2\in\R$, and
$\alpha_\theta=(1-\theta)\alpha_1+\theta\alpha_2$.
The following real interpolation identities hold with equivalent
norms.

\smallskip
\noindent$\mathrm{(i)}$ If $1<p,q_1,q_2,u<\infty$, $W\in\mathcal A_p$, and
$\alpha_1\ne\alpha_2$, then
\begin{align}
 (\dot B^{\alpha_1}_{p,q_1}(W),
  \dot B^{\alpha_2}_{p,q_2}(W))_{\theta,u}
 &=\dot B^{\alpha_\theta}_{p,u}(W),
 \label{eq:B-real-smoothness}\\
 (\dot F^{\alpha_1}_{p,q_1}(W),
  \dot F^{\alpha_2}_{p,q_2}(W))_{\theta,u}
 &=\dot B^{\alpha_\theta}_{p,u}(W).
 \label{eq:F-real-smoothness}
\end{align}
\smallskip
\noindent$\mathrm{(ii)}$ If $1<p_1,p_2<\infty$, $p_1\ne p_2$,
$W_i\in\mathcal A_{p_i}$ for $i\in\{1,2\}$, and
$p_\theta,W_\theta$ are as in
\eqref{eq:intermediate-weight}, then
\begin{align}
 (\dot B^{\alpha_1}_{p_1,p_1}(W_1),
  \dot B^{\alpha_2}_{p_2,p_2}(W_2))_{\theta,p_\theta}
 &=\dot B^{\alpha_\theta}_{p_\theta,p_\theta}(W_\theta),
 \label{eq:B-real-diagonal}\\
 (\dot F^{\alpha_1}_{p_1,p_1}(W_1),
  \dot F^{\alpha_2}_{p_2,p_2}(W_2))_{\theta,p_\theta}
 &=\dot F^{\alpha_\theta}_{p_\theta,p_\theta}(W_\theta).
 \label{eq:F-real-diagonal}
\end{align}
\end{corollary}

\begin{proof}
Take $r=2$ in Lemma~\ref{lem:BT-retraction}.
Parts $\mathrm{(i)}$ and $\mathrm{(ii)}$ follow from
Theorems~\ref{lem:Z-real-smoothness} and~\ref{thm:Z-real},
respectively, using \cite[Theorem C.3.3]{HNVW16}.
\end{proof}

\medskip
\noindent\textbf{Acknowledgments.}
The authors acknowledge the use of AI tools during the exploratory stage
of this project. All mathematical arguments and proofs in the final
manuscript were checked and written by the authors.

\bigskip

\noindent\begin{minipage}{\textwidth}
\noindent Tengfei Bai, Pengfei Guo

\medskip

\noindent School of Mathematics and Statistics, Hainan Normal University,\newline
Haikou, Hainan 571158, China

\medskip

\noindent Jingshi Xu (Corresponding author)

\medskip

\noindent School of Mathematics and Computing Science, Guangxi Colleges and Universities Key Laboratory of Data Analysis and Computation, Guilin University of Electronic Technology, Guilin, 541004, China

\smallskip

\noindent Center for Applied Mathematics of Guangxi (GUET), Guilin, 541004, China

\smallskip

\noindent {\it E-mails}:
\texttt{202311070100007@hainnu.edu.cn}, \texttt{baitengfei1@126.com} (T. Bai)

\noindent\phantom{{\it E-mails:}}
\texttt{guopf999@163.com} (P. Guo)

\noindent\phantom{{\it E-mails:}}
\texttt{jingshixu@126.com} (J. Xu)
\end{minipage}
\end{document}